\documentclass[10pt,a4paper,oneside,psamsfonts]{report}
\usepackage{amsmath}
\usepackage{ngerman} 
\usepackage{amsfonts}
\usepackage{mathrsfs}
\usepackage{amssymb}
\usepackage{tabularx}
\usepackage{amssymb}
\usepackage{graphicx}
\usepackage[latin1]{inputenc}

\usepackage{geometry}
\newcommand{\C}{{\mathbb C}}
\newcommand{\N}{{\mathbb N}}

\newcommand{\R}{{\mathbb R}}
\newcommand{\Z}{{\mathbb Z}}
\newcommand{\pz}{{\mathbb P}}

\usepackage[linktocpage={false}]{hyperref}

\usepackage{amsthm}

\newenvironment{beweis}{\begin{proof}[Beweis]}{\end{proof}}

\begin{document}

\theoremstyle{plain}
\newtheorem{theorem}{Theorem}[chapter]
\newtheorem{lemma}[theorem]{Lemma}

\theoremstyle{remark}
\newtheorem*{bemerkung}{Bemerkung}
\theoremstyle{definition}
\newtheorem*{bedingung}{Bedingung}

\begin{center}
\vspace*{3cm}
 {\Huge Über die Linniksche Konstante }
\\
\vspace*{4cm}
{\Large Diplomarbeit}\linebreak \\
{\Large von}\linebreak \\
{\Large Triantafyllos Xylouris}\\
\vspace*{8cm}

{
\begin{tabular}{ll}
Betreuer und Erstprüfer: & PD Dr. B. Z. Moroz, Universität Bonn\\
Zweitprüfer: & Prof. Dr. S. Wewers, Universität Hannover\\
Abgabetermin: & 22. April 2009 \\
Datum dieser Version: & 12. Juni 2009 (weniger Druckfehler \\
& gegenüber der Version vom 22. April 2009)\\
\end{tabular}
}
\end{center}

\newpage

\tableofcontents

 \newpage

\chapter{Einführung}\label{sec:Test-chapter}
\section{Abstract}
Seien $a$ und $q$ teilerfremde natürliche Zahlen. 1944 bewies Yu. V. Linnik (\cite{Lin44-1}, \cite{Lin44-2}), dass die kleinste Primzahl $p \; (\equiv a \; mod \; q)$ in einer arithmetischen Progression kleiner als $C q^L$ ist. Dabei sind $C$ und $L$ positive effektiv berechenbare Konstanten. Seitdem wurde der zulässige Wert für die Konstante $L$ oft verbessert, zuletzt 1992 durch D. R. Heath-Brown auf $L=5.5$ (\cite{Hea92}). Im selben Artikel gibt Heath-Brown verschiedene kleine Potentiale zur Verbesserung seiner Arbeit an. Mit diesen Potentialen verbessern wir die Zwischenresultate aus \cite{Hea92} und schließlich die Konstante auf $L=5.2$.

\section{Notationen}\label{sec:Einfuehrung-Notationen}
\noindent Wir benutzen die folgenden Notationen, die in der analytischen Zahlentheorie üblich sind. Dazu seien $x \in \R,\,  m, n, a, q \in \N$:
\begin{itemize}
 \item $\N = \{1,2,3,\ldots \}$ und $\N_0 = \{0,1,2,3, \ldots\}$.
 \item Die Menge $\{a+jq\,|\,j \in \N_0 \}$ nennen wir arithmetische Progression $a \; mod \; q$.
 \item $\pz = \{2,3,5,7,11,\ldots\}$ (Primzahlen), $p$ bzw. $p_i$ steht immer für eine Primzahl,
 \item $\Lambda(n)=\log p$, wenn $n=p^m$ und $0$ sonst,
 \item $\pi (x) = \sum\limits_{p \leq x} 1 , \; \; \psi (x) = \sum\limits_{n \leq x} \Lambda(n), \; \; \pi(x;q,a)= \sum\limits_{p \leq x,\atop p \equiv a \; mod \; q } 1 $,
 \item $li(x)=\int_2^{x} \frac{dt}{\log t}$,
 \item $(m,n)=ggT(m,n)=\max\{d \in \N \,| \, d|m \text{ und } d|n\}$, es wird keine Gefahr der Verwechslung mit einem Intervall bestehen,
 \item $kgV(m,n)=\min \{d \in \N \,|\, m | d \text{ und } n | d \}$,
 \item $\varphi(q)$ ist die Eulersche $\varphi$-Funktion,
 \item $\mu(d)$ ist die Möbiussche $\mu$-Funktion,
 \item $x = [x] + \{x\}$ mit $[x]=\max\{n \in \Z \,|\, n \leq x\}$.
 \item Seien $X \subset \C$ und $g : X \rightarrow \C,\, f:X \rightarrow [0,\infty)$ Funktionen. Dann schreiben wir $g(x)=O(f(x))$ bzw. gleichbedeutend $g(x) \ll f(x)$  (auch $f(x) \gg g(x)$), falls $\exists \, C \in \R \; \forall \, x \in X : |g(x)|\leq C f(x)$. Dabei wird $X$ in der Regel aus dem Kontext ersichtlich sein. Die Konstante $C$ nennen wir auch implizite Konstante.  Meist werden die Funktionen $g$ und $f$ von verschiedenen Parametern abhängen. Dann hängt die implizite Konstante in der Regel auch von diesen Parametern ab. Manchmal werden wir dies durch die Schreibweise $O_A (f(x))$ bzw. $\ll_A f(x)$ kennzeichnen, was bedeutet, dass $C$ von $A$ abhängt.\\Weiterhin schreiben wir $g(x)=o(f(x))$, falls $\frac{g(x)}{f(x)} \rightarrow 0$ für $x \rightarrow \infty$.
 \item Wenn wir für ein $s \in \C$ schreiben $s=\sigma+it$, so gilt die stille Vereinbarung, dass $\sigma,t \in \R$. Wir schreiben $Re\{s\}=\sigma$ und $Im\{s\}=t$. Wenn wir sagen, dass $s$ komplex ist, so bedeutet dies, dass $s \notin \R$.
 \item Wir vereinbaren folgende Abkürzungen für $x,y>0$ : $\log^y x :=(\log x)^y$ und $\log xy :=\log (xy)$.
 \item $\chi$ steht für einen Dirichlet-Charakter ($mod \; q$) und $\chi_0$ für den Hauptcharakter ($mod \; q$). Unter einem Charakter $\chi \neq \chi_0$ verstehen wir einen Dirichlet-Charakter $\chi$, der ungleich dem Hauptcharakter ist. Unter einer Dirichletschen L-Funktion verstehen wir die Funktion $$L(s,\chi)=\sum_{n=1}^{\infty} \chi(n) n^{-s},$$ abhängig von der komplexen Veränderlichen $s$.\\
     Mit $ord \; \chi$ bezeichnen wir die Ordnung des Charakters $\chi$ in der Gruppe der Dirichlet-Charaktere $mod \; q$. Ist $ord \; \chi = 2$, so ist $\chi(n)$ reellwertig, wir sagen auch $\chi$ ist ein reeller Charakter. Ist $ord \; \chi \geq 3$, so nennen wir $\chi$ einen komplexen Charakter.
\end{itemize}

\noindent Zusätzlich gelten die folgenden Vereinbarungen:
\begin{itemize}
 \item Wir benutzen die Abkürzung $\mathscr{L}:=\log q$.
 \item $\varepsilon$ steht immer für eine reelle und positive Konstante. Diese ist keineswegs immer gleich und kann von Zeile zu Zeile variieren.
\item Wir beweisen die meisten Resultate nur unter der Bedingung $q \geq q_0$, wobei $q_0$ jeweils eine absolute hinreichend große Konstante ist. Das werden wir nicht immer explizit erwähnen. Weiterhin weisen wir darauf hin, dass jegliche implizite Konstanten in den Kapiteln 2 bis 4 dieser Arbeit effektiv berechenbar sind.
\end{itemize}

\noindent Schließlich geben wir hier eine Liste von ausgewählten Bezeichnungen an, die wir im Laufe der Arbeit definieren und benutzen werden:
\begin{itemize}
  \item $p(a,q), \; P(q)$ (S.\pageref{Definition-p-P})
    \item $K(s,\chi)$ (S.\pageref{Definition-K-s-chi})
  \item $\phi=\phi(\chi)$, Bedingung 1 und 2, Laplace-Transformierte $F$ von $f$ (S.\pageref{Definition-phi-chi}-\pageref{Definition-Bedingungen})
  \item $R_0(K),\; L_1 $ (S.\pageref{Definition-R0})
  \item $R, \; A_1, \; A_2$ (S.\pageref{hb-lemma-5-2-extended})
  \item $\chi_j, \; \rho_j,  \; \beta_j, \; \gamma_j, \; \lambda_j, \; \mu_j, \; \rho',\; \beta',\; \gamma',\; \lambda',\; \mu'$ (S.\pageref{Definition-rho-k})
  \item $\gamma,\; g(t), \; f(t),\; F(z)$ (S.\pageref{explizites-f}-\pageref{explizites-F})
  \item $\Sigma, \; h(x), \; H(z), \; \sum_{\rho}{'}, \; R_1, \; H_2(z)$ (S.\pageref{Definition-Sigma}-\pageref{Definition-H-2})
  \item $A(s_1,s_2,t),\; A_{sup}$ (S.\pageref{Definition-A})
  \item $K_0(s,\chi)$ (S.\pageref{Definition-K0})
  \item Tabellen 2-3 (S.\pageref{Tabelle2-3}), Tabellen 4-10 (S.\pageref{Tabelle-4}-\pageref{Tabelle-10}), Tabellen 12-13 (S.\pageref{Tabelle-12})
  \item $\chi^{(j)},\; \rho^{(j)},  \; \beta^{(j)}, \; \gamma^{(j)}, \; \lambda^{(j)}, \; \mu^{(j)}$, $N(\lambda)$ (S.\pageref{Definition-rho-j-oben})
\end{itemize}

\section{Einführung in die Problematik}\label{Kapitel1-Problematik}
Im Folgenden seien immer $a,q \in \N$ mit $(a,q)=1$. Dirichlet bewies 1837, dass die arithmetische Progression $\{a+jq\,|\, j \in \N_0 \}$ unendlich viele Primzahlen enthält. Dies wurde später durch den Primzahlsatz für arithmetische Progressionen (1896 durch de la Vallée-Poussin) präzisiert, nämlich
\begin{equation}\label{PNT-festes-q} \pi (x;q,a)=\frac{x}{ \varphi(q) \log x} +o_q\left(\frac{x}{\log x}\right).\end{equation}
Viel Interesse hat die Frage erzeugt, inwiefern (\ref{PNT-festes-q}) gleichmäßig in $q$ ist (siehe z.B. \cite[Kapitel 17]{IwaKow04}). Man kann zwar keine genau Formel im Stile von (\ref{PNT-festes-q}) für $x \approx q$ erwarten, da die betrachtete arithmetische Progression in $[1,q]$ nur ein Element hat. Eine nicht-triviale Formel dafür würde bedeuten, dass man angeben kann, ob die Zahl $a$ prim ist oder nicht, was natürlich so ohne weiteres nicht geht.\\
Ist jedoch $x$ groß genug in Abhängigkeit von $q$, so ist die Aussage des Siegel-Walfisz Theorems (1936) nicht-trivial: Für $x \geq 2$ und $A>3$ gilt (siehe z.B. \cite[S.124]{IwaKow04})
\begin{equation}\label{PNT-SW} \pi (x;q,a) = \frac{li(x)}{\varphi(q)} + O_A \left(\frac{x}{\log^A x}\right),\end{equation}
wobei die implizite Konstante von $A$ abhängt, aber nicht effektiv berechenbar ist. Man beachte, dass für $q \ll_A \log^{A-2} x$ der Hauptterm in letzterer Gleichung von größerer Ordnung ist als der Fehlerterm. Da nun der Hauptterm unabhängig von $a$ ist und wir $A$ beliebig groß wählen können, bedeutet dies, dass für $e^{(q^{\varepsilon})} \ll_{\varepsilon} x$ eine "`Gleichverteilung"' der Primzahlen innerhalb der $\varphi(q)$ verschiedenen arithmetischen Progressionen $a \; mod \; q$ vorliegt.\\
Aus der Verallgemeinerten Riemannschen Vermutung für die Dirichletschen L-Funktionen\footnote{Diese Vermutung besagt, dass für einen Charakter $\chi \; mod \; q$ die  nicht-trivialen Nullstellen von $L(s,\chi)$ auf der Geraden $Re\{s\}=\frac{1}{2}$ liegen.} folgt noch mehr, nämlich für $x \geq 2$ (\cite[S.426]{MonVau07})
\begin{equation}\label{PNT-RV} \pi (x;q,a) = \frac{li(x)}{\varphi(q)} + O (x^{\frac{1}{2}} \log x) .\end{equation}
Dies liefert "`Gleichverteilung"' der Primzahlen in den verschiedenen arithmetischen Progressionen bereits ab $q^{2+\varepsilon} \ll_{\varepsilon} x$.\\

Verwandt mit der Frage, ab wann in den obigen Gleichungen der Hauptterm von größerer Ordnung ist als der Fehlerterm, ist die Frage, wie groß die \emph{erste} Primzahl in einer arithmetischen Progression ist. Mit anderen Worten: Sei $q$ fest. Wie groß muss man $x$ wählen, damit für alle $a$ mit $(a,q)=1$ gilt, dass $\pi(x;q,a)>0$  ? Wir führen folgende Notationen ein:\label{Definition-p-P}
$$p(a,q):=\min \{p \in \pz \,|\, p \equiv a \; mod \; q\}$$ und
\begin{eqnarray*}
  P(q)&:=&\max \{p(a,q)\,|\, a \in \N, \; a \leq q,\; (a,q)=1\} \\
      & =&\min \{ p \in \pz \,|\,\pi(p;q,a)>0 \text{ für } a \in \N, \, a\leq q, (a,q)=1 \}.
\end{eqnarray*}
 Wir interessieren uns nun für die Größenordnung von $P(q)$. Eine untere Schranke bekommt man sofort aus der Ungleichung $$\pi(P(q)) \geq \# \{a \in \N \,|\, a \leq q, (a,q)=1 \}=\varphi(q).$$
Daraus folgt für $q \rightarrow \infty$, dass
\begin{equation}\label{Paq-lower}P(q) \geq (1+o(1))\varphi(q) \log q,\end{equation}
da man sonst einen Widerspruch zum Primzahlsatz bekäme\footnote{Wir erinnern daran, dass für $n \geq 2$ gilt $$\frac{n}{4\log n} \leq \varphi(n) \leq n-1 .$$}.
Was obere Schranken angeht, so liefert (\ref{PNT-SW}) sofort $\pi(x;q,a)>0$ für $e^{(q^{\varepsilon})} \ll_{\varepsilon} x$ unabhängig von $a$, also
$$P(q) \ll_{\varepsilon} e^{(q^{\varepsilon})} .$$
Setzt man die Verallgemeinerte Riemannsche Vermutung für $L(s,\chi)$ voraus, bzw. (\ref{PNT-RV}), so erhält man
$$P(q)\ll_{\varepsilon} q^{2+\varepsilon} .$$
Man vermutet sogar eine noch stärkere Version von (\ref{PNT-RV}), welche
$$P(q) \ll_{\varepsilon} q^{1+\varepsilon} $$
zur Folge hat (siehe z.B. \cite[S.419]{IwaKow04}, \cite[S.1]{Bal07}). S.S.Wagstaff (1979, \cite{Wag79}) vermutet letztendlich aufgrund von heuristischen Überlegungen, dass  $$P(q) \approx \varphi(q) \log^2 q.$$
\vspace*{0.1cm}

\noindent In diesem Zusammenhang ist nun das folgende Ergebnis von Yu. V. Linnik bemerkenswert (1944 - \cite{Lin44-1}, \cite{Lin44-2}), nämlich
\begin{equation}\label{LinniksTheorem}P(q) \leq C q^L\end{equation}
mit effektiv berechenbaren positiven Konstanten $C$ und $L$. Die Zahl $L$, die hier auftaucht, wird manchmal "`Linniksche Konstante"' genannt und ist natürlich nicht eindeutig, da letztere Ungleichung für verschiedene $L$ gelten kann\footnote{Man kann die "`Linniksche Konstante"' auch als das Infimum von denjenigen Konstanten $L$ definieren, die in (\ref{LinniksTheorem}) zulässig sind. Dann hat man Eindeutigkeit.}. Wir werden in dieser Arbeit das Resultat (\ref{LinniksTheorem}) auch als "`Theorem von Linnik"' bezeichnen.\\
\noindent Linnik selbst gab keinen konkreten zulässigen Wert für $L$ an. C. D. Pan holte dies 1957 nach. Der zulässige Wert für $L$ wurde seitdem mehrmals verbessert. Es folgt eine Liste von bisher bewiesenen Verbesserungen.\\

\begin{center}
\begin{tabular}[t]{|llll|}
\multicolumn{4}{c}{\textbf{Tabelle 1. Zulässige Werte für $L$}}\\\hline
$L$ & Jahr & Autor &  Referenz  \\ \hline
10000& 1957 & Pan &  \cite{Pan57}\\
5448& 1958 & Pan &   \cite{Pan58}\\
777& 1965 & Chen &   \cite{Che65}\\
630& 1971 & Jutila &   \cite[S.370]{Tur71}\\
550& 1970 & Jutila &   \cite{Jut70}\\
168& 1977 & Chen &   \cite{Che77}\\
80& 1977 & Jutila &   \cite{Jut77}\\
36& 1977 & Graham &   \cite{Gra77}\\
20& 1981 & Graham &   \cite{Gra81}\\
17& 1979 & Chen &   \cite{Che79}\\
16& 1986 & Wang &   \cite{Wan86}\\
13.5& 1989 & Chen und Liu& \cite{CheLiu89}\\
5.5& 1992 & Heath-Brown & \cite{Hea92}\\
\hline
\end{tabular}
\end{center}
\vspace*{0.5cm}

Wie wird nun das Theorem von Linnik bewiesen und wie kann man den zulässigen Wert für $L$ verbessern? Die Behandlungen dieses Theorems basieren in der Regel auf den folgenden drei Prinzipien. Das erste Prinzip ist dabei allgemein bekannt (\cite[Satz 3.1.4]{Bru95}). Prinzip 2 und 3 gehen zurück auf Linnik (\cite{Lin44-1}, \cite{Lin44-2}):\\

\textbf{Prinzip 1 ("`Nullstellenfreie Region"')} \emph{Es gibt eine effektiv berechenbare positive Konstante $c_1$, so dass
\begin{equation}\label{prod-chi} \prod_{\chi \; mod \; q} L(s,\chi)\end{equation}
höchstens eine Nullstelle in der Region
$$\sigma \geq 1-\frac{c_1}{\log q (2+|t|)} \left(:= \{\sigma+it \in \C \,| \, \sigma \geq 1-\frac{c_1}{\log (q (2+|t|))} \}\right)$$
hat. Wenn diese Ausnahmenullstelle (auch genannt: "`Siegel-Nullstelle"' oder "`Siegel-Landau-Nullstelle"') existiert, so ist sie reell, einfach und gehört zu einem reellen Charakter $\chi \neq \chi_0$.}\\

\textbf{Prinzip 2 ("`Deuring-Heilbronn Phänomen"')} \emph{Es gibt eine effektiv berechenbare positive Konstante $c_2$, so dass, wenn die Ausnahmenullstelle in Prinzip 1 existiert und gleich $1-\lambda_1(\log q)^{-1}$ ist, dann hat die Funktion (\ref{prod-chi}) keine weiteren Nullstellen in der Region
$$\sigma \geq 1-\frac{c_2 \log (\lambda_1^{-1})}{\log q (2+|t|)}.$$ }

\noindent Letzteres Phänomen wird auch als "`Nullstellen-Abstoßung"' bezeichnet, weil die Existenz einer Ausnahmenullstelle nah bei $s=1$ erzwingt, dass alle weiteren Nullstellen "`viel weiter links"' von $s=1$ liegen müssen. Letztere werden sozusagen "`abgestoßen"'. \\

\textbf{Prinzip 3 ("`Log-Freie Abschätzung der Nullstellendichte"')} \emph{Es gibt zwei effektiv berechenbare positive Konstanten $c_3$ und $c_4$, so dass für $T \geq 1$ und $\sigma \in [\frac{1}{2},1]$ gilt:
$$\sum_{\chi \; mod \; q} N(\sigma,T,\chi) \leq c_3 (qT)^{c_4(1-\sigma)} .$$
Dabei ist $N(\sigma,T,\chi)=\# \{\rho \in \C  \,| \, L(\rho,\chi)=0,\, Re \{\rho\} \geq \sigma, \,|Im \{\rho\}| \leq T \}$.}\\

\noindent Der Wert der Konstanten $L$, für die man das Theorem von Linnik beweisen kann, hängt direkt von den Konstanten $c_1$, $c_2$, $c_3$ und $c_4$ ab. Dabei ist es durchaus möglich und später auch der Fall, dass zu obigen Prinzipien verwandte Theoreme zu Verbesserungen für die zulässige Konstante $L$ führen. Ein solches verwandtes Theorem wäre z.B. der Nachweis einer bis auf zwei Nullstellen nullstellenfreien Region.\\

\noindent Unter gewissen Voraussetzungen kann man den zulässigen Wert für $L$ verbessern. Hat $q$ beschränkten \label{Definition-Kubikanteil} Kubikanteil\footnote{Sei $q=u v^3$ mit $u,v \in \N$. Weiterhin gelte $p^3 \nmid u$ für alle $p \in \pz$. Dann nennen wir $v^3$ (eindeutig) den Kubikanteil von $q$.}, so ist $$L=4.5$$ zulässig (Z. Meng, \cite{Men01}). Seien $c_1$ und $\lambda_1$ so definiert wie in Prinzip 1 und 2. Existiert die Ausnahmenullstelle und ist  $\lambda_1 \leq \lambda(\varepsilon)$ für ein hinreichend kleines von $q$ unabhängiges $\lambda(\varepsilon)>0$, so ist $$L=3+\varepsilon$$ zulässig (Heath-Brown, \cite{Hea90}). Gilt Prinzip 1 für beliebig großes $c_1$ (und ab $q \geq q_0(c_1)$) und existiert die Ausnahmenullstelle nicht, so bekommt man (siehe Bemerkung in \cite[S.269]{Hea92}) $$L=2.4+\varepsilon.$$ Den gleichen Wert für $L$ erhält H. Iwaniec (\cite{Iwa74}) auch für jene $q$, deren Primteiler aus einer festen endlichen Menge stammen. Schließlich erinnern wir noch einmal an die weiter oben erwähnte Aussage, dass unter der Voraussetzung der Verallgemeinerten Riemannschen Vermutung für $L(s,\chi)$ man
$$L=2+\varepsilon$$
bekommt.\\

Nach diesen Vorbemerkungen möchten wir jetzt Ziel und Inhalt der vorliegenden Diplomarbeit benennen. Die letzte Verbesserung der zulässigen Linnikschen Konstante auf $L=5.5$ wurde durch Heath-Brown im Jahr 1992 bewiesen. Im besagten Artikel werden am Ende (\cite[S.332-337]{Hea92}) neun Verbesserungspotentiale angegeben, mit denen die Argumente des Artikels leicht verbessert bzw. verfeinert werden können. Wir gehen in dieser Diplomarbeit die neun Potentiale durch und verbessern damit die Zwischenergebnisse in \cite{Hea92}, sowie den zulässigen Wert für die Konstante $L$.\\

An dieser Stelle möchten wir noch erwähnen, dass die Konstanten, die Heath-Brown in \cite{Hea92} für die verschiedenen nullstellenfreien Regionen bekommt, im Allgemeinen deutlich besser sind, als die analogen Resultate anderer Autoren. Eine Ausnahme bildet der (recht spezielle) Fall, dass man Prinzip 1 nur für reelle Charaktere $\chi \neq \chi_0$ und reelle Nullstellen betrachtet. Für diesen Fall beweist J. Pintz mit einer anderen Vorgehensweise die Konstante $c_1=4+o(1)$ (\cite{Pin77}), was schon recht besser ist als die Abschätzung mit der Konstanten $2.427\ldots \;$ aus \cite[Lemma 8.2]{Hea92}. Heath-Brown weist darauf hin (\cite[S.269 oben]{Hea92}), dass dies vermuten lässt, dass es weiterhin Raum für Verbesserungen gibt für die aus seiner Methode erhaltenen nullstellenfreien Regionen.\\
In diesem Zusammenhang weisen wir daraufhin, dass wir in dieser Arbeit zwar viele Zwischenergebnisse aus \cite{Hea92} verbessern konnten. Dies trifft aber nicht auf den eben genannten Wert $c_1=2.427\ldots \;$ zu.\\

\section{Zusammenfassung der Ergebnisse}
Wir verbessern die Konstante von $L=5.5$ auf $L=5.2$:

\begin{theorem}\label{haupttheorem}
  Es gilt $$P(q) \ll q^{5.2}$$ mit einer effektiv berechenbaren impliziten Konstanten.
\end{theorem}

\noindent Die bestmögliche Konstante, die wir erreichen könnten, ist dabei nicht die "`glatte"' Zahl $L = 5.2$, sondern eine "`nichtglatte"' Zahl leicht unterhalb von $5.2$. Außerdem weisen wir darauf hin, dass man die Anzahl der Computerberechnungen in dieser Arbeit hätte exzessiv erhöhen können. Aufgrund gewisser Untersuchungen vermuten wir, dass man dann die zulässige Konstante bis auf $L=5.14$ hätte verbessern können. Dies haben wir jedoch nicht bewiesen.

Weiterhin sind in dieser Arbeit mehrere kleine Zwischenergebnisse enthalten, die die entsprechenden Zwischenergebnisse aus \cite{Hea92} geringfügig verbessern. An dieser Stelle möchten wir nur das Folgende hervorheben, bei dem wir eine Verbesserung der auftretenden Konstanten von $0.348$ (\cite[Theorem 1, S.268]{Hea92}) auf $0.440$ bekommen haben:\\

\begin{theorem}\label{theorem-1-von-HB} Es gibt eine effektiv berechenbare Konstante $q_0$ in $\R$, so dass für $q \geq q_0$ die Funktion (\ref{prod-chi}) höchstens eine Nullstelle in der Region $$\sigma \geq 1-\frac{0.440}{\log q}, \; |t| \leq 1 $$
hat. Wenn diese Ausnahmenullstelle existiert, so ist sie reell, einfach und gehört zu einem reellen Charakter $\chi \neq \chi_0$.
\end{theorem}

\begin{bemerkung}
  Durch eine kleine Variation der Schlussweise von Heath-Brown beweisen Liu und Wang (1998, \cite[S.345-346]{LiuWan98}) Theorem \ref{theorem-1-von-HB}  mit der Konstanten $0.364$.
\end{bemerkung}

\section{Kommentare zur Diplomarbeit}

\subsubsection{Allgemeines}

Diese Diplomarbeit basiert vollständig auf dem Artikel \cite{Hea92} von Heath-Brown. Jedoch haben wir versucht die Arbeit so zu erstellen, dass sie möglichst ohne die Hinzunahme des hiesigen Artikels nachvollziehbar ist. Zu diesem Zweck stellen wir in Kapitel 2 einige wesentliche Lemmata und Methoden aus \cite{Hea92} zusammen. In Kapitel 3 diskutieren wir einen Teil der Verbesserungspotentiale und verbessern oder erweitern damit einige Zwischenresultate aus \cite{Hea92}. In Kapitel 4 benutzen wir die Resultate aus Kapitel 3, um Theorem \ref{haupttheorem} zu beweisen. Jene Verbesserungspotentiale, aus denen wir keine "`nennenswerten Verbesserungen"' für die Endkonstante $L$ bekommen, diskutieren wir abschließend kurz in Kapitel 5.

An gewissen Stellen dieser Arbeit mag sich der Leser fragen, warum wir mal sehr genau und mal recht grob vorgegangen sind. Da die gesamte Arbeit sehr rechenlastig ist, haben wir in der Regel versucht nur so viele Rechnungen und Fallunterscheidungen durchzuführen, als für unser optimales Endergebnis (also $L=5.2$) notwendig ist. Die Entscheidung an welchen Stellen wir genauer und an welchen ungenauer vorgehen würden ergab sich in der Regel nicht durch Zufall, sondern resultierte aus bereits gemachten genaueren Rechnungen.

An dieser Stelle möchten wir noch einmal auf die in der Notation vereinbarte Abkürzung
$$\mathscr{L} =\log q $$
hinweisen. Weiterhin beweisen wir jegliche Resultate meist nur für $q \geq q_0$, wobei $q_0$ immer eine effektiv berechenbare hinreichend große Konstante ist.\\

\subsubsection{Zum Beitrag der Verbesserungspotentiale für das Endresultat "`$L=5.2$"'}

Von den neun Verbesserungspotentialen in \cite[S.332-337]{Hea92}, nennen wir sie VB1 bis VB9, lieferten VB1, VB3 und VB6 Verbesserungen, die "`nicht nennenswert"' waren\footnote{Unter "`nicht nennenswerten"' Verbesserungen verstehen wir jene, die zu weniger als ca. $0.005$ Verbesserung für die zulässige Endkonstante $L$ führten. Selbstverständlich müssen wir dabei betonen, dass diese Potentiale, so wie \emph{wir} sie verwendet haben, keine hilfreichen Verbesserungen lieferten.}.  VB8 lieferte zwar eine Verbesserung von "`mehreren Hundertstel"' (wir beziehen uns immer auf die Konstante $L$), aber dies nur für einen Teil der auftretenden Fälle, für den wir keine weitere Verbesserung nötig hatten. Ähnliches gilt für VB4. Deswegen haben wir VB1, VB3, VB4, VB6 und VB8 am Ende nicht benutzt.

Die Anwendung von VB2 in §\ref{para-ldash-abschaetzungen}, sowie in Verbindung mit VB5 in §\ref{para-l2-abschaetzungen} und §\ref{para9-l3bounds},  liefert "`ca. $0.15$-$0.20$ Verbesserung"'. Weiterhin benutzen wir VB2 in §\ref{para10analysis} für eine etwas erweiterte Analyse von \cite[§10]{Hea92}, was "`ca. $0.05$ Verbesserung"' liefert. VB7 (§\ref{para11-neu}) liefert eine Verbesserung von ca. $0.02$ und VB9 (§\ref{para12-neu}) nochmal um ca. $0.05$.\\
\noindent Man beachte, dass sämtliche genannten Verbesserungswerte nur als ungefähre Größenordnungen zu verstehen sind. Eine genaue Zuordnung, welches Potential wieviel zum Endergebnis $L=5.2$ beiträgt, ist aufgrund der Zusammenhänge untereinander, sowie aufgrund der vielen auftretenden Fälle, ohne weiteres nicht möglich.\\

Außer den neun genannten Potentialen, haben wir zahlreiche weitere (zum Teil recht naive) Ansätze ausprobiert, um die Zwischenergebnisse und damit die Endkonstante $L$ zu verbessern. Interessanterweise schlugen alle fehl (bis auf einen Ansatz, der leichte Verbesserungen versprach; jedoch hätte man dann nicht VB2 anwenden können und damit lohnte sich dieser auch nicht). Es blieb also nur das, was Heath-Brown schon vorgeschlagen hatte.\\

\subsubsection{Zu den Berechnungen}

Weiterhin müssen wir betonen, dass diese Arbeit - genauso wie der Artikel \cite{Hea92} (vergleiche Bemerkung in letzterem Artikel auf Seite 269 unten) - sich stark auf Computerberechnungen stützt. Dies ist notwendig, da wir ja die best-mögliche Konstante erzielen möchten. Deswegen darf man die Terme in den auftretenden Ungleichungen nicht zu sehr vereinfachen, sondern rechnet mit den komplizierten Termen weiter.

Die Berechnungen wurden mit dem Computer-Algebra-Programm Maple 12.0 auf einem Laptop mit einem Intel Core Duo 2 GHz-Prozessor und 2 GB Ram durchgeführt. Auf Anfrage stellen wir gerne jegliche benutzten Prozeduren und durchgeführten Rechnungen zur Verfügung.

Im Laufe der Arbeit war es nötig (und zwar immer dann, wenn die Methode aus §\ref{para-rechenkapitel} angewendet wurde) eine große Anzahl von gewissen Integralen zu berechnen. Einerseits berechneten wir eine explizite Form dieser Integrale (was angesichts einer existierenden Stammfunktion möglich war) und ließen Maple damit die Werte der Integrale berechnen. Andererseits verwendeten wir zur Kontrolle die normale zugehörige Maple-Prozedur, welche numerisch integriert. Die Werte, die aus den beiden verschiedenen Methoden resultierten, stimmten dabei bis zur ca. neunten Stelle\footnote{Die Nummerierung der Stellen beginnen wir dabei ab der ersten Stelle, die ungleich Null ist.} überein. Diese Fehlertoleranz ist nicht überraschend, da wir im Falle der numerischen Integrationen mit 10 Stellen Rechengenauigkeit rechneten. Allgemein ließen wir Maple mit mindestens 10 Stellen Rechengenauigkeit rechnen, was eine höhere Genauigkeit darstellt, als wir für unsere Resultate benutzten.

\section{Danksagung}
Ich möchte Herrn PD Dr. B. Z. Moroz für die Betreuung und konstante Unterstützung während der gesamten Bearbeitungszeit dieser Diplomarbeit eingehend danken.

\chapter{Die Arbeit von Heath-Brown}
Wie in der Einführung erwähnt, basiert diese Diplomarbeit vollständig auf dem Artikel \cite{Hea92}.
Deswegen möchten wir in diesem Kapitel gewisse Lemmata und Methoden aus \cite{Hea92} vorstellen, die wir später benötigen. Bei zwei Lemmata werden wir die kurzen Beweise miteinfügen, weil die entsprechenden Sachverhalte in \cite{Hea92} eher beiläufig erklärt werden und nicht einwandfrei zitierbar sind.\\
Schließlich werden wir auch den Zusammenhang zwischen den Nullstellen der Dirichletschen L-Funktionen und dem Theorem von Linnik skizzieren.

\section{Nullstellenfreie Regionen}
\subsection{Verbesserung der Standardmethode durch \cite{Hea92}}\label{para-methoden-fuer-nullstellenfreie-regionen}

Wir möchten hier erst die übliche Methode ("`Standardmethode"') vorstellen, mit der man Prinzip 1 aus §\ref{sec:Test-chapter} beweist, um dann aufzuzeigen wie diese in \cite{Hea92} verbessert wird. Wir betrachten also für einen Charakter $\chi$ die Funktion $L(s,\chi)$. Der Einfachheit halber setzen wir voraus, dass $\chi^2 \neq \chi_0$ und dass $q$ groß genug ist. Außerdem betrachten wir nur Nullstellen $\rho$ mit $|Im\{\rho\}|\leq \mathscr{L}\,(=\log q)$. Man beachte, dass der Beweis für die verbleibenden Fälle, als auch der Beweis für Prinzip 2, völlig analog abläuft.\\

\noindent Nun also zur Standardmethode. Man betrachtet die reellwertige Funktion
$$ w(\sigma+it):= Re  \left\{-\frac{L'}{L} (\sigma + it,\chi) \right\}=\sum_{n=1}^{\infty} \frac{\Lambda(n)}{n^{\sigma}} Re \left\{\frac{\chi(n)}{n^{it}}\right\} \;\;\; (\sigma>1, \, t\in \R).$$
Aus der Ungleichung $$0 \leq 2(1+\cos \theta)^2=3+4\cos \theta + \cos 2\theta$$ folgert man dann (setze $\theta=\arg(\chi(n))-t \log n$, wobei $e^{i \arg(\chi(n))}=\chi(n)$, multipliziere mit $\Lambda(n) n^{-\sigma}$ und summiere über $n$)
\begin{equation}\label{zeta-ungleichung} 0 \leq -3 \frac{\zeta'}{\zeta}(\sigma) -4Re  \left\{\frac{L'}{L}(\sigma + it,\chi)\right\}-Re  \left\{\frac{L'}{L}(\sigma + 2it,\chi^2)\right\} \;\;\; ( \sigma >1, \,t \in \R) .\end{equation}
\newpage
\noindent Die Funktion $(L'/L)(s,\chi)$ ist meromorph und hat bei den Nullstellen $\rho$ von $L(s,\chi)$ Pole erster Ordnung. Also gilt\footnote{Wähle dazu eine Nullstelle $\rho$ von $L(s,\chi)$, setze $s-\rho=r e^{i \theta} \; (r \geq 0,\, \theta \in \R)$ und $k$ gleich der Nullstellenordnung von $\rho$. Aus der Poleigenschaft von $\rho$ folgt, dass es eine in $\rho$ holomorphe Funktion $w_1(s)$ gibt, so dass in einer Umgebung von $\rho$ gilt
$$ w(s)= Re  \{w_1(s)\}+ Re \frac{-k}{s-\rho} = O(1) - k \frac{ \cos \theta}{r}. $$Letzteres geht für $r \rightarrow 0$ und $\cos \theta \geq 0$ gegen $-\infty$. Diese beiden Bedingungen sind aber genau die genannten Bedingungen $s \rightarrow \rho$ und $Re  \{s-\rho\} \geq 0$.} für eine Nullstelle $\rho$
\begin{equation*}\lim_{s \rightarrow \rho \atop Re  \{s-\rho\} \geq 0} w(s) =- \infty .\end{equation*}
Die Idee ist nun grob gesagt die Folgende. Dazu lassen wir die Voraussetzung $\sigma>1$ in (\ref{zeta-ungleichung}) mal kurz außer Acht: Sei $\rho_0=\sigma_0 + it_0$ eine nicht-triviale Nullstelle von $L(s,\chi)$. Setzen wir $t=t_0$ in (\ref{zeta-ungleichung}) und lassen wir dort das $\sigma$ von rechts gegen $\sigma_0$ gehen, so wird der zweite Term in (\ref{zeta-ungleichung}) wegen der Poleigenschaft gegen "`$-4 \infty$"' gehen, der erste Term geht gegen etwas, dass "`$\leq+3 \infty$"' ist und der dritte Term ist nach oben beschränkt. Insgesamt bekämen wir "`$0\leq - \infty + C$"', also einen Widerspruch, wenn $\sigma$ hinreichend nah bei $\sigma_0$ ist. Das Ergebnis: Die Nullstelle $\rho_0$ darf nicht "`zu nah"' an irgendeinem Wert $\sigma+it$ liegen, der in die Ungleichung (\ref{zeta-ungleichung}) eingesetzt werden darf.\\

\noindent Wir möchten die eben genannte Idee präzisieren. Dazu brauchen wir konkrete Abschätzungen der drei Terme auf der rechten Seite von (\ref{zeta-ungleichung}). Diese erhält man im Wesentlichen aus der Partialbruchzerlegung von $(L'/L) (s,\chi)$ (vergleiche z.B. \cite[S.273]{Hea92}, \cite[Solution 6.5.9]{Mur08}). Für $1< \sigma < 2, \; |t|\leq \mathscr{L}, \; \varepsilon>0$ und $q \geq q_0(\varepsilon)$ gilt
\begin{equation}\label{zeta-pbz-1}- \frac{\zeta'}{\zeta} (\sigma) < \frac{1}{\sigma-1} + \varepsilon \mathscr{L}, \end{equation}
\begin{equation}\label{allgemein-Ldash-L-abschaetzung}- Re  \left\{ \frac{L'}{L} (\sigma + it,\chi) \right\}<-\sum_{\rho} Re \frac{1}{s-\rho} + (\frac{1}{2}+\varepsilon) \mathscr{L}.\end{equation}
Wegen $Re\{1/(s-\rho)\}\geq 0$ folgt aus der letzten Ungleichung
\begin{equation}\label{zeta-pbz-2}- Re  \left\{ \frac{L'}{L} (\sigma + it_0,\chi) \right\}<-\frac{1}{\sigma-\sigma_0} + (\frac{1}{2}+\varepsilon) \mathscr{L}, \end{equation}
\begin{equation}\label{zeta-pbz-3}- Re  \left\{ \frac{L'}{L} (\sigma + 2it_0,\chi^2) \right\}<(\frac{1}{2}+\varepsilon) \mathscr{L}. \end{equation}
Aus (\ref{zeta-ungleichung})-(\ref{zeta-pbz-3}) folgt dann für $1<\sigma<2$, dass
$$0 \leq \frac{3}{\sigma-1}- \frac{4}{\sigma-\sigma_0} +(\frac{5}{2}+8\varepsilon)\mathscr{L} $$
und daraus für eine Konstante $C>0$, dass $\sigma_0< 1- C \mathscr{L}^{-1}$.\\

Welche Veränderungen nimmt nun Heath-Brown an diesem Standardargument vor, um die nullstellenfreie Region für $L(s,\chi)$ zu verbessern?
\begin{itemize}
\item Offensichtlich spielt die Konstante $\frac{1}{2}$ in (\ref{allgemein-Ldash-L-abschaetzung}) eine wichtige Rolle für die Weite der resultierenden nullstellenfreien Region bzw. für den Wert der Konstanten $C$. Indem Heath-Brown die meisten Nullstellen in der Summe in (\ref{allgemein-Ldash-L-abschaetzung}) weglässt, ist er in der Lage die Konstante auf $\frac{1}{6}$
    (bzw. in manchen Spezialfällen $\frac{1}{8}$) zu verbessern (\cite[Lemma 3.1]{Hea92}).
\item Eine erste spezielle Methode, die in \cite{Hea92} benutzt wird, besteht darin zusätzlich zur Funktion $(L'/L)(s,\chi)$ auch ihre $k$-te Ableitung $$\left(\frac{L'}{L}\right)^{(k)} (s,\chi) $$ zu verwenden (\cite[§4]{Hea92}). Man erhält dann verschiedene neue Ungleichungen von ähnlichem Typ wie (\ref{zeta-ungleichung})-(\ref{allgemein-Ldash-L-abschaetzung}). Durch geschicktes Kombinieren solcher Ungleichungen erhält man letztendlich Verbesserungen der nullstellenfreien Region. Diese Methode wird in unserer Arbeit jedoch gar nicht vorkommen. Der Grund ist, dass sich für uns die nachfolgend genannte zweite Methode in Verbindung mit Verbesserungspotential Nr. 2 als überlegen erwies.
\item Eine zweite, sehr hilfreiche Methode aus \cite{Hea92} liegt darin, anstelle von
$$-Re \left\{\frac{L'}{L}(s,\chi) \right\}=\sum_{n=1}^{\infty} \Lambda(n) Re \left\{ \frac{\chi(n)}{n^s}\right\}$$ die Variante \label{Definition-K-s-chi}
$$K(s,\chi):=\sum_{n=1}^{\infty} \Lambda(n) Re \left\{\frac{\chi(n)}{n^s}\right\}f(\mathscr{L}^{-1} \log n) $$ zu betrachten. Dabei wird $f$ aus einer gewissen Menge von Funktionen gewählt, die ab einer Stelle $x_0>0$ gleich Null sind und die eine zu $Re \{1/(s-\rho)\}\geq 0$ analoge Bedingung erfüllen.\\
Ein Vorteil dieser Methode ist, dass man $f$ aus einer größeren Menge von Funktionen wählen kann, insbesondere so wählen kann, dass die resultierenden Konstanten bestmöglich werden. Außerdem wird man, da $f$ größtenteils gleich Null ist, in den zu (\ref{zeta-ungleichung}) analogen Ungleichungen auch Werte für $s$ einsetzen können, die links von $Re\{s\}=1$ liegen, was wiederum vorteilhaft ist.\\
Im folgenden Abschnitt möchten wir die beiden für $K(s,\chi)$ resultierenden Analoga zu (\ref{zeta-pbz-1}) und (\ref{allgemein-Ldash-L-abschaetzung}) zitieren.
\end{itemize}

\subsection{Zwei wichtige Lemmata (\cite[§5]{Hea92})}

Sei $\chi$ ein Charakter $mod \; q$ mit $\chi \neq \chi_0$. Heath-Brown ordnet dem Charakter $\chi$ den Wert $\phi(\chi)$ zu durch (vergleiche \cite[Lemma 2.5]{Hea92}) \label{Definition-phi-chi}
$$\phi = \phi(\chi)=\left\lbrace \begin{array}{ll}
\frac{1}{4} & \text{ falls} \; \text{ $q$ kubik-frei\footnotemark[2] oder } ord \; \chi \leq \mathscr{L}, \\
\\
\frac{1}{3} & \text{ sonst.}\\
\end{array} \right. $$
\footnotetext[2]{Das bedeutet, dass für alle Primzahlen $p$ gilt $p^3 \nmid q$.}
 Dieses $\phi$ stammt aus einer Abschätzung für die Funktion $L(s,\chi)$, welche eben für gewisse Spezialfälle verbessert werden kann.\\
Es sei erwähnt, dass nahezu sämtliche Zwischenergebnisse und damit auch die letztendlich zulässige Linniksche Konstante $L$ von der Größe des $\phi$ abhängen. Je kleiner das $\phi$, desto besser die Ergebnisse. In diese Richtung zeigt Meng in \cite{Men01} durch Durchgehen des Artikels \cite{Hea92} mit $\phi=\frac{1}{4}$, dass für $q$ mit beschränktem Kubikanteil (Definition von Kubikanteil in der Fußnote auf S.\pageref{Definition-Kubikanteil}) gilt: $P(q) \ll q^{4.5}$. In dieser Arbeit werden wir jedoch wie in \cite{Hea92} meist $\phi=\frac{1}{3}$ nehmen müssen.\\
Wir definieren noch die beiden folgenden Bedingungen an eine Funktion $f$.\\

\noindent \textbf{Bedingung 1} (\cite[S.280]{Hea92})  \emph{Sei $x_0>0$ und sei $f:[0,\infty) \rightarrow \R$ eine stetige Funktion mit $f(x)=0$ für $x \geq x_0$. Weiterhin sei $B>0$ eine Konstante, so dass für alle $x \in (0,x_0)$ die Funktion $f$ zwei Mal stetig differenzierbar ist und $|f''(x)|\leq B$.}\\

\noindent
Wenn nun $f$ Bedingung 1 erfüllt, so ist $f(t)$ als stetige Funktion, die für $t \geq x_0$ gleich Null ist, auch beschränkt. Aus der Beschränktheit von $f''(t)$ auf $(0,x_0)$ folgt außerdem, dass auch $f'(t)$ auf $(0,x_0)$ beschränkt ist (siehe \cite[(5.1)]{Hea92}) und dass die beiden Grenzwerte $\lim\limits_{x \rightarrow 0^+} f'(x)$ und $\lim\limits_{x \rightarrow x_0^-} f'(x)$ in $\R$ existieren.\\

\noindent \textbf{Bedingung 2} (\cite[S.286]{Hea92})  \emph{Es sei $f(t)\geq 0$ für $t \in [0,\infty)$. Die Laplace-Transformierte
 $$F(z):=\int_{0}^{\infty} e^{-zt}f(t) \,dt $$
von $f$ erfülle für $Re \{z\}\geq 0$ die Ungleichung $Re  \{F(z)\}\geq 0$.}\\ \label{Definition-Bedingungen}

\noindent Allgemein gilt für die gesamte Arbeit, dass mit $F(z)$ immer die Laplace-Transformierte der Funktion $f(t)$ gemeint ist, wobei letztere Funktion aus dem Zusammenhang ersichtlich sein wird.\\
In §\ref{hb-para-7} werden wir konkrete Beispiele für Funktionen $f$ vorstellen, die diese beiden Bedingungen erfüllen. Wir kommen jetzt zu den beiden Lemmata.

\begin{lemma}[Lemma 5.3 aus \cite{Hea92}]\label{hb-lemma-5-3}
Sei $s=\sigma+it \in \C$ mit
$$ |\sigma-1|\leq \frac{(\log \mathscr{L})^{\frac{1}{2}}}{\mathscr{L}},\; |t|\leq \mathscr{L} .$$
Weiterhin sei $f$ eine Funktion, die Bedingung 1 erfüllt. Dann gilt
$$ \sum_{n=1}^{\infty} \Lambda(n) \frac{\chi_0(n)}{n^{s}} f(\mathscr{L}^{-1} \log n)
= \mathscr{L} (F((s-1)\mathscr{L}))+O\left(\frac{\mathscr{L}}{\log \mathscr{L}}\right),$$
wobei die implizite Konstante beim "`Groß O"' nur von $f$ abhängt.
\end{lemma}
\vspace*{0.3cm}
\begin{lemma}[Lemma 5.2 aus \cite{Hea92}]\label{hb-lemma-5-2}
Sei $\chi \neq \chi_0$ ein Charakter $mod \; q$ und $\phi=\phi(\chi)$ so definiert wie zu Beginn dieses Unterabschnitts. Ferner sei $s=\sigma+it \in \C$ mit
\begin{equation} |\sigma-1|\leq \frac{(\log \mathscr{L})^{\frac{1}{2}}}{\mathscr{L}}, \; |t|\leq \mathscr{L} .\end{equation}
Weiterhin sei $f$ eine Funktion, die Bedingung 1 erfüllt und für die $f(0) \geq 0$ gilt. Dann gibt es für jedes $\varepsilon>0$ ein $\delta \in (0,1)$, das von $f$, aber nicht von $\chi,\, q$ oder $s$ abhängt, und ein $q_0=q_0(f,\varepsilon)$ in $\N$, so dass für alle $q \geq q_0$ gilt
$$\sum_{n=1}^{\infty} \Lambda(n) Re  \left\{\frac{\chi(n)}{n^{s}}\right\}f(\mathscr{L}^{-1} \log n)
\leq -\mathscr{L} \sum_{|1+it-\rho|\leq \delta} Re  \{F((s-\rho)\mathscr{L})\} + \frac{f(0)}{2} \phi \mathscr{L} + \varepsilon\mathscr{L} .$$
Dabei erstreckt sich die Summation über jene nicht-trivialen Nullstellen $\rho$ von $L(s,\chi)$, welche die Bedingung $|1+it-\rho|\leq \delta$ erfüllen. Jede Nullstelle taucht gemäß Ihrer Vielfachheit auf.
\end{lemma}

\noindent In \cite{Hea92} wird in der Aussage von Lemma \ref{hb-lemma-5-2} anstelle von $\frac{1}{2} f(0) \phi \mathscr{L} + \varepsilon\mathscr{L}$ der Term $f(0)(\frac{1}{2} \phi + \varepsilon)\mathscr{L}$ angegeben. Betrachtet man den Beweis des Lemmas, so handelt es sich offensichtlich um einen Druckfehler. Man beachte, dass diese kleine Änderung jedoch praktisch keinen Unterschied macht: Ist $f(0)>0$, so sind beide Terme "`äquivalent"', was durch eine Anpassung des $\varepsilon$ folgt. Einen Unterschied gibt es nur für $f(0)=0$, dann wäre $f(0)(\frac{1}{2}\phi+\varepsilon)\mathscr{L}=0$, was nicht zulässig ist. Die benutzten Funktionen $f$ werden später immer die Bedingung $f(0)>0$ erfüllen.\\
Man beachte, dass das Lemma trivialerweise auch dann gilt, wenn man $\phi$ durch eine größere Zahl ersetzt. Wenn wir später also nicht wissen, ob $\phi=\frac{1}{3}$ oder $\phi=\frac{1}{4}$ ist, so werden wir den größeren Wert $\frac{1}{3}$ nehmen.

\subsection{Die Anwendung der beiden Lemmata (\cite[§6]{Hea92})}\label{hb-para-6}
Für den Beweis des Theorems von Linnik werden am Ende nur jene Nullstellen $\rho$ relevant sein, die $Re \{\rho\} \geq 1- C \mathscr{L}^{-1}$ und $|Im\{\rho\}| \leq C \mathscr{L}^{-1}$ erfüllen, wobei $C$ eine absolute positive Konstante ist. Deswegen reicht es aus nur Nullstellen im folgenden Viereck zu betrachten, welches von einem $K \in \N$ abhängt
\begin{equation}\label{Definition-R0} R_0(K):=\{\sigma+it \in \C \,| \, 1- \frac{\log \log \mathscr{L}}{3\mathscr{L_1}} \leq \sigma \leq 1, \; |t| \leq K\} .
\end{equation}\\
Für $K$ werden wir später die von $q$ abhängige Konstante $L_1=L_1(q)$ aus folgendem Lemma wählen.
\begin{lemma}[Lemma 6.1 aus \cite{Hea92}]\label{hb-lemma-6-1}
Es gibt ein $q_0$ in $\N$ und eine von $q$ abhängige positive ganze Zahl $L_1$ mit $L_1 \leq \frac{1}{10}\mathscr{L}$, so dass für $q \geq q_0$ die Funktion
 $$\prod_{\chi \; mod \; q} L(s,\chi)$$ keine Nullstellen hat in den beiden Vierecken
$$R_0(10 L_1) \backslash R_0(L_1). $$
\end{lemma}

Mit Hilfe von Lemma \ref{hb-lemma-6-1} formulieren wir jetzt Lemma \ref{hb-lemma-5-2} insoweit neu, wie wir es später für die Anwendung brauchen werden. Dabei speisen wir die Tatsache mit ein, dass man bei Lemma \ref{hb-lemma-5-2} gewisse endlich viele Nullstellen $\rho$, die nicht der Bedingung $|1+it-\rho|\leq \delta$ genügen, mit in die Summe reinnehmen darf, weil die entsprechenden zusätzlichen Terme dann vernachlässigbar klein sind. Dieser Sachverhalt wird in \cite[S.287 Mitte]{Hea92} kurz angesprochen, aber nicht explizit ausformuliert.\\
Ein weiterer zentraler Aspekt, der in die Neuformulierung miteinfließen wird, ist der Fakt, dass man, wenn $f$ die Bedingung 2 erfüllt, jegliche Nullstellen $\rho$ mit $Re \{\rho\} \leq Re  \{s\}$ in Lemma \ref{hb-lemma-5-2} weglassen darf.

\begin{lemma}["`Arbeitsversion"' von Lemma 5.2 aus \cite{Hea92}]\label{hb-lemma-5-2-extended}
Sei $\chi \neq \chi_0$ ein Charakter $mod \; q$ und seien $L_1$ und $R:=R_0(L_1)$ (abhängig von $q$) so definiert wie in Lemma \ref{hb-lemma-6-1} und (\ref{Definition-R0}). Weiterhin sei $s \in R_0 (9 L_1)$ und die Anzahl der Nullstellen\footnote{Damit ist die Anzahl der Nullstellen gezählt mit ihrer Vielfachheit gemeint. Allgemein gilt für die gesamte Arbeit, dass wir eine Nullstelle $\rho$ der Vielfachheit $k\in \N$ so behandeln, als handelte es sich um $k$ separate Nullstellen, die eben alle den gleichen Wert haben.} von $L(s,\chi)$ in $R$, die rechts von $s$ liegen, sei maximal $10$, genauer:
$$\#{'} A_1 := \#{'} \{\rho \in R \,| \, L(\rho,\chi)=0, \; Re \{s\} < Re  \{\rho\} \} \leq 10, $$
wobei das "`$\#{'}$"' die Zählung mit Vielfachheit andeutet. Weiter sei $A_2$ eine beliebige Menge mit
$$A_2 \subseteq \{\rho \in R \,| \, L(\rho,\chi)=0, \; Re  \{\rho\} \leq Re \{s\} \} \; \text{ und } \;\#{'} A_2 \leq 10 .$$
Nun sei $f$ eine Funktion, die Bedingung 1 und 2 erfüllt. Dann gilt für $\varepsilon >0$ und $q \geq q_0(f,\varepsilon)$, dass
$$ K(s,\chi) \leq -\mathscr{L} \sum_{\rho \in A_1 \cup A_2} Re  \{F((s-\rho)\mathscr{L})\} + \frac{f(0)}{2} \phi \mathscr{L}+ \varepsilon \mathscr{L} ,$$
  wobei
  $$K(s,\chi)=\sum_{n=1}^{\infty} \Lambda (n) Re  \left\{\frac{\chi(n)}{n^s} \right\} f(\mathscr{L}^{-1}\log n).$$
\end{lemma}
\begin{beweis}(vergleiche \cite[§6]{Hea92})\\
Sei $\varepsilon>0$ und $s \in R_0(9 L_1)$. Wir benutzen die Aussage in Lemma \ref{hb-lemma-5-2} mit $\frac{\varepsilon}{2}$ anstelle von $\varepsilon$ und bekommen feste Parameter $\delta=\delta(f)$ und $q_1=q_1(f,\varepsilon)$. In der Summe müssen wir die Bedingung $\rho \in \{w \in \C \,| \, |1+it-w| \leq \delta \}$ durch $\rho \in A_1 \cup A_2$ ersetzen. Der Beweis der folgenden zwei Ungleichungen liefert also die Aussage des Lemmas.
\begin{eqnarray*}
   -\sum_{|1+it-\rho| \leq \delta} Re  \{F((s-\rho)\mathscr{L})\} &\leq& -\sum_{|1+it-\rho| \leq \delta,\, \rho \in R, \atop (Re \{\rho\} \leq Re \{s\} \Rightarrow \rho \in A_2) }Re  \{F((s-\rho)\mathscr{L})\} \\
   &\leq& -\sum_{\rho \in R, \atop (Re \{\rho\} \leq Re \{s\} \Rightarrow \rho \in A_2)}Re  \{F((s-\rho)\mathscr{L})\} +\frac{\varepsilon}{2}\\
   &=& -\sum_{\rho \in A_1 \cup A_2}Re  \{F((s-\rho)\mathscr{L})\}  +\frac{\varepsilon}{2} .
\end{eqnarray*}
\textbf{Zur ersten Ungleichung.} Sei $\rho$ gegeben mit $|1+it-\rho| \leq \delta$. Wegen der Voraussetzung $s \in R_0(9 L_1)$ und $\delta < 1 \leq L_1$ folgt $|Im \{\rho\}| \leq |t|+\delta < 10 L_1$. Aus Lemma \ref{hb-lemma-6-1} folgt dann, dass
$$Re  \{\rho\} < 1- \frac{\log \log \mathscr{L}}{3 \mathscr{L}} \; \text{ oder } \; \rho \in R .$$
Tritt der zweite Fall auf, also $\rho \in R$, so müssen wir nichts weiter machen. Sollte der erste Fall eintreten, dann gilt $Re \{s-\rho\} \geq 0$ und nach Bedingung 2 folgt $-Re  \{F((s-\rho)\mathscr{L})\} \leq 0$, so dass wir nach oben abschätzen können, indem wir diesen Term weglassen.\\
Zusätzlich lassen wir jegliche $\rho$ weg, welche die beiden Bedingungen $Re \{\rho\} \leq Re \{s\}$ und $\rho \notin A_2$ erfüllen, was wieder nach Bedingung 2 geht. Es folgt die Ungleichung.\\

\noindent \textbf{Zur zweiten Ungleichung.} Sei $\rho \in R$ eine Nullstelle von $L(s,\chi)$ mit $|1+it-\rho| > \delta$. Ist $Re \{\rho\} \leq Re \{s\}$, so erfülle $\rho$ die Bedingung $\rho \in A_2$. Wir müssen zeigen, dass die zu diesen Nullstellen $\rho$ gehörigen Terme einen Beitrag von maximal $\frac{\varepsilon}{2}$ liefern.
Zunächst haben wir mit der Dreiecksungleichung
$$|Im\{s-\rho\}|=|i(t-Im\{\rho\})| \geq |1+it-\rho| - |1- Re \{\rho\}| > \delta - \frac{\log \log \mathscr{L}}{3\mathscr{L}} > \frac{\delta}{2} $$
für $q \geq q_2=q_2(\delta)$. Damit ist $|Im\{(s-\rho)\mathscr{L}\}|\geq \frac{\delta}{2}\mathscr{L}$ und $|Re \{(s-\rho)\mathscr{L}\}| \leq \frac{1}{3} \log \log \mathscr{L} $. Schreiben wir nun $(s-\rho)\mathscr{L}=u+iv$, dann folgt mit partieller Integration\footnote{Normalerweise wird für die partielle Integration vorausgesetzt, dass die zu differenzierende Funktion auf dem Intervall \emph{einschließlich} der Integrationsgrenzen definiert und differenzierbar ist. Unsere Art die partielle Integration anzuwenden, obwohl $f'$ nicht in $0$ und $x_0$ definiert ist, ist trotzdem korrekt, was man schnell aus Bedingung 1 folgern kann (der gleiche Sachverhalt wird in \cite[S.280-281]{Hea92} benutzt). Eine Begründung wäre $f$ außerhalb von $[0,x_0]$ passend zu einer Funktion $\overline{f}$ fortzusetzen, so dass $\overline{f}$ auf $[0,x_0]$ differenzierbar und $\overline{f}'=f'$ auf $(0,x_0)$ ist.}
 (siehe auch \cite[S.280-281]{Hea92})
\begin{eqnarray}
  Re  \{F((s-\rho)\mathscr{L})\}&=& \int_0^{x_0} f(t) e^{-u t} \cos(vt) \,dt \nonumber \\
  &=& f(x_0) e^{-u x_0} \frac{\sin(v x_0)}{v}- f(0) e^0 \frac{\sin (0)}{v}  \nonumber \\
  & & -\int_0^{x_0} \Big{(}f'(t) e^{-u t}- u f(t) e^{-u t} \Big{)}\frac{\sin (v t)}{v} \,dt \label{partielle-integration-von-Re-F} \\
  &\ll_{f,\delta} & \mathscr{L}^{-1} (\log^{x_0} \mathscr{L}) \log \log \mathscr{L} \nonumber .
\end{eqnarray}
Dabei haben wir benutzt, dass $f'$ auf $(0,x_0)$ beschränkt ist, was weiter oben erwähnt wurde. Jede anfängliche Nullstelle $\rho$ liefert also für $q\rightarrow \infty$ einen Beitrag von $o_{f,\delta}(1)$. Beachte, dass $f$ und $\delta$ fest und unabhängig von $q$ sind. Da wegen der Voraussetzungen an $A_1$ und $A_2$ nur maximal $20$ solche $\rho$ existieren, gibt es ein $q_3=q_3(f,\delta,\varepsilon)$, so dass der Beitrag dieser Nullstellen insgesamt $\leq \frac{\varepsilon}{2}$ ist, wenn $q \geq q_3$. Für das Lemma wählen wir $q_0=\max \{q_1,q_2,q_3\}$ und wir sind fertig.
\end{beweis}

\noindent Die Zahl $10$ im vorigen Lemma war willkürlich gewählt. Für später hätte es auch $3$ getan. Der Beweis hätte mit jeder festen Zahl $N \in \N$ geklappt.\\

Wir benennen jetzt einige der Nullstellen in diesem Viereck $R$ so, wie es Heath-Brown in \cite[S.285 und 287]{Hea92} macht. Man beachte dabei, dass wenn wir im Folgenden von einer speziellen Nullstelle $\rho$ einer Funktion $L(s,\chi)$ sprechen, so tun wir das unter der impliziten Voraussetzung, dass diese Nullstelle überhaupt existiert.\\
Zu einem festen $q$ betrachte man sämtliche in $R$ liegenden Nullstellen $\rho$ der Funktion\footnote{Beachtet man die nullstellenfreie Region der Riemannschen $\zeta$-Funktion, welche sich bekanntlich auf $L(s,\chi_0)$ überträgt, so folgt, dass $L(s,\chi_0)$ keine Nullstellen in $R$ hat, wenn $q$ groß genug ist.}
\begin{equation}\label{produkt-nullstellen}P(s):=\prod_{\chi \; mod \; q, \atop \chi \neq \chi_0} L(s,\chi) . \end{equation}
Sei $\rho_1$ eine Nullstelle von $P(s)$ in $R$ mit maximalem Realteil und $\chi_1$ ein dazugehöriger Charakter (also $L(\rho_1,\chi_1)=0$). Wähle nun eine Nullstelle $\rho_2$ von $$P(s) (L(s,\chi_1) L(s,\overline{\chi_1}))^{-1}$$ in $R$ mit maximalem Realteil und schreibe $\chi_2$ für den zugehörigen Charakter\footnote{Zusätzlich zu den Nullstellen von $L(s,\chi_1)$ wurden auch jene von $L(s,\overline{\chi_1})$ "`herausgenommen"'. Das liegt daran, dass $\rho \in R$ genau dann eine Nullstelle von $L(s,\chi)$ ist, wenn $\overline{\rho}$ eine Nullstelle von $L(s,\overline{\chi})$ ist. Kennt man also die Nullstellen von $L(s,\chi_1)$, dann auch jene von $L(s,\overline{\chi_1})$.}. Mache so weiter, bis es in $R$ keine weiteren Nullstellen mehr gibt. Mit anderen Worten betrachtet man im ($k+1$)-ten Schritt die Nullstellen von $$P(s) (L(s,\chi_1) L(s,\overline{\chi_1}) \ldots L(s,\chi_k) L(s,\overline{\chi_k}))^{-1} $$ in $R$ und wählt $\rho_{k+1}$ unter diesen Nullstellen so, dass $Re\{\rho_{k+1}\}$ maximal ist. Wir halten die zentralen Eigenschaften der eben gewählten Nullstellen und Charaktere fest. Das wäre einerseits
$$\chi_i \neq \chi_j,\overline{\chi_j} \text{ für } i \neq j$$ und andererseits wegen Lemma \ref{hb-lemma-6-1} die Aussage für $q\geq q_0$, dass wenn $\chi \neq \chi_0, \chi_i, \overline{\chi_i}$ ist für $1\leq i < k$ und $L(\rho,\chi)=0$, dann gilt
\begin{equation}\label{Eigenschaft1-Nullstellen}Re \{\rho\} \leq Re \{\rho_k\} \; \text{ oder }\; |Im\{\rho\}|\geq 10 L_1 . \end{equation}
Für die spätere Anwendung setzen wir \label{Definition-rho-k}$$\rho_k=\beta_k + i \gamma_k, \; \beta_k=1-\mathscr{L}^{-1}\lambda_k, \; \gamma_k=\mathscr{L}^{-1} \mu_k .$$
Wir betrachten noch eine weitere Nullstelle. Angenommen $L(s,\chi_1)$ hat eine weitere Nullstelle $\rho' \neq \rho_1$ in $R$ oder die Nullstellenordnung der Nullstelle $\rho_1$ von $L(s,\chi_1)$ ist $\geq 2$. Dann wählen wir folgendermaßen eine Nullstelle $\rho' \in R$ von $L(s,\chi_1)$:
\begin{itemize}
  \item Fall 1: Ist die Nullstellenordnung von $\rho_1$ größer oder gleich $2$, so wählen wir $\rho'=\rho_1$.
  \item Fall 2: Sind wir nicht in Fall 1 und ist $\chi_1$ reell und $\rho_1$ komplex, so wählen wir $\rho'$ unter den Nullstellen in $R \backslash \{ \rho_1, \overline{\rho_1}\}$ so, dass $Re\{\rho'\}$ maximal ist.
  \item Fall 3: Sind wir nicht in Fall 1 oder Fall 2, so wählen wir $\rho'$ unter den Nullstellen in $R \backslash \{\rho_1 \}$ so, dass $Re\{\rho'\}$ maximal ist.
\end{itemize}
In Analogie zu vorher setzen wir
$$\rho'=\beta' + i \gamma', \; \beta'=1-\mathscr{L}^{-1}\lambda', \;\gamma'=\mathscr{L}^{-1} \mu' .$$

In \cite{Hea92} werden Abschätzungen für die Realteile der Nullstellen $\rho_1, \rho_2, \rho_3$ und $\rho'$ bewiesen. Wie in §\ref{para-methoden-fuer-nullstellenfreie-regionen} angeschnitten wurde, wird dies durch zwei verschiedene Methoden erreicht. In dieser Arbeit werden wir nur die zweite in §\ref{para-methoden-fuer-nullstellenfreie-regionen} erwähnte Methode benutzen. Wir möchten die Struktur dieser Methode exemplarisch hier vorstellen, indem wir den Beweis eines Teils von Prinzip 1 damit skizzieren:\\

\noindent \textbf{Lemma}\emph{
 Angenommen es ist $\chi_1^2 \neq \chi_0$ und $q \geq q_0$. Dann gibt es eine Konstante $C>0$, so dass $$\prod_{\chi \; mod \; q, \atop \chi \neq \chi_0} L(s,\chi) $$ keine Nullstellen hat in $\sigma \geq 1-\mathscr{L}^{-1}C, \; |t| \leq 1 . $}

\noindent\textbf{Beweis.} Da für hinreichend großes $q$ die Beziehung $$R = R_0(L_1) \supseteq \{\sigma+it \in \C \,| \, \sigma \geq 1-\mathscr{L}^{-1}C, |t| \leq 1 \}$$
gilt, folgt das Lemma aus der Aussage $\lambda_1 > C$. Dies gilt es nun zu beweisen.\\

\noindent \textbf{1. Schritt (Anfangsungleichung)}\\ Man wähle eine "`scharfe Anfangsungleichung"' aus, beispielsweise
\begin{equation}\label{trigonometric-inequality}0 \leq 2(1+\cos (\theta))^2 = 3+4\cos \theta + \cos 2\theta .\end{equation}
Mit $\chi_1(n)=e^{i \arg (\chi_1(n))}$ folgt $Re  \{\chi_1(n) n^{-i \gamma_1}\}=\cos (\arg(\chi_1(n))-\gamma_1 \log n)$. Setzt man also $\theta=\arg(\chi_1(n))- \gamma_1 \log n$ dann bekommt man
\begin{equation}\label{trigonometric-inequality-2}0 \leq 3 + 4 Re  \left\{\frac{\chi_1(n)}{n^{i\gamma_1}} \right\} + Re  \left\{\frac{\chi_1^2(n)}{n^{2i \gamma_1}} \right\}.  \end{equation}
Sei jetzt $f$ eine Funktion, die Bedingung 1 und 2 erfüllt und $\beta \in \R$. Dann multipliziere die letzte Ungleichung mit $\chi_0(n) \Lambda(n) f(\mathscr{L}^{-1}\log n) n^{-\beta}$, summiere über $n$ und erhalte das Analogon zu (\ref{zeta-ungleichung}), nämlich (für die Definition von $K(s,\chi)$ siehe Lemma \ref{hb-lemma-5-2-extended})
\begin{equation}\label{trig-inequality-corollary}0 \leq 3 K(\beta,\chi_0) + 4 K(\beta +i\gamma_1,\chi_1) + K(\beta+2i\gamma_1,\chi_1^2). \end{equation}

Später bleibt die Aufgabe die Anfangsungleichung optimal zu wählen. Beispielsweise hätte man
$0\leq (k+\cos \theta)^2 (l+\cos \theta)^2=\ldots$ wählen können, wobei $k,l \geq 0$ Parameter sind, die im Nachhinein optimal gewählt werden. Ein Nachteil dieser komplizierteren Ungleichung ist die erhöhte Anzahl der resultierenden "`$K$-Terme"', welche dann in der zu (\ref{trig-inequality-corollary}) analogen Ungleichung auftauchen. Dies zieht eine notwendige Diskussion mehrerer Fälle nach sich und führt nicht notwendigerweise zu besseren Ergebnissen (in diesem speziellen Fall aber später doch).\\

\noindent \textbf{2. Schritt ($\beta$ wählen und mittels Lemmata Ungleichung aufstellen)}\\
Wir müssen einen Wert für $\beta$ festlegen und wählen $\beta=\beta_1$ (später wird $\beta=\max \{\beta_2, \beta'\}$ besser sein, dies erfordert aber Zusatzüberlegungen). Aus Lemma \ref{hb-lemma-5-3} folgt für ein $\varepsilon>0$ und $q$ groß genug, dass
$$K(\beta_1,\chi_0) \leq \mathscr{L} F(-\lambda_1)+\varepsilon \mathscr{L}. $$
Was die Abschätzung von $K(\beta_1 +i \gamma_1,\chi_1)$ bzw. $K(\beta_1 +2i\gamma_1,\chi_1^2)$  mittels Lemma \ref{hb-lemma-5-2-extended} angeht (hier brauchen wir $\chi_1^2 \neq \chi_0$), so ist dies recht leicht, da wir $\beta=\beta_1$ gewählt haben. Berücksichtigen wir die Definition von $\rho_1$, so ist nämlich im letzteren Lemma die auftauchende Menge $A_1$ beides Mal leer und wir wählen $A_2=\{\rho_1 \}$ bzw. $A_2=\emptyset$. Mit $\phi \leq \frac{1}{3}$ folgt
$$K(\beta_1 +i \gamma_1,\chi_1) \leq (-F(0)+\frac{f(0)}{6}+\varepsilon )\mathscr{L}$$
bzw.
$$K(\beta_1 +2i \gamma_1,\chi_1^2)\leq (\frac{f(0)}{6}+\varepsilon )\mathscr{L} .$$
Nach Anpassung des $\varepsilon$ erhält man insgesamt
\begin{equation}\label{beispiel-l1-abschaetzung}0 \leq 3F(-\lambda_1)-4F(0)+ \frac{5}{6}f(0)+\varepsilon.\end{equation}\\

\noindent \textbf{3. Schritt (Ausrechnen und optimieren)}\\ Wäre nun $\lambda_1 \leq \lambda_{12}$ für einen konkreten Wert $\lambda_{12}>0$, dann würde wegen der Monotonie von $F(-\lambda_1)$ folgen, dass
$$-\varepsilon \leq 3F(-\lambda_{12})-4F(0)+ \frac{5}{6}f(0). $$
Fänden wir nun eine Funktion $f$, die Bedingung 1 und 2 erfüllt, und für die die letzte rechte Seite negativ ist, so ist dies für hinreichend kleines $\varepsilon>0$ und $q \geq q_0(\varepsilon)$ ein Widerspruch. Also haben wir unter diesen Voraussetzungen bewiesen, dass
$\lambda_1 > \lambda_{12}=:C$. Wählt man speziell die zu $\gamma=1.9$ gehörige Funktion $f$ (genaueres dazu im nächsten Abschnitt), so bekommt man $C= 0.144$. Die Optimierung nach der Funktion $f$ erfolgt dabei am besten mittels Computer. Dazu müssen wir aber erst wissen, welche Funktionen $f$ überhaupt zur Verfügung stehen. Deswegen beschäftigen wir uns im nächsten Abschnitt mit

\subsection{Funktionen, die Bedingung 1 und 2 erfüllen (\cite[§7]{Hea92})}\label{hb-para-7}
Zur Anwendung des zentralen Lemmas \ref{hb-lemma-5-2-extended} benötigt man Funktionen $f$, die Bedingung 1 und 2 erfüllen. Bedingung 1 stellt keine Hürden, Bedingung 2 aber sehr wohl. Heath-Brown kreiert in \cite[§7]{Hea92} passende Funktionen $f$, indem er
\begin{equation}\label{hb-para7-ansatz-fuer-f}f(t)=\int_{-\infty}^{\infty} g(x) g(t-x) \,dx \end{equation}
wählt. Dabei sei die Funktion $g(x): \R \rightarrow \R$ stetig und für $x \in \R$ erfülle sie $g(x)=g(-x)$, sowie $g(x)\geq 0$. Außerdem sei $g(x)=0$ für $|x| \geq \gamma$, wobei $\gamma$ eine positive Konstante sei. \\
Um zu zeigen, dass $f$ Bedingung 2 erfüllt, wird folgendes funktionentheoretisches Lemma benötigt, welches ein Korollar vom Maximumsprinzip holomorpher Funktionen ist.

\begin{lemma}[Lemma 4.1 aus \cite{Hea92}]\label{hb-lemma-4-1}
Seien $F_1(z)$ und $F_2(z)$ holomorphe Funktionen in $\mathscr{H}=\{z \in \C\,| \, Re \{z\} \geq 0 \}$ und sei $Re \{F_1(z)\} \geq |F_2(z)|$ auf $Re\{z\}=0$. Weiterhin seien $F_1$ und $F_2$ in $\mathscr{H}$ gleichmäßig konvergent gegen $0$ für $|z|\rightarrow \infty$. Dann ist $Re \{F_1(z)\} \geq |F_2(z)|$ auf ganz $\mathscr{H}$.
\end{lemma}
\noindent Damit beweist man

\begin{lemma}[S.289 in \cite{Hea92}]\label{hb-lemma-7-0}
  Sei $f$ so definiert wie in (\ref{hb-para7-ansatz-fuer-f}) und $g$ erfülle die obigen Bedingungen. Zusätzlich erfülle $f$ Bedingung 1. Dann erfüllt $f$ Bedingung 2.
\end{lemma}
\begin{beweis}
Sei $F$ die Laplace-Transformierte von $f$. Mit den Eigenschaften von $g$ und der Gleichung $\cos x + \cos y = 2 \cos \frac{x+y}{2} \cos \frac{x-y}{2}$ schlussfolgern wir für $y \in \R$
\begin{eqnarray*}
Re  \{F(iy)\} &=&  \int_{0}^{\infty} f(t) \cos(yt) \,dt\\
&=& \int_{0}^{\infty} \int_{-\infty}^{\infty} g(x)g(t-x)\cos(yt)\,dx \,dt \\
&=& \int_{0}^{\infty} \int_{0}^{\infty} g(x)g(t-x)\cos(yt)\,dx \,dt +\int_{0}^{\infty} \int_{-\infty}^{0} g(x)g(t-x)\cos(yt)\,dx \,dt \\
&=& \int_{0}^{\infty} \int_{0}^{\infty} g(x)g(t-x)\cos(yt)\,dt \,dx +\int_{0}^{\infty} \int_{0}^{\infty} g(x)g(t+x)\cos(yt)\,dt \,dx \\
&=& \int_{0}^{\infty} \int_{-x}^{\infty} g(x)g(l)\cos(y(l+x))\,dl \,dx +\int_{0}^{\infty} \int_{x}^{\infty} g(x)g(l)\cos(y(l-x))\,dl \,dx \\
&=& \int_{0}^{\infty} \int_{0}^{\infty} g(x)g(l)\cos(y(l+x))\,dl \,dx +\int_{0}^{\infty} \int_{0}^{\infty} g(x)g(l)\cos(y(l-x))\,dl \,dx\\
&=& 2\left(\int_{0}^{\infty} g(t) \cos(t y) \,dt \right)^2 \geq 0.
\end{eqnarray*}
In Lemma \ref{hb-lemma-4-1} wählen wir nun $F_1=F$ und $F_2=0$ und beachten, dass $F(z)$ in $\mathscr{H}$ gleichmäßig gegen $0$ geht, wenn $|z| \rightarrow \infty$. Letzteres folgt z.B. aus \cite[(5.3),(5.4)]{Hea92} (beachte, dass $f$ Bedingung 1 erfüllt). Damit gilt das Lemma.
  \end{beweis}

In der Menge der oben genannten $f$ findet Heath-Brown optimale Elemente durch einen "`Variationsansatz"': Ziel ist es ja Ausdrücke wie die rechte Seite von (\ref{beispiel-l1-abschaetzung}) bei festem $\lambda_1$ durch geschickte Wahl eines $f$ zu minimieren. Geht man davon aus, dass die Funktion $$f_1(t)=\int_{-\infty}^{\infty} g_1(x) g_1(t-x) \,dx$$ in diesem Sinne optimal ist, so bedeutet dies, dass kleine Änderungen an $g_1$ dazu führen, dass sich die rechte Seite von (\ref{beispiel-l1-abschaetzung}) vergrößert. Aus einem solchen Sachverhalt lässt sich u.a. extrahieren, dass $g_1$ eine gewisse Differentialgleichung 2. Ordnung erfüllen sollte (\cite[(7.10)]{Hea92}). Die Analyse führt letztendlich zu Funktionen $f$, die aus den folgenden Funktionen $g$ resultieren\footnote{Man beachte dabei, dass wenn man für ein $\theta \in \R \backslash \{0\}$ die Funktion $g_2(x)=\theta g_1(x)$ anstelle von $g_1(x)$ verwendet, dies zu $f_2(t)= \theta^2 f_1(t)$ führt. Es folgt für die Laplace-Transformierten $F_2(z)=\theta^2 F_1(z)$. Also sind die beiden Ungleichungen vom Typ (\ref{beispiel-l1-abschaetzung}), die man aus diesen beiden Funktionen erhält, nach Anpassung der $\varepsilon$ identisch.}  ($x \in (-\gamma,\gamma)$)
\begin{eqnarray}
 g(x)&=& \gamma^2-x^2, \nonumber \\
 g(x)&=& \cos(c_1 x)- \cos(c_1 \gamma) , \label{g-version-2} \\
 g(x)&=& \cosh(c_1 \gamma)- \cosh(c_1 x), \label{g-version-3}
\end{eqnarray}
wobei $g(x)=0$ gesetzt wird für $x \notin (-\gamma,\gamma)$ und $c_1, \, \gamma$ gewisse Parameter sind (vergleiche \cite[Lemma 7.1-7.5]{Hea92}). Für unsere Untersuchungen werden wir nun nur diejenigen $f$ benutzen, die aus $g(x)=\gamma^2-x^2$ resultieren. Für diese Funktionen $f$ gilt $f(t)=0$ für $t \geq 2 \gamma$ und für $t \in [0,2\gamma)$
\begin{eqnarray}
  f(t)&=&\int_{t-\gamma}^{\gamma} g(x) g(t-x) \, dx =-\frac{1}{30}t^5 +\frac{2 \gamma^2}{3}t^3 - \frac{4 \gamma^3}{3} t^2 + \frac{16\gamma^5}{15} . \label{explizites-f}
\end{eqnarray}
Es folgt, dass $f$ Bedingung 1 erfüllt, also nach obigem Lemma auch Bedingung 2. Für die zugehörige Laplace-Transformierte $F$ gilt $F(0)=\frac{8 \gamma^6}{9}$ und für $z \in \C \backslash \{0\}$ folgt mittels partieller Integration (vergleiche \cite[S.312]{Hea92})
\begin{eqnarray}
  F(z)&=&\frac{16 \gamma^5}{15} z^{-1}-\frac{8 \gamma^3}{3}z^{-3} + 4\gamma^2(1+ e^{-2 \gamma z})z^{-4}+ 4z^{-6} (-1+e^{-2 \gamma z} +2 \gamma z e^{-2 \gamma z})  .  \nonumber \\
  & & \label{explizites-F}
\end{eqnarray}
Wir weisen daraufhin, dass für die gesamte restliche Diplomarbeit eine Funktion $f$ immer für die Funktion in (\ref{explizites-f}) stehen wird für ein gewisses $\gamma>0$. An gewissen Stellen werden wir später nur das $\gamma>0$ anstelle des $f$ angeben, wodurch letzteres dann ja eindeutig bestimmt ist.\\

\noindent Der größte Vorteil der obigen speziellen Wahl für $g$ sind die später viel schneller durchführbaren Computerberechnungen, welche aufgrund der Einfachheit von $f(t)$ und $F(z)$ resultieren. Zusätzlich lässt sich für großes $|Im\{z\}|$ die Funktion $F(z)$ viel leichter nach oben abschätzen. Dies wird vorteilhaft sein.
Außerdem zeigen ein paar numerische Proberechnungen, dass mit denjenigen Funktionen $f$, die aus (\ref{g-version-2}) oder (\ref{g-version-3}) resultieren, höchstens nur minimal bessere Werte erzielbar sind, als mit unserer Wahl für die Funktion $g$. Dabei muss man beachten, dass die Graphen aller drei Versionen von $g$ "`qualitativ recht identisch aussehen"'. Ein paar weitere Bemerkungen zu Funktionen, die Bedingung 1 und 2 erfüllen, sind in §\ref{para-vb1} zu finden.

\subsection{Abschätzungen für $\lambda_2$, $\lambda'$, wenn $\rho_1$, $\chi_1$ reell (\cite[§8]{Hea92})}
In diesem Abschnitt möchten wir zwei Tabellen aus \cite[§8]{Hea92} vorstellen. Generalvoraussetzung ist dabei, dass $\chi_1$ und $\rho_1$ beide reell sind (insbesondere existieren). Weiterhin wird vorausgesetzt, dass $\lambda_1 \geq \lambda(\varepsilon)$ für ein gewisses $\lambda(\varepsilon)>0$ (ansonsten könnte man nach \cite{Hea90} im Theorem von Linnik $L=3+\varepsilon$ nehmen). \\

\noindent Aus \cite[§8]{Hea92} werden wir später u.a. die beiden Tabellen "`Tabelle H4"':=\cite[Table 4 (§8)]{Hea92} und "`Tabelle H7"':=\cite[Table 7 (§8)]{Hea92} benötigen. Wir möchten hier zusammenfassen, was diese Tabellen aussagen.\\

\textbf{Erklärung der Tabelle H4:} Diese Tabelle ist ein Korollar von \cite[Lemma 8.3]{Hea92}, welches unter der Voraussetzung gilt, dass $\rho'$ komplex ist. Ist jedoch $\rho'$ reell, so liefert \cite[Lemma 8.2]{Hea92} bessere Werte, so dass also letztendlich diese Tabelle für beide Fälle gilt. Was die Tabelle aussagt, wird kurz vor der Tabelle erklärt: Beispielsweise besagt der Eintrag $0.7,1.724,4.5,0.85$, dass wenn $\lambda_1 \leq 0.70$, dann ist $\lambda' \geq 1.724$. Zuerst wurde dabei letztere Aussage nur für den Fall bewiesen, dass $\lambda_1 \in [0.65,0.70]$. Da aber für kleinere $\lambda_1$ noch bessere Abschätzungen gelten, gilt letztendlich die Abschätzung $\lambda' \geq 1.724$ immer wenn $\lambda_1 \leq 0.70$ (man beachte, dass die oben genannten Generalvoraussetzungen stets angenommen werden).\\

\textbf{Erklärung der Tabelle H7:} Hier gilt Analoges wie bei der vorigen Tabelle. Diesmal wurde diese Tabelle unter der Voraussetzung $\chi_2^4 \neq \chi_0$ und $\lambda_2 \leq \lambda'$ bewiesen. Für den Fall $\chi_2^4 = \chi_0$ bekommt man vermöge \cite[Table 6 (§8)]{Hea92} bessere Werte (bis auf den Fall $\lambda_1 \leq 0.10$). Genauso bekommt man bessere Werte durch \cite[Table 2-4 (§8)]{Hea92}, wenn $\lambda_2>\lambda'$ ist. Insgesamt folgt wieder, dass diese Tabelle in allen Fällen gültig ist (ausgenommen wenn $\lambda_1\leq 0.10$).\\
Im Gegensatz zu Tabelle H4 wurde für diese sofort bewiesen, dass z.B. (2. Zeile) aus $\lambda_1 \leq 0.12$ folgt, dass $\lambda_2 \geq 2.56$.\\

\noindent Die in diesen Tabellen auftauchenden Werte für $\lambda,\, a$ und $K$ sind Parameter, die bei den jeweiligen Rechnungen benutzt wurden. Sie sind nicht Inhalt der letztendlich bewiesenen Aussage.

Wir weisen auch darauf hin, dass die Werte in Tabelle H4 mit Hilfe von Verbesserungspotential Nr. 2 (\cite[S.332]{Hea92}) hätten verbessert werden können. Dies ist für den Beweis von Theorem \ref{haupttheorem} aber nicht notwendig, also werden wir dies nicht tun.

\section{Das Theorem von Linnik und die Nullstellen von $L(s,\chi)$}\label{para-beweisdeshaupttheorems-einfuehrung}
Ein faszinierender Aspekt der analytischen Zahlentheorie ist sicher die Tatsache, dass man "`diskrete"' Phänomene, wie die natürlichen Zahlen und spezieller die Primzahlen, mit "`kontinuierlichen"' Funktionen behandeln kann. Wir stellen hier den Zusammenhang zwischen den Primzahlen in arithmetischen Progressionen und den Nullstellen der Dirichletschen L-Funktionen dar. Dies wird in drei wesentlichen Schritten geschehen. Wir folgen dabei präzise der Herleitung in \cite[§13]{Hea92}. Dabei fügen wir keine Beweise ein, sondern zitieren und kommentieren nur die wichtigen Zwischenresultate. \\

\noindent Seien $a,q \in \N$ und $(a,q)=1$. Man beginnt mit der Summe \label{Definition-Sigma} $$\Sigma=\sum_{p \equiv a \; mod \; q} \frac{\log p}{p}h(\mathscr{L}^{-1}\log p),$$
wobei das Gewicht $h(x)$ definiert ist als
$$ h(x):= \left\lbrace \begin{array}{ll}
0 & \text{ für } \; x\leq L-2K, \\
x-(L-2K) & \text{ für } \; L-2K \leq x \leq L-K, \\
L - x & \text{ für } \; L-K \leq x \leq L,\\
0 & \text{ für } \; x\geq L \\
\end{array} \right.
$$
mit Konstanten $L,\,K>0$ und $L-2K>0$. Wir müssen zeigen, dass für $q \geq q_0$ die Summe $\Sigma$ unabhängig vom gewählten $a$ positiv ist. Daraus folgt $P(q) \leq q^L$ für $q \geq q_0$ und daraus wiederum $P(q) \ll q^L$ für alle $q$.

Um mit der Summe $\Sigma$ besser arbeiten zu können, ist es vorteilhaft, sie über die Primzahlpotenzen gehen zu lassen anstatt nur über die Primzahlen. Dies ist mit einem vernachlässigbaren Fehler möglich:
$$\Sigma=\sum_{n \equiv a \; mod \; q} \frac{\Lambda(n)}{n} h(\mathscr{L}^{-1}\log n) +O( q^{-\frac{L-2K}{2}}) . $$

\noindent Nun kommt der \textbf{erste} wesentliche Schritt, nämlich die Einführung der Charaktere $\chi$, wodurch man die Summation über eine einzige Restklasse $a \; mod \; q$ in den Griff bekommt. Diese beeindruckende Methode geht zurück auf Dirichlet. Man benutzt die bekannte Gleichung ($n \in \N$)
$$\frac{1}{\varphi (q)} \sum_{\chi \; mod \; q} \overline{\chi}(a) \chi(n)=\left\lbrace \begin{array}{ll}
1 & \text{ falls } n \equiv a \; mod \; q, \\
0 & \text{ sonst. }
\end{array} \right. $$
Es folgt
\begin{equation}\label{para13-Sigma-2}\Sigma=\frac{1}{\varphi(q)} \sum_{\chi \; mod \; q} \overline{\chi}(a) \sum_{n=1}^{\infty} \frac{\Lambda(n) \chi(n)}{n}h(\mathscr{L}^{-1}\log n) + O(q^{-\frac{L-2K}{2}}) .\end{equation}

\noindent Der \textbf{zweite} wesentliche Schritt ist die Umschreibung der inneren Summen in (\ref{para13-Sigma-2}) als komplexe Linienintegrale. Dies ist der eigentliche (und einzige) Moment, in dem die Verbindung zwischen der Zahlentheorie und der Funktionentheorie hergestellt wird. Bekanntlich gilt es sehr oft, dass die Verbindung dadurch hergestellt wird, dass man Summen als komplexe Linienintegrale aufschreibt. Das klassischste Beispiel hierfür ist wohl die Perronsche Formel, aus der man eine Summe $\sum\limits_{n \leq x} b(n)$ leicht in Zusammenhang setzt mit der zugehörigen Dirichletschen Reihe $$\sum_{n=1}^{\infty} b(n) n^{-s},$$ ohne dass diese in $s=0$ konvergieren muss (vergleiche \cite[S.26f]{Bru95}).

Mit vernachlässigbarem Fehler ersetzt man nun in der inneren Summe von (\ref{para13-Sigma-2}) den Charakter $\chi$ durch den primitiven Charakter $\chi^{\star}$, der $\chi$ induziert. Dies wird gemacht, da man die gleich in Erscheinung tretenden Funktionen $(L'/L)(s,\chi)$ für primitives $\chi$ besser handhaben kann. Man hat schließlich für $\chi \neq \chi_0$ und den zu $\chi$ gehörigen primitiven Charakter $\chi^{\star}$, dass
\begin{equation}\label{para13-Sigma-3}\mathscr{L}^{-1} \sum_{n=1}^{\infty} \frac{\Lambda (n) \chi^{\star} (n)}{n}h(\mathscr{L}^{-1} \log n)=\frac{1}{2\pi i} \int_{2-i\infty}^{2+i\infty} \left(-\frac{L'}{L}(s,\chi^{\star})\right) H((1-s)\mathscr{L}) \,ds ,\end{equation}
wobei für die Laplace-Transformierte $H$ von $h$ gilt
\begin{equation}\label{para13-definition-F}H(z)=e^{-(L-2K)z} \left(\frac{1-e^{-Kz}}{z}\right)^2 .\end{equation}
Weiterhin hat man für $\chi=\chi_0$, dass $\chi^{\star}$ immer $1$ ist, also
\begin{equation}\label{para13-Sigma-4}\mathscr{L}^{-1} \sum_{n=1}^{\infty} \frac{\Lambda(n)}{n} \chi^{\star}(n) h(\mathscr{L}^{-1} \log n)=H(0) + O (\mathscr{L}^{-1}).  \end{equation}
Letzteres beweist man mit partieller Summation und dem Primzahlsatz.\\
Kombiniert man (\ref{para13-Sigma-2}), (\ref{para13-Sigma-3}) und (\ref{para13-Sigma-4}), so hat man $\Sigma$ durch die Kombination verschiedener recht kompliziert anmaßender Terme ersetzt. Bisher ist auch noch nichts gewonnen, da man die Integranden wegen des Verhaltens von $H((1-s)\mathscr{L})$ auf $Re \{s\}=2$ keineswegs vernünftig abschätzen kann.\\

\noindent Damit kämen wir zum \textbf{dritten} Schritt. Bei den $\varphi(q)-1$ Integralen verschiebt man die Integrationslinie von $Re \{s\}=2$ auf $Re \{s\}=-\frac{1}{2}$. Mit dem Residuensatz und leichter Abschätzung des verschobenen Integrals kann man mit vernachlässigbarem Fehler das anfängliche Integral durch eine gewichtete Summe über die nicht-trivialen Nullstellen von $L(s,\chi^{\star})$ ersetzen. Diese sind die selben, wie die nicht-trivialen Nullstellen von $L(s,\chi)$. Man ist am Ziel und hat damit $\Sigma$ mit Hilfe der Nullstellen der Dirichletschen L-Funktionen abgeschätzt:

\begin{lemma}[Lemma 13.1 aus \cite{Hea92}]
Seien $L,\,K>0$ mit $L \geq 2+2K$ und $\Sigma,\, h$ und $H$ definiert wie oben. Dann gilt
$$\left|\mathscr{L}^{-1} \Sigma - \frac{H(0)}{\varphi(q)} \right|\leq \frac{1}{\varphi(q)} \left(\sum_{\chi \neq \chi_0} \sum_{\rho} |H((1-\rho)\mathscr{L})| + O(\mathscr{L}^{-1}) \right).$$
Die Summe geht dabei über die nicht-trivialen Nullstellen $\rho$ von $L(s,\chi)$.
\end{lemma}

\noindent Nun hat $H(z)$ solche Gestalt, dass in letzterer Summe der Beitrag der Nullstellen, die zu weit weg von $s=1$ sind, extrem gering ist. Benutzt man die Abschätzung \cite[(1.4)]{Hea92} für die Nullstellendichte der Dirichletschen L-Funktionen in gewissen Regionen, so folgt, dass die meisten Nullstellen in der letzten Summe mit vernachlässigbarem Fehler weggelassen werden können, genauer:

\begin{lemma}[Lemma 13.2 aus \cite{Hea92}]\label{hb-lemma-13-2}
Es gelten die gleichen Bezeichnungen wie oben. Weiterhin sei $\varepsilon>0$ und $L>3+2K$. Dann existieren positive \label{Definition-R-1}Konstanten $R_1=R_1(\varepsilon)$ und $q_0=q_0(\varepsilon)$, so dass für $q\geq q_0$ gilt
$$\left|\mathscr{L}^{-1} \Sigma - \frac{H(0)}{\varphi(q)} \right|\leq \frac{1}{\varphi(q)} \left(\sum_{\chi \neq \chi_0} \sum_{\rho}{'} |H((1-\rho)\mathscr{L})| + \varepsilon \right).$$
Die Summe $\sum\limits_{\rho}{'}$ geht dabei über die Nullstellen $\rho=\beta+i\gamma$ von $L(s,\chi)$ im Viereck
$$1-\mathscr{L}^{-1}R_1 \leq \beta \leq 1,\; \; |\gamma| \leq \mathscr{L}^{-1}R_1 . $$
\end{lemma}
\noindent Aus letzterem Lemma folgt
\begin{equation}\label{para13-Sigma-abzuschaetzender-Term}\frac{\varphi(q)}{\mathscr{L}}\Sigma \geq H(0)-\sum_{\chi \neq \chi_0} \sum_{\rho}{'} |H((1-\rho)\mathscr{L})| - \varepsilon  .\end{equation}
Um also das Theorem von Linnik mit der Konstante $L$ (diese steckt in $H$ drin) zu beweisen, muss man jetzt zeigen, dass die letzte rechte Seite positiv ist (für $q\geq q_0(\varepsilon)$).\\

\begin{bemerkung}
  Mit Hilfe von Prinzip 3 aus §\ref{Kapitel1-Problematik} können wir das Theorem von Linnik für hinreichend großes $L$ bereits jetzt schon beweisen:\\ Wenn $\rho_1$ nicht existiert, d.h. in $R$ liegen keine Nullstellen von (\ref{produkt-nullstellen}), dann erhalten wir sofort nach Wahl eines hinreichend kleinen $K>0$ die zulässige Konstante $L=3+\varepsilon$. Wenn andererseits $\rho_1$ existiert und sehr nah bei $s=1$ liegt, nämlich $\lambda_1 \leq \lambda_1(\varepsilon)$ für ein gewisses hinreichend kleines $\lambda(\varepsilon)>0$, dann können wir wieder $L=3+\varepsilon$ wählen nach \cite{Hea90}. Also sei $\lambda_1 \geq \lambda(\varepsilon)$. Es folgt mit der für alle $z \in \C$ gültigen Abschätzung $|H(z)| \leq H(Re\{z\})$ und Prinzip 3, dass
  \begin{eqnarray*}
  \frac{\varphi(q)}{\mathscr{L}}\Sigma &\geq& H(0)-\sum_{\chi \; mod \; q} N(1-R_1 \mathscr{L}^{-1}, 1, \chi) \, H((1-(1-\lambda(\varepsilon) \mathscr{L}^{-1})\mathscr{L}) - \varepsilon \\
  &\geq&K^2- c_3 e^{c_4 R_1}H(\lambda(\varepsilon)) - \varepsilon  . \end{eqnarray*}
  Nun ist $\lambda(\varepsilon)$ eine feste positive Zahl, analog sind $c_3,\, c_4$ und $R_1$ fest. Wählt man also feste Parameter $K$ und $\varepsilon$ mit $K^2-\varepsilon>0$, dann folgt wegen $\lim\limits_{L \rightarrow \infty} H(\lambda(\varepsilon))=0$ für hinreichend großes $L$ die Behauptung.
\end{bemerkung}
\vspace*{0.25cm}
Um einen besseren zulässigen Wert für $L$ zu beweisen, sind offensichtlich Aussagen zu der Lage der Nullstellen nahe $s=1$ (Stichwort: Nullstellenfreie Regionen) und Abschätzungen zu der Anzahl derselbigen (Stichwort: Abschätzungen für die Nullstellendichte) hilfreich.\\

Es erweist sich als nutzvoll die verbliebenen Nullstellen in Lemma \ref{hb-lemma-13-2} durch eine interessante Anwendung von Lemma \ref{hb-lemma-5-2} und mit Hinblick auf \cite[Theorem 11.1]{Hea92}
zu behandeln. Dies führt zu \cite[Lemma 13.3]{Hea92}, welches wir für spätere Zwecke hier zitieren.

\begin{lemma}[Lemma 13.3 aus \cite{Hea92}]\label{hb-lemma-13-3}
  Sei $\chi \neq \chi_0$ ein Charakter $mod \; q$ und $\lambda>0$. $L(s,\chi)$ habe keine Nullstellen im Viereck $$1-\mathscr{L}^{-1}\lambda < \sigma \leq 1, \; \; |t|\leq 1.$$ Seien darüberhinaus
  \begin{equation}\label{Definition-H-2} H_2(z):=\left(\frac{1-e^{-Kz}}{z}\right)^2\end{equation} und $\eta>0$. Dann gibt es ein $q_0=q_0(\eta)$, so dass für $q\geq q_0$ gilt
  $$\sum_{\rho}{'} |H_2((1-\rho)\mathscr{L})| \leq \frac{\phi}{2} \left(\frac{1-e^{-2K\lambda}}{\lambda}\right)+\frac{2K\lambda-1+e^{-2K\lambda}}{2\lambda^2} + \eta .$$
  Die Summation in $\sum\limits_{\rho}{'}$ ist dabei so definiert wie in Lemma \ref{hb-lemma-13-2}.
\end{lemma}

\noindent Eine weitere aufmerksame Analyse führt schließlich zu \cite[Lemma 15.1]{Hea92}. Letzteres liefert eine Abschätzung der rechten Seite von (\ref{para13-Sigma-abzuschaetzender-Term}) nach unten durch eine Konstante $C$ und ist genau so geschaffen, dass man nur noch die Zwischenresultate aus \cite{Hea92} einspeisen muss, um einen konkreten Wert für $C$ zu bekommen. Da wir \cite[Lemma 15.1]{Hea92} geringfügig verändern werden, zitieren wir es hier nicht, sondern kommen später darauf zurück.

\chapter{Verbesserungen der Nullstellenabschätzungen}

\section{Abschätzung gewisser Suprema gemäß \cite[§10]{Hea92}}\label{para-rechenkapitel}
Sei $f$ die Funktion in (\ref{explizites-f}) mit einem $\gamma \geq \frac{1}{2}$ und $F$ die dazugehörige Laplace-Transformierte. Später, bei der Anwendung von Verbesserungspotential Nr. 2, wird es von zentraler Notwendigkeit sein folgenden Term nach oben abzuschätzen: \label{Definition-A}
$$A(s_1,s_2,t):= Re  \Big{\{} k_1 F(-s_1+it)-k_2 F(-(s_1-s_2)+it) - k_3 F(it) \Big{\}}.$$
Dabei gelte $s_1 \in [s_{11},s_{12}]$, $s_2 \in [s_{21},s_{22}]$, $t \in \R$, $s_2 \leq s_1$, $0 \leq s_{11} \leq s_{12} \leq 4$, $0 \leq s_{21} \leq s_{22}$ und $k_i \geq 0$ für $i \in \{1,2,3\}$. Die $s_{ij}$ und $k_i$ werden später feste reelle Konstanten sein. Wir definieren noch
$$s_3:=s_1-s_2 \in \big{[}\max \{0,\, s_{11}-s_{22}\}, s_{12}-s_{21}\big{]} =: [s_{31},s_{32}].$$ Konkret werden Abschätzungen
\begin{equation}\label{rechenkapitel-zentrale-ungleichung}A_{sup}:= \sup_{s_1 \in [s_{11},s_{12}], \atop {
{s_2 \in [s_{21},s_{22}], \atop s_2 \leq s_1,\, t \in \R}}} A(s_1,s_2,t) \leq C \end{equation}
mit explizitem $C \in \R$ benötigt. Heath-Brown beweist in \cite[S.312-313]{Hea92} eine Abschätzung der Form (\ref{rechenkapitel-zentrale-ungleichung}) für ein konkretes $f$ und $k_1=1, \; k_2=0, \; k_3=2, \; s_{11}=0, \; s_{12}= (\frac{7}{6} + 2\varepsilon)^{-1}$. Wegen $k_2=0$ ist es dabei nicht nötig einen speziellen Wert für $s_2$ festzulegen, man kann z.B. $s_2=0=:s_{21}=:s_{22}$ wählen. Für allgemeine Parameter $k_i$ und $s_{ij}$ gehen wir völlig analog vor, wobei wir an manchen Stellen etwas schärfere Abschätzungen wählen.\\

Zunächst gilt
\begin{eqnarray}
  Re  \{F(z)\}&=&\int_0^{2\gamma} f(x) e^{-Re \{z\} x} \cos(-Im\{z\} x) \,dx \nonumber \\
  &=&\int_0^{2\gamma} f(x) e^{-Re \{z\} x} \cos(Im\{z\} x) \,dx  \nonumber \\
  &=& Re  \{F(\overline{z})\}. \label{Im-z-Vorzeichen-in-Re-z-unwichtig}
\end{eqnarray}
Damit ist $A(s_1,s_2,t)=A(s_1,s_2,-t)$ und wir können im folgenden $t \geq 0$ annehmen. Weiterhin gilt
\begin{equation}\label{rechenkapitel-A-grosser-0}A_{sup} \geq 0,\end{equation}
da für beschränktes $\sigma$ aus (\ref{partielle-integration-von-Re-F}) folgt, dass $$\lim_{t \rightarrow \infty} Re \{F(\sigma+it)\}=0.$$
Diese Konvergenz wird in den Anwendungen später recht schnell ablaufen. Aus diesem Grund macht es Sinn einerseits die Funktion $A(s_1,s_2,t)$ für große $t$, sagen wir $t\geq x_1$, mit groben Mitteln abzuschätzen. Für kleine $t \in [0,x_1]$ geht man dann präzise vor, indem man $A(s_1,s_2,t)$ konkret an sehr vielen Gitterpunkten der Menge $[s_{11},s_{12}] \times [s_{21},s_{22}] \times [0,x_1]$ berechnet. Es ergibt sich dann eine Abschätzung, wenn man das Maximum der Gitterpunkte nimmt und dazu einen Fehler addiert, der aus der Größe der partiellen Ableitungen von $A$ resultiert. Wir beginnen mit der

\subsection{Abschätzung für große $t \; (\geq x_1)$}
Da $f$ Bedingung 2 erfüllt und $k_3 \geq 0$, gilt
\begin{eqnarray*}
  & &Re  \Big{\{} k_1 F(-s_1+it) - k_2 F(-(s_1-s_2)+it) - k_3 F(it) \Big{\}} \\
  & &\;\;\;\;\;\;\;\;\;\;\;\; \leq Re  \big{\{} k_1 F(-s_1+it)-k_2 F(-s_3+it) \big{\}} \\
  & &\;\;\;\;\;\;\;\;\;\;\;\; =: \tilde{A} (s_1,s_3,t) .
\end{eqnarray*}
Wir werden jetzt $\tilde{A}(s_1,s_3,t)$ für $t \geq x_1 >0$ nach oben abschätzen. Zuerst erinnern wir daran, dass $F(z)$ explizit gegeben ist durch (\ref{explizites-F}), also
\begin{equation}\label{explizites-F-nochmal}
  F(z)= \frac{16 \gamma^5}{15} z^{-1}-\frac{8 \gamma^3}{3}z^{-3} + 4\gamma^2(1+ e^{-2 \gamma z})z^{-4} + 4z^{-6} \left(-1+e^{-2 \gamma z} +2 \gamma z e^{-2 \gamma z}\right).
\end{equation}
Setze  $$\tilde{A} (s_1,s_3,t)=\tilde{A_1}(s_1,s_3,t)+\tilde{A_2}(s_1,s_3,t)+\tilde{A_3}(s_1,s_3,t)+\tilde{A_4}(s_1,s_3,t)$$ gemäß den vier Termen in (\ref{explizites-F-nochmal}). Also ist beispielsweise
$$\tilde{A_3}(s_1,s_3,t)= Re \Big{\{} k_1 4\gamma^2(1+ e^{-2 \gamma (-s_1+it)})(-s_1+it)^{-4} -k_2 4\gamma^2(1+ e^{-2 \gamma (-s_3+it)})(-s_3+it)^{-4} \Big{\}}. $$
Um die einzelnen $\tilde{A_i}(s_1,s_3,t) \;\; (i \in \{1,2,3,4\})$ leichter abschätzen zu können, setzen wir voraus, dass \begin{equation}\label{rechenbedingung-1}t \geq x_1 \geq 4 .\end{equation}

\noindent \textbf{Abschätzung von $\tilde{A_1}$.} Wegen $\max \{s_{12},s_{32}\} \leq 4 \leq t$ ist  $s_i/(s_i^2+t^2)$ monoton wachsend in $s_i \; (i \in \{1,3\})$. Aus $$Re  \left\{(-s_i+it)^{-1}\right\} =\frac{-s_i}{s_i^2+t^2}$$ folgt dann
\begin{eqnarray}
   \tilde{A_1}(s_1,s_3,t) &\leq& \frac{16 \gamma^5}{15} \left(\frac{s_{32} k_2}{s_{32}^2+t^2}-\frac{s_{11}k_1}{s_{11}^2+t^2} \right)  \nonumber\\
   &\leq&\frac{16 \gamma^5}{15} \cdot \frac{ t^2 \max \{0, \; s_{32}k_2 - s_{11}k_1\} + s_{11} s_{32} \max \{0, \; s_{11}k_2 - s_{32}k_1 \}}{(s_{32}^2+t^2)(s_{11}^2+t^2)}. \nonumber \\
   & & \label{rechenkapitel-A1}
\end{eqnarray}

\noindent \textbf{Abschätzung von $\tilde{A_2}$.} Wegen $t \geq \max \{s_{12},s_{32}\}$ gilt $3t^2-s_i^2 \geq 0 \; (i \in \{1,3\})$. Zusammen mit $$Re  \left\{(-s_i+it)^{-3}\right\}= \frac{s_i(3t^2-s_i^2)}{(s_i^2+t^2)^3}$$
folgt
\begin{eqnarray}
\tilde{A_2}(s_1,s_3,t) &\leq& \frac{8 \gamma^3}{3} \frac{k_2 s_{32} (3t^2-s_{31}^2)}{(s_{31}^2+t^2)^3} \leq \frac{8 \gamma^3 k_2 s_{32} t^2}{(s_{31}^2+t^2)^3}. \label{rechenkapitel-A2}\end{eqnarray}

\noindent \textbf{Abschätzung von $\tilde{A_3}$.} Es ist $(1+e^{2\gamma s})/((s^2+t^2)^2)$ monoton wachsend in $s$. Das folgt aus $\gamma \geq \frac{1}{2},\; t \geq 4$ und der Standardmethode, dass die Ableitung dann nicht-negativ ist. Durch triviales Abschätzen folgt dann ($|z^{-4}|=|z|^{-4} = (Re\{z\}^2+Im\{z\}^2)^{-2}$)
\begin{equation}\label{rechenkapitel-A3} |\tilde{A_3}(s_1,s_3,t)| \leq 4 \gamma^2 k_1 \frac{(1+e^{2\gamma s_{12}})}{(s_{12}^2+t^2)^2}+ 4 \gamma^2 k_2 \frac{(1+e^{2\gamma s_{32}})}{(s_{32}^2+t^2)^2}. \end{equation}
Damit verbleibt noch die \\

\noindent \textbf{Abschätzung von $\tilde{A_4}$.} Auf Grund der hohen negativen Potenz $z^{-6}$ genügt es diesen Term gänzlich trivial abzuschätzen:
\begin{equation}\label{rechenkapitel-A4} |\tilde{A_4}(s_1,s_3,t)| \leq 4 k_1 \frac{1+e^{2\gamma s_{12}}+2\gamma \sqrt{s_{12}^2+t^2} e^{2\gamma s_{12}}}{t^6}+4 k_2 \frac{1+e^{2\gamma s_{32}}+2\gamma \sqrt{s_{32}^2+t^2} e^{2\gamma s_{32}}}{t^6} .\end{equation}

\noindent Der Ausdruck $t/(s_{ij}^2+t^2)$ ist monoton fallend in $t$, da $t \geq \max \{s_{12},s_{32}\}$. Auch sind Produkte und Summen nicht-negativer monoton fallender Funktionen wieder monoton fallend. Damit sind die rechten Seiten von (\ref{rechenkapitel-A1})-(\ref{rechenkapitel-A4}) alle monoton fallend in $t$. Es folgt
\begin{lemma}\label{rechenkapitel-lemma-B1}
Sei $x_1 \geq 4$. Dann gilt unter den Bedingungen am Anfang von §\ref{para-rechenkapitel} für alle $t \geq x_1$
$$A(s_1,s_2,t) \leq A_1(x_1)+A_2(x_1)+A_3(x_1)+A_3(x_1),$$
wobei $A_1(t)$, $A_2(t)$, $A_3(t)$ bzw. $A_4(t)$ durch die rechte Seite von (\ref{rechenkapitel-A1}), (\ref{rechenkapitel-A2}), (\ref{rechenkapitel-A3}) bzw. (\ref{rechenkapitel-A4}) definiert ist.
\end{lemma}
\vspace*{0.25cm}

\subsection{Abschätzung für kleine $t \; (\in [0,x_1])$}
Wir weisen darauf hin, dass wir für die Abschätzung in diesem Abschnitt der Einfachheit halber nicht die zusätzliche Voraussetzung $s_2 \leq s_1$ verwenden.
Seien $\Delta s_1,\, \Delta s_2,\, \Delta t$ und $x_1$ beliebige positive Konstanten. Definiere ein Gitter $$G \subseteq M:=[s_{11}, s_{12}] \times [s_{21}, s_{22}] \times [0,x_1]$$ durch
\begin{eqnarray}
  G:= \{(s_1, s_2 ,t) &|& s_1=\min\{s_{11} + j_1 \Delta s_1, s_{12}\}, \; j_1=0,\ldots,\left[\frac{s_{12}-s_{11}}{\Delta s_1} \right]+1, \nonumber \\
  & & s_2=\min\{s_{21} + j_2 \Delta s_2,s_{22}\}, \; j_2=0,\ldots,\left[\frac{s_{22}-s_{21}}{\Delta s_2} \right]+1, \nonumber \\
  & & t=\min \{j_3 \Delta t, x_1\},\; j_3=0,\ldots,\left[\frac{x_1}{\Delta t}\right]+1 \} \label{rechenkapitel-Gitter}.
\end{eqnarray}
Im Falle von $s_{i1}=s_{i2}$ für ein $i \in \{1,2\}$ lassen wir auch $\Delta s_i=0$ zu, wobei wir dann in der Definition des Gitters den Term $[(s_{i2}-s_{i1})/\Delta s_i]+1$ durch $0$ ersetzen. Wir schreiben weiterhin
\begin{equation} \label{rechenkapitel-M0-B1}M_0=\max_{(s_1,s_2,t) \in G} A(s_1,s_2,t). \end{equation}
Die Funktion $A(s_1,s_2,t)$ ist in $M$ nach allen drei Variablen differenzierbar. Sei $D_1$ bzw. $D_2$ bzw. $D_3$ eine obere Schranke für den Betrag der partiellen Ableitung nach $s_1$ bzw. $s_2$ bzw. $t$ von $A(s_1,s_2,t)$ in $M$ (genaueres später). Sei nun $(s_1,s_2,t) \in M$ beliebig gewählt und fest. Das Gitter ist gerade so konstruiert, dass es dann ein $(a,b,c) \in G$ gibt, so dass
$$|s_1-a|\leq \frac{\Delta s_1}{2}, \;|s_2-b|\leq \frac{\Delta s_2}{2}, \;|t-c|\leq \frac{\Delta t}{2}.$$
 Zusammen mit dem Mittelwertsatz der Differentialrechnung folgt
\begin{eqnarray*}
  |A(s_1,s_2,t)-A(a,b,c)| &\leq& |A(s_1,s_2,t)-A(a,b,t)| + |A(a,b,t)-A(a,b,c)| \\
  &\leq& |A(s_1,s_2,t)-A(a,s_2,t)| +|A(a,s_2,t)-A(a,b,t)|\\
  & & + \; |A(a,b,t)-A(a,b,c)| \\
  &\leq& \frac{\Delta s_1}{2} D_1 + \frac{\Delta s_2}{2} D_2  + \frac{\Delta t }{2} D_3.
\end{eqnarray*}
Also folgt
$$  \sup_{(s_1,s_2,t) \in M} A(s_1,s_2,t) \leq M_0 + \frac{\Delta s_1}{2} D_1 + \frac{\Delta s_2}{2} D_2  + \frac{\Delta t }{2} D_3 .$$
Es verbleibt die Abschätzung der $D_i$. Aus der Definition der Laplace-Transformierten $F$ von $f$ folgt sofort
\begin{equation*}
A(s_1,s_2,t)= \int_{0}^{2 \gamma} f(x) \cos(tx) e^{s_1 x}(k_1  - k_2 e^{-s_2 x}) \,dx - k_3 \int_0^{2 \gamma} f(x) \cos(tx) \,dx.
\end{equation*}
Da die Funktionen im Inneren der Integrale (abhängig von den beiden Variablen $t,x$ bzw. $s_1,x$ bzw. $s_2,x$) stetig sind, können wir unter dem Integral differenzieren (Standardresultat aus der Analysis, siehe z.B. \cite[S.235]{Lan83}) und erhalten für $(s_1,s_2,t) \in M$
\begin{equation}\label{rechenkapitel-ableitung-s1}
\left|\frac{d A(s_1,s_2,t)}{d s_1}\right| \leq d_0 \int_{0}^{2 \gamma} x f(x) e^{s_{12} x} \,dx =:D_1,
\end{equation}
\begin{equation}\label{rechenkapitel-ableitung-s2}
\left|\frac{d A(s_1,s_2,t)}{d s_2}\right| \leq k_2 \int_{0}^{2 \gamma} x f(x) e^{s_{32} x} \,dx =:D_2, \end{equation}
und
\begin{equation}\label{rechenkapitel-ableitung-t}
\left|\frac{d A(s_1,s_2,t)}{d t}\right| \leq d_0 \int_{0}^{2 \gamma} x f(x) e^{s_{12} x} \,dx + k_3 \int_0^{2\gamma}x f(x) \,dx =:D_3,
\end{equation}
wobei $d_0=\sup\limits_{x \in [0,2\gamma]} |k_1  - k_2 e^{-s_2 x}|$. Da $h(x)=k_1  - k_2 e^{-s_2 x}$ monoton in $x$ ist, nimmt $|h(x)|$ sein Supremum in $[0,2\gamma]$ an einem der beiden Randpunkte $x=0$ oder $x=2\gamma$ an. Es folgt
\begin{eqnarray}
d_0&=& \max \{|k_1  - k_2|,|k_1  - k_2 e^{-s_{22} 2 \gamma}| \} \nonumber  \\
&=&\max \{k_2-k_1,k_1  - k_2 e^{-s_{22} 2 \gamma} \} \label{rechenkapitel-d0} .\end{eqnarray}
Die letzte Gleichung folgt aus der Überlegung, dass, wenn $k_2-k_1<0$, so ist $d_0=k_1  - k_2 e^{-s_{22} 2 \gamma}$ und analog, wenn $k_1  - k_2 e^{-s_{22} 2 \gamma}<0$, so ist $d_0=k_2-k_1$. Wir fassen alles nochmal im folgenden Lemma zusammen.

\begin{lemma}\label{rechenkapitel-lemma-B2}
  Seien die Voraussetzungen vom Beginn von §\ref{para-rechenkapitel} gegeben. Ferner seien $\Delta s_1, \; \Delta s_2, \; \Delta t$ und $x_1$ positive Konstanten. Im Falle von $s_{i1}=s_{i2}$ für ein $i \in \{1,2\}$ ist auch $\Delta s_i=0$ zulässig. Mit den Notationen aus (\ref{rechenkapitel-Gitter})-(\ref{rechenkapitel-d0}) folgt dann für $t \in [0,x_1]$, dass
  $$ A(s_1,s_2,t) \leq M_0 + \frac{\Delta s_1}{2} D_1 + \frac{\Delta s_2}{2} D_2  + \frac{\Delta t }{2} D_3 .$$
  Zusammen mit Lemma \ref{rechenkapitel-lemma-B1} folgt also
  $$A_{sup} \leq \max \left\{A_1(x_1)+A_2(x_1)+A_3(x_1)+A_4(x_1), \, M_0 + \frac{\Delta s_1}{2} D_1 + \frac{\Delta s_2}{2} D_2  + \frac{\Delta t }{2} D_3 \right\} .$$
\end{lemma}

Die Abschätzung in Lemma \ref{rechenkapitel-lemma-B1} ist umso "`schärfer"', je größer man das $x_1$ wählt und diejenige in Lemma \ref{rechenkapitel-lemma-B2} je kleiner man die Parameter $\Delta s_1, \; \Delta s_2$ und $\Delta t$ wählt. In diesem Zusammenhang sei darauf hingewiesen, dass man mit der Methode in diesem Abschnitt den Term $A_{sup}$ "`beliebig scharf"' nach oben abschätzen kann, wenn man nur $x_1$ groß genug und $\Delta s_1,\; \Delta s_2, \; \Delta t$ klein genug wählt, bzw. den Computer nur lange genug rechnen lässt.\\

\noindent Für die konkrete Berechnung bliebe noch die Wahl der Parameter $x_1$, $\Delta s_1$, $\Delta s_2$ und $\Delta t$, so dass $A_{sup}$ bei geringem Rechenaufwand möglichst gut abgeschätzt wird. Was man dabei u.a. beachten kann ist, die Parameter $\Delta s_1$, $\Delta s_2$ und $\Delta t$ so zu wählen, dass $$D_1 \Delta s_1 \approx D_2 \Delta s_2 \approx D_3 \Delta t .$$
Durch eine solche Wahl wird nämlich die Anzahl der Rechenschritte (= Anzahl der Gitterpunkte  $\approx C_1(\Delta s_1 \Delta s_2 \Delta t)^{-1}$) bei fast gleichbleibender Abschätzung minimiert: Man kann das in die Analysisaufgabe übersetzen, die Funktion $h(x,y,z)=C_1(xyz)^{-1}$ zu minimieren unter der Nebenbedingung $C_2=x+y+z$.\\

\section{$\lambda'$-Abschätzungen}\label{para-ldash-abschaetzungen}
Für diesen Abschnitt setzen wir voraus, dass mindestens eines der beiden Elemente $\chi_1$, $\rho_1$ komplex ist (und damit insbesondere auch existiert). Auch setzen wir voraus, dass $\rho'$ existiert. Ab jetzt werden wir darauf verzichten zu erwähnen, dass wir voraussetzen, dass eine gewisse Nullstelle überhaupt existiert. Sobald wir nämlich diese Nullstelle aufschreiben und damit arbeiten, versteht sich diese implizite Annahme von selbst. Weiterhin setzen wir voraus, dass eine Funktion $f$ gegeben ist, die Bedingung 1 und 2 erfüllt. Am Ende wird dies die Funktion in (\ref{explizites-f}) sein.\\
In diesem Abschnitt verbessern wir \cite[Table 8 (§9)]{Hea92}. Dieses "`Table 8"' liefert Aussagen der Form
$$\lambda_1 \leq C_1 \Rightarrow \lambda' \geq C_2. $$
 Mit Hilfe von Verbesserungspotential Nr. 2 kann man die Konstante $C_2$ in diesen Abschätzungen klar verbessern. Weiter unten werden wir an passender Stelle erklären, worin dieses Verbesserungspotential überhaupt liegt.\\
Wir sollten noch erwähnen, dass für den Beweis der letztendlich zulässigen Linnikschen Konstante $L$ das Einbringen der hier verbesserten $\lambda'$-Abschätzungen vergleichsweise wenig Beitrag liefert. Der eigentliche Vorteil liegt darin, dass wir diese in späteren Abschnitten verwenden können zur Herleitung der wichtigen $\lambda_2$- und $\lambda_3$-Abschätzungen.

\subsection{Anfangsungleichung}

Wir beginnen mit der Ungleichung in \cite[S.302 oben]{Hea92}, nämlich
\begin{eqnarray}
  0 &\leq& \chi_0(n) \left(1+ Re  \left\{\frac{\chi_1(n)}{n^{i\gamma'}} \right\} \right) \left(k+ Re  \left\{\frac{\chi_1(n)}{n^{i\gamma_1}} \right\} \right)^2 \nonumber \\
  &=&(k^2+\frac{1}{2}) \left(\chi_0(n)+ Re  \left\{\frac{\chi_1(n)}{n^{i\gamma'}} \right\}\right) \nonumber \\
  & & + k \left(Re  \left\{\frac{\chi_0(n)}{n^{i(\gamma_1-\gamma')}}\right\} +2Re  \left\{\frac{\chi_1(n)}{n^{i\gamma_1}}\right\} + Re  \left\{\frac{\chi_1(n)^2}{n^{i(\gamma_1+\gamma')}}\right\} \right) \nonumber \\
  & &+\frac{1}{4}\left(Re  \left\{\frac{\chi_1(n)}{n^{i(2\gamma_1-\gamma')}}\right\} +2Re  \left\{\frac{\chi_1(n)^2}{n^{2i\gamma_1}}\right\} + Re  \left\{\frac{\chi_1(n)^3}{n^{i(2\gamma_1+\gamma')}}\right\} \right), \label{para9-l-dash-startinequality}
\end{eqnarray}
wobei $k \geq 0$ eine zu wählende Konstante ist.

\begin{bemerkung}
Wir haben hier und an späterer Stelle leichte Variationen der von Heath-Brown benutzten Ungleichungen getestet. Beispielsweise probierten wir statt der obigen Ungleichung die Ungleichung $$ 0 \leq \chi_0(n) \left(l+Re  \left\{\frac{\chi_1(n)}{n^{i \gamma'}}\right\}\right)^2\left(k+Re  \left\{\frac{\chi_1(n)}{n^{i \gamma_1}}\right\}\right)^2=\ldots \; , $$
welche von den beiden Parametern $l,\,k \geq 0$ abhängt. Doch diese variierten Ungleichungen lieferten im Allgemeinen keine oder nur sehr minimal bessere Werte.
\end{bemerkung}
\vspace*{0.25cm}

Während Heath-Brown jetzt mit der in §\ref{para-methoden-fuer-nullstellenfreie-regionen} kurz erwähnten "`ersten Methode"' (Verwendung von $(L'/L)^{(k)} (s,\chi)$) weitermacht, verwenden wir die "`zweite Methode"' mit den Lemmata \ref{hb-lemma-5-3} und \ref{hb-lemma-5-2-extended}, welche in §\ref{hb-para-6} angesprochen wurde. Das Problem bei letzterer Methode ist, dass die resultierenden Ungleichungen (vom Typ (\ref{beispiel-l1-abschaetzung})) auf der rechten Seite Terme vom Typ $$\sup_{t \in \R} Re  \big{\{} k_1 F(-s_1+it)-k_2 F(-s_3+it)-k_3 F(it) \big{\}}$$ beinhalten können. Verbesserungspotential Nr. 2 (\cite[S.332]{Hea92}) liegt nun darin, sich von solchen Termen nicht abschrecken zu lassen, sondern sie konsequent und genau abzuschätzen, und zwar so wie es in \cite[§10]{Hea92} gemacht wird bzw. wie wir es im letzten Abschnitt nochmal allgemein vorgeführt haben.\\

Sei jetzt $\lambda^{\star}>0$ mit
\begin{equation}\label{para9-l-dash-bedingung1}\lambda^{\star} \leq \min \{\lambda_2, \lambda' \}.\end{equation} Im Falle, dass $\rho_2$ nicht existiert, setzen wir dabei nur $\lambda^{\star} \leq \lambda'$ voraus. Setze $$\beta^{\star}=1-\lambda^{\star} \mathscr{L}^{-1}$$ und erhalte aus (\ref{para9-l-dash-startinequality}) mit der Standardmethode (also Multiplikation mit $ n^{-\beta^{\star}} \Lambda(n) f(\mathscr{L}^{-1}\log n)$, dann Summation über $n$) die Ungleichung
\begin{eqnarray}
  0 &\leq& (k^2+\frac{1}{2})K(\beta^{\star},\chi_0)+(k^2+\frac{1}{2}) K(\beta^{\star}+i \gamma',\chi_1) \label{para9-l-dash-inequality} \\
  & &+k K(\beta^{\star} + i(\gamma_1-\gamma'),\chi_0) +2k K(\beta^{\star}+i \gamma_1,\chi_1) +k K(\beta^{\star} + i(\gamma_1+\gamma'),\chi_1^2) \nonumber \\
  & &+ \frac{1}{4} K(\beta^{\star} + i(2\gamma_1-\gamma'),\chi_1) +\frac{1}{2} K(\beta^{\star}+2 i \gamma_1,\chi_1^2) +\frac{1}{4} K(\beta^{\star} + i(2\gamma_1+\gamma'),\chi_1^3), \nonumber
\end{eqnarray}
wobei
$$K(s,\chi)=\sum_{n=1}^{\infty} \Lambda(n) Re  \left\{\frac{\chi(n)}{n^{s}}\right\} f(\mathscr{L}^{-1} \log n) .$$
Wir unterscheiden zwei Fälle. Zuerst zu

\subsection{Fall 1: $ord \; \chi_1 \geq 5$}
Wir beweisen
\begin{lemma}\label{lemma-lstrich-ordchi-geq-5}
  Sei $ord \; \chi_1 \geq 5$ und $f$ eine Funktion, die Bedingung 1 und 2 erfüllt. Weiterhin sei $\lambda^{\star}$ eine Konstante mit  $0<\lambda^{\star} \leq \min \{\lambda',\lambda_2 \}$. Wenn $\rho_2$ nicht existiert, so setzen wir nur $0<\lambda^{\star} \leq \lambda'$ voraus. Seien ferner $k$ und $\varepsilon$ positive Konstanten. Dann gilt für $q \geq q_0(f,k,\varepsilon)$
\begin{eqnarray}
  0 &\leq& (k^2+\frac{1}{2})(F(-\lambda^{\star})-F(\lambda' - \lambda^{\star}))-2k  F(\lambda_1-\lambda^{\star}) + \frac{f(0)}{6} (k^2+3k+\frac{3}{2}) + \varepsilon \nonumber \\
  & &+ \sup_{t \in \R} Re  \left\{k F(-\lambda^{\star}+i t)- (k^2+\frac{3}{4})F(\lambda_1-\lambda^{\star} + i t) \right\} . \label{para9-l-dash-hauptfall}
\end{eqnarray}
\end{lemma}

\noindent Für den Beweis beginnen wir mit der Ungleichung (\ref{para9-l-dash-inequality}). In dieser schätzen wir die acht verschiedenen "`$K$-Terme"' mit Hilfe der Lemmata \ref{hb-lemma-5-3} und \ref{hb-lemma-5-2-extended} ab. Zuerst bekommen wir mit Lemma \ref{hb-lemma-5-3} die Ungleichungen
\begin{eqnarray}
(k^2+\frac{1}{2})\mathscr{L}^{-1} K(\beta^{\star},\chi_0)&\leq&  (k^2+\frac{1}{2})( F(-\lambda^{\star}) +\varepsilon), \label{para9-l-dash-help0} \\
k\mathscr{L}^{-1} K(\beta^{\star} + i(\gamma_1-\gamma'),\chi_0)&\leq&  k (Re \{ F(-\lambda^{\star}+i(\mu_1-\mu'))\} + \varepsilon) .\label{para9-l-dash-help1}
\end{eqnarray}
Wir sparen uns im Allgemeinen den Hinweis, dass die jeweiligen Ungleichungen nur für $q \geq q_0$ gelten. Dabei ist es wichtig zu bemerken, dass wir nur endlich viele solche Schlüsse in dieser Diplomarbeit ziehen werden, d.h. wir bekommen nur endliche viele $q_0$. Da es endlich viele sind, kann man das Maximum $q_{max}$ all dieser $q_0$ nehmen, so dass dann sämtliche Aussagen der Arbeit für $q \geq q_{max}$ gelten.

Wegen der Voraussetzung $ord \; \chi_1 \geq 5$ gilt $\chi_1^2 \neq \chi_0, \chi_1, \overline{\chi_1}.$ Letztere Schreibweise verwenden wir als Synonym für $\chi_1^2 \notin \{ \chi_0, \chi_1, \overline{\chi_1} \}$. \\
Da $\chi_1^2 \neq \chi_0$ können wir zur Abschätzung des Terms
$K(\beta^{\star}+i(\gamma_1+\gamma'),\chi_1^2)$ das Lemma \ref{hb-lemma-5-2-extended} verwenden. Betrachte für diesen Fall, also $s=\beta^{\star}+i(\gamma_1+\gamma') \in R(2 L_1) \subseteq R(9 L_1)$ und $\chi=\chi_1^2$, die beiden Mengen $A_1$ und $A_2$ aus Lemma \ref{hb-lemma-5-2-extended}. Mit $\chi_1^2 \neq \chi_1, \overline{\chi_1}$ und $\beta^{\star} \geq Re\{\rho_2\}$ folgt aus (\ref{Eigenschaft1-Nullstellen}), dass $A_1=\emptyset$. Falls $\rho_2$ nicht existiert kann man genauso schließen. Da wir außerdem keine weiteren Informationen zu dem Charakter $\chi_1^2$ bzw. den Nullstellen von $L(s,\chi_1^2)$ in $R$ haben, wählen wir $A_2=\emptyset$. Das Lemma liefert also mit $\phi(\chi_1^2) \leq \frac{1}{3}$
\begin{eqnarray}
k \mathscr{L}^{-1}K(\beta^{\star} + i(\gamma_1+\gamma'),\chi_1^2) &\leq& k (\frac{f(0)}{6}+\varepsilon). \label{para9-l-dash-help5}
\end{eqnarray}
Völlig analog folgen
\begin{eqnarray}
  \frac{1}{2}\mathscr{L}^{-1} K(\beta^{\star} + 2i\gamma_1,\chi_1^2) &\leq& \frac{1}{2}(\frac{f(0)}{6}+\varepsilon), \label{para9-l-dash-help6}\\
  \frac{1}{4}\mathscr{L}^{-1} K(\beta^{\star} + i(2\gamma_1+\gamma'),\chi_1^3) &\leq& \frac{1}{4}(\frac{f(0)}{6}+\varepsilon) . \label{para9-l-dash-help7}
\end{eqnarray}
Zur Abschätzung von $K(\beta^{\star}+i \gamma',\chi_1)$ benutzen wir wieder Lemma \ref{hb-lemma-5-2-extended}. Diesmal hat jedoch die Funktion $L(s,\chi_1)$ die Nullstelle $\rho_1$, für die durchaus gelten kann, dass $\beta^{\star}<Re\{\rho_1\}$. Ist dies der Fall, so folgt unter Berücksichtigung von $\beta^{\star} \geq Re\{\rho'\}$, dass $A_1=\{\rho_1 \}$. Zusätzlich wählen wir dann $A_2=\{\rho'\}$. Ist jedoch $Re\{\rho_1\} \leq \beta^{\star} $, so ist $A_1=\emptyset$ und wir wählen $A_2=\{\rho_1,\rho' \}$. In beiden Fällen ist also $A_1 \cup A_2 = \{\rho_1, \rho' \}$ und wir bekommen
\begin{eqnarray}
  (k^2+\frac{1}{2}) \mathscr{L}^{-1}K(\beta^{\star}+i \gamma',\chi_1) &\leq& (k^2+\frac{1}{2})(-Re  \{F(\lambda_1-\lambda^{\star} + i (\mu'-\mu_1))\}  \nonumber \\
  & &-F(\lambda'-\lambda^{\star}) + \frac{f(0)}{6} + \varepsilon ). \label{para9-l-dash-help2}
\end{eqnarray}
Für die Abschätzung der beiden verbleibenden Terme $K(\beta^{\star}+i \gamma_1,\chi_1)$ und $K(\beta^{\star} + i(2\gamma_1-\gamma'),\chi_1)$ benutzen wir die Nullstelle $\rho'$ nicht, wir wählen also $A_2$ in Lemma \ref{hb-lemma-5-2-extended} so, dass $A_1 \cup A_2 = \{\rho_1 \}$ und es folgt
\begin{eqnarray}
2k\mathscr{L}^{-1} K(\beta^{\star}+i \gamma_1,\chi_1) &\leq&  2k (-F(\lambda_1-\lambda^{\star}) +\frac{f(0)}{6} + \varepsilon), \label{para9-l-dash-help3} \\
\frac{1}{4} \mathscr{L}^{-1}K(\beta^{\star} + i(2\gamma_1-\gamma'),\chi_1) &\leq& \frac{1}{4}(-Re  \{F(\lambda_1-\lambda^{\star} + i(\mu_1-\mu'))\} +\frac{f(0)}{6}+ \varepsilon). \nonumber \\
& & \label{para9-l-dash-help4}
\end{eqnarray}
Schließlich folgt noch aus (\ref{Im-z-Vorzeichen-in-Re-z-unwichtig}), dass
$$Re  \{F(\lambda_1-\lambda^{\star} + i (\mu'-\mu_1))\}= Re  \{F(\lambda_1-\lambda^{\star} + i (\mu_1-\mu'))\} .$$
Dies verwenden wir in (\ref{para9-l-dash-help2}). Benutzt man nun in (\ref{para9-l-dash-inequality}) die Abschätzungen (\ref{para9-l-dash-help0})-(\ref{para9-l-dash-help4}), so folgt das Lemma (das $t$ im Lemma entspricht dabei dem $\mu_1-\mu'$; außerdem müssen wir natürlich das $\varepsilon$ vorher passend verkleinern).\\

Wir weisen noch daraufhin, dass wir an manchen Stellen auch weitere Nullstellen (nämlich $\rho'$) hätten hinzunehmen können. Dies hätte jedoch die nachfolgenden Rechnungen bei nur minimaler Verbesserung zu sehr strapaziert.

Das nächste Mal, wenn wir mit Hilfe von Lemma \ref{hb-lemma-5-2-extended} einen "`$K$-Term"' abschätzen, werden wir nur angeben, welche Nullstellen $\rho$ wir berücksichtigen (müssen). Wir werden also die Menge $A_1 \cup A_2$ angeben. Die leichten Details zur Zulässigkeit der entsprechenden Wahl von $A_1 \cup A_2$ lassen wir dann in der Regel weg.

\subsection{Fall 2: $ord \; \chi_1 \in \{2,3,4\}$}
In diesem Unterabschnitt beweisen wir
\begin{lemma}\label{lemma-lstrich-ordchi-234}
  Sei $ord \; \chi_1 \leq 4$ und im Falle, dass $\chi_1$ reell ist, setzen wir $\rho_1$ komplex voraus. Sei $f$ eine Funktion, die Bedingung 1 und 2 erfüllt und seien $k$ und $\varepsilon$ positive Konstanten. Dann gilt für $q \geq q_0(f,k,\varepsilon)$, dass
\begin{eqnarray}
  0 &\leq& (k^2+\frac{1}{2})(F(-\lambda_1)-F(\lambda'-\lambda_1))-2kF(0) + \frac{f(0)}{8}(k^2+3k+\frac{3}{2}) +\varepsilon \label{para9-l-dash-ungleichung-nebenfall} \\
  & &+ 2\sup_{t \in \R} Re  \left\{kF(-\lambda_1+it) - (k^2+\frac{3}{4})F(it)\right\} + \sup_{t\in \R} Re  \left\{\frac{1}{2}F(-\lambda_1+it)-2kF(it)\right\}. \nonumber
\end{eqnarray}
\end{lemma}

\noindent Zum Beweis des Lemmas verändern wir erstmal die Ungleichung (\ref{para9-l-dash-inequality}), indem wir dort $\beta^{\star}$ durch $\beta_1$ ersetzen. Dies ist offensichtlich möglich. Die daraus resultierende Ungleichung nennen wir (\ref{para9-l-dash-inequality})'.  Wir möchten nun analog zum Fall $ord \; \chi_1 \geq 5$ die acht verschiedenen "`$K$-Terme"' in (\ref{para9-l-dash-inequality})' abschätzen. Dazu notieren wir vorher ein paar generelle Bemerkungen.\\

Im Fall $ord \; \chi_1=2$ muss man beachten, dass aus $L(\rho_1,\chi_1)=0$ folgt, dass $L(\overline{\rho_1},\chi_1)=0$. Dabei gilt $\overline{\rho_1} \in R$ und $\overline{\rho_1} \neq \rho_1, \rho'$, denn nach Voraussetzung ist dann $\rho_1$ komplex und $\overline{\rho_1} \neq \rho'$ nach Definition. Diese zusätzliche Nullstelle $\overline{\rho_1}$ werden wir im Fall $ord \; \chi_1=2$ manchmal verwenden.

Analog ist immer $\overline{\rho_1} \in R$ eine Nullstelle von $L(s,\overline{\chi_1})$. Diese Nullstelle werden wir einmal benutzen und zwar im Fall $ord \; \chi_1 =3$.

Man beachte weiterhin, dass für sämtliche Nullstellen $\rho \in R$ der Funktionen $L(s,\chi) \; (\chi \neq \chi_0)$ gilt, dass $Re \{\rho\} \leq Re \{\rho_1\}$. Wegen der Wahl von $\beta_1$ anstelle von $\beta^{\star}$ haben wir also freie Wahl welche Nullstellen bzw. Terme "`$-Re\{F((s-\rho)\mathscr{L})\}$"' wir bei der Abschätzung der Ausdrücke $K(\beta^{\star}+it,\chi) \; (\chi \neq \chi_0$) miteinbeziehen. Mit anderen Worten: Die Menge $A_1$ in Lemma \ref{hb-lemma-5-2-extended} wird immer leer sein.

Eine letzte Bemerkung ist, dass für $ord \; \chi_1 \leq 4$ gilt, dass die Ordnungen von $\chi_1,\; \chi_1^2$ und $\chi_1^3$ auch alle $\leq 4$ sind. Also haben wir immer $\phi=\frac{1}{4}$.\\

\noindent Wir verwenden die Abkürzung \label{Definition-K0}
$$K_0(s,\chi):=\mathscr{L}^{-1} K(s,\chi)-\varepsilon. $$

\noindent Zunächst haben wir mit Lemma \ref{hb-lemma-5-3}
\begin{eqnarray*}
  K_0(\beta_1,\chi_0) &\leq& F(-\lambda_1) \;\; \;\; \text{ für } ord \; \chi_1 \in \{2,3,4\}, \\
  K_0(\beta_1+i(\gamma_1-\gamma'),\chi_0) &\leq& Re\{F(-\lambda_1 + i(\mu_1-\mu'))\} \;\;\;\;   ord \; \chi_1 \in \{2,3,4\} .
\end{eqnarray*}\\

\noindent Aus Lemma \ref{hb-lemma-5-2-extended} folgt
\begin{eqnarray*}
  K_0(\beta_1+i\gamma',\chi_1) &\leq& \left\lbrace \begin{array}{ll}
-F(\lambda'-\lambda_1)-Re\{F(i(\mu'-\mu_1))\} \\\;\;\;\;\;\;\;\;\;\;\;\;\;\;\;\; -Re\{F(i(\mu'+\mu_1))\}+\frac{f(0)}{8} &\; ord \; \chi_1 = 2, \\
&\\
-F(\lambda'-\lambda_1)-Re\{F(i(\mu'-\mu_1))\}+\frac{f(0)}{8} &\; ord \; \chi_1 \in \{3,4\},
 \end{array} \right.
\end{eqnarray*}
wobei für $ord \; \chi_1=2$ wir $A_1 \cup A_2=\{\rho', \rho_1, \overline{\rho_1} \}$ wählten und im Fall $ord \; \chi_1 \in \{3,4\}$ wählten wir $A_1 \cup A_2=\{\rho', \rho_1\}$.\\

\noindent Weiterhin folgt mit $A_1 \cup A_2=\{\rho_1, \overline{\rho_1} \}$ bzw. $\{\rho_1 \}$ für den Fall $ord \; \chi_1=2$ bzw. $ord \; \chi_1 \in \{3,4\}$
\begin{eqnarray*}
K_0(\beta_1+i\gamma_1,\chi_1) &\leq& \left\lbrace \begin{array}{ll}
-F(0)-Re\{F(2i \mu_1)\}+\frac{f(0)}{8}  &\; ord \; \chi_1 = 2, \\
-F(0)+\frac{f(0)}{8}  &\; ord \; \chi_1 \in \{3,4\} .
  \end{array} \right.
\end{eqnarray*}\\

\noindent Es ist $\chi_1^2=\chi_0$ bzw. $\overline{\chi_1}$ für $ord \; \chi_1=2$ bzw. $3$. Verwendet man also Lemma \ref{hb-lemma-5-3} für den Fall $ord \; \chi_1=2$ und Lemma \ref{hb-lemma-5-2-extended} mit $A_1 \cup A_2 =\{\overline{\rho_1}\}$ bzw. $\emptyset$ für $ord \; \chi_1=3$ bzw. $4$, so folgt
\begin{eqnarray*}
  K_0(\beta_1+i(\gamma_1+\gamma'),\chi_1^2) &\leq& \left\lbrace \begin{array}{ll}
Re \{F(-\lambda_1+i(\mu_1+\mu'))\}  &\; ord \; \chi_1 = 2, \\
-Re\{F(i(2\mu_1+\mu'))\} + \frac{f(0)}{8} &\; ord \; \chi_1 = 3, \\
\frac{f(0)}{8}  &\; ord \; \chi_1 = 4 .
 \end{array} \right.
\end{eqnarray*}\\

\noindent Mit $A_1 \cup A_2 =\{\rho_1 \}$ in allen Fällen folgt
\begin{eqnarray*}
   K_0(\beta_1+i(2\gamma_1-\gamma'),\chi_1) &\leq& -Re\{F(i(\mu_1-\mu'))\}+\frac{f(0)}{8}  \;\; \;\;   ord \; \chi_1 \in \{2,3,4\}.
\end{eqnarray*}\\

\noindent Außerdem folgt mit der Wahl $A_1 \cup A_2=\emptyset$ in den Fällen $ord \; \chi_1 \in \{3,4\}$, dass
\begin{eqnarray*}
  K_0(\beta_1+2 i \gamma_1,\chi_1^2) &\leq& \left\lbrace \begin{array}{ll}
Re\{F(-\lambda_1+2i\mu_1)\}   &\; ord \; \chi_1 = 2, \\
\frac{f(0)}{8} &\; ord \; \chi_1 \in \{3,4\}
 \end{array} \right.\\
\end{eqnarray*}
und mit $A_1 \cup A_2=\{\rho_1 \}$ bzw. $\emptyset$ im Fall $ord \; \chi_1=2$ bzw. $4$ folgt
\begin{eqnarray*}
 K_0(\beta_1+i (2\gamma_1+\gamma'),\chi_1^3) &\leq& \left\lbrace \begin{array}{ll}
-Re\{F(i(\mu_1+\mu'))\}+\frac{f(0)}{8}   &\; ord \; \chi_1 = 2, \\
Re\{F(-\lambda_1+i(2\mu_1+\mu'))\}   &\; ord \; \chi_1 = 3, \\
\frac{f(0)}{8} &\; ord \; \chi_1 =4 .
 \end{array} \right.
\end{eqnarray*}\\

\noindent Setzen wir die obigen Abschätzungen in (\ref{para9-l-dash-inequality})' ein und berücksichtigen wir (\ref{Im-z-Vorzeichen-in-Re-z-unwichtig}), so folgt für $ord \; \chi_1=2$
\begin{eqnarray*}
  0 &\leq& (k^2+\frac{1}{2})(F(-\lambda_1)-F(\lambda'-\lambda_1))-2kF(0) + \frac{f(0)}{8}(k^2+2k+1) + \varepsilon \\
  & &+ 2\sup_{t \in \R} Re  \left\{kF(-\lambda_1+it) - (k^2+\frac{3}{4})F(it)\right\} + \sup_{t\in \R} Re  \left\{\frac{1}{2}F(-\lambda_1+it)-2kF(it)\right\},
\end{eqnarray*}
für  $ord \; \chi_1 =3$ bekommen wir
\begin{eqnarray*}
  0 &\leq& (k^2+\frac{1}{2})(F(-\lambda_1)-F(\lambda'-\lambda_1))-2kF(0) + \frac{f(0)}{8}(k^2+3k+\frac{5}{4}) + \varepsilon\\
   & &+ \sup_{t \in \R} Re  \left\{kF(-\lambda_1+it) - (k^2+\frac{3}{4})F(it)\right\} + \sup_{t\in \R} Re  \left\{\frac{1}{4}F(-\lambda_1+it)-k F(it)\right\}
\end{eqnarray*}
und für  $ord \; \chi_1 =4$ schließlich
\begin{eqnarray*}
  0 &\leq& (k^2+\frac{1}{2})(F(-\lambda_1)-F(\lambda'-\lambda_1))-2kF(0) + \frac{f(0)}{8}(k^2+3k+\frac{3}{2}) + \varepsilon\\
   & &+ \sup_{t \in \R} Re  \left\{kF(-\lambda_1+it) - (k^2+\frac{3}{4})F(it)\right\} .
\end{eqnarray*}
Aus diesen drei Ungleichungen folgt nach Berücksichtigung von (\ref{rechenkapitel-A-grosser-0}) die Aussage des Lemmas.

\subsection{Ergebnisse}
Aus der Ungleichung (\ref{para9-l-dash-hauptfall}) bzw. (\ref{para9-l-dash-ungleichung-nebenfall}) gilt es nun Abschätzungen für Fall 1 bzw. Fall 2 zu folgern. Wir besprechen die dazugehörige Vorgehensweise exemplarisch an einem Beispiel.\\

\noindent \textbf{Beispiel: $\lambda_1 \in [0.34,0.36]$ in Fall 1 ($ord \; \chi_1 \geq 5$)}\\
Angenommen wir sind im Fall 1 und $\lambda_1 \in [0.34,0.36]=:[\lambda_{11},\lambda_{12}]$. Nehme nun an, dass $\lambda' \leq 2.06$ (die Konstante $2.06$ ist im Nachhinein bestmöglich). Ziel ist es einen Widerspruch herzuleiten, womit wir dann unter den gemachten Voraussetzungen $\lambda' > 2.06$ bewiesen hätten.

Die Tabelle \cite[Table 8 (§9)]{Hea92} besagt $\lambda' \geq 1.309$. Die Tabelle \cite[Table 10/Lemma 9.4 (§9)]{Hea92} gibt $\lambda_2 \geq 0.903$. Also können wir $\lambda^{\star}=0.903 \leq \min \{\lambda',\lambda_2\}$ wählen.

Sei $A$ das Supremum auf der rechten Seite von (\ref{para9-l-dash-hauptfall}). Mittels §\ref{para-rechenkapitel} bekommen wir die Abschätzung $A \leq C=0.017\ldots$, wobei wir dafür die folgenden Parameter wählten:
\begin{eqnarray}
 & &\gamma=1.13-\frac{\lambda_{12}}{5}, \; \;k=0.75+\frac{\lambda_{12}}{7}, \label{para9-l-dash-vor0}\\
 & & k_1=k,\; k_2=k^2+\frac{3}{4},\; k_3=0, \label{para9-l-dash-vor1}\\
 & & s_1=\lambda^{\star},\; s_{11}=s_{12}=\lambda^{\star}, \;s_2=\lambda_1, \; s_{21}=\lambda_{11},\; s_{22}=\lambda_{12}, \label{para9-l-dash-vor2}\\
 & &\Delta s_1=0,\; \Delta s_2=0.004, \; \Delta t =0.004,\; x_1=15 . \label{para9-l-dash-vor3}
\end{eqnarray}
Die Wahl des $\gamma$ und $k$ erfolgte dabei durch Experimentieren mittels Computer.

Nun ist die rechte Seite von (\ref{para9-l-dash-hauptfall}) ohne das Supremum monoton wachsend in $\lambda_1$ und $\lambda'$, was sofort aus der Definition der Laplace-Transformierten $F(z)$ folgt. Also folgt mit Hilfe der Abschätzung $A \leq C$, dass
\begin{equation} \label{para9-l-dash-endungleichung}- \varepsilon \leq (k^2+\frac{1}{2})(F(-\lambda^{\star})-F(2.06 - \lambda^{\star}))-2k  F(\lambda_{12}-\lambda^{\star}) + \frac{f(0)}{6} (k^2+3k+\frac{3}{2}) + C. \end{equation}
Schließlich berechnen wir die rechte Seite der letzten Ungleichung mit dem Computer und erhalten einen negativen Wert. Für hinreichend kleines $\varepsilon>0$ ist dies ein Widerspruch. Also haben wir für $q \geq q_0(f,k,\varepsilon)$ bewiesen, dass aus $\lambda_1 \in [0.34,0.36]$ die Abschätzung
$$ \lambda' > 2.06$$
folgt.\\

\noindent \textbf{Bemerkungen zur Tabelle für Fall 1:}\\
Völlig analog zum eben genannten Beispiel erstellen wir eine Tabelle für Fall 1 für die Intervalle
$$\lambda_1 \in [\lambda_{11},\lambda_{12}]=[0.34,0.36],[0.36,0.38],\ldots,[0.80,0.82],[0.82,0.827].$$
Dazu müssen wir jeweils ein $\lambda^{\star}$ wählen. Für die Fälle, in denen $\lambda_1 \leq 0.68$ ist, lesen wir ein $\lambda^{\star}$ aus \cite[Table 8, 10 (§9))]{Hea92} ab. Dazu bemerken wir, dass in \cite[Table 8 (§9)]{Hea92} in der Zeile mit $\lambda_1\leq 0.66$ ein Zahlendreher vorliegt. Die Abschätzung $\lambda' \geq 0.783$ müsste dort $\lambda' \geq 0.738$ lauten. Wir benutzen der Einfachheit halber (damit letzteres nicht erst überprüft werden muss) die Abschätzung $\lambda' \geq 0.714$ für $\lambda_1 \leq 0.68$, die ja insbesondere auch für $\lambda_1 \leq 0.66$ gültig ist.\\
Ist $\lambda_1 \geq 0.68$, so nehmen wir $\lambda^{\star}=\lambda_1$. Wir halten fest, dass die rechte Seite von (\ref{para9-l-dash-hauptfall}) auch mit dieser Wahl von $\lambda^{\star}$ monoton wachsend in $\lambda_1$ und $\lambda'$ ist, denn
$$F(-\lambda_1)-F(\lambda'-\lambda_1)=\int_0^{2\gamma}f(x) e^{\lambda_1 x}(1-e^{-\lambda'x})\, dx .$$
Weiterhin wählen wir in den Fällen mit $\lambda_1 \leq 0.68$ die Parameter (\ref{para9-l-dash-vor0}) - (\ref{para9-l-dash-vor3}) und für $\lambda_1 \geq 0.68$ wählen wir (\ref{para9-l-dash-vor0}) und
\begin{eqnarray*}
 & & k_1=k,\; k_2=0,\; k_3=k^2+\frac{3}{4},\\
 & & s_1=\lambda_1, \; s_{11}=\lambda_{11}, \; s_{12}=\lambda_{12}, \;(s_2=0,) \\
 & &\Delta s_1=0.004,\; \Delta s_2=0, \; \Delta t =0.004,\; x_1=15 .
\end{eqnarray*}\\

\noindent \textbf{Bemerkungen zur Tabelle für Fall 2:}\\
Analog zum obigen Beispiel erstellen wir aus der Ungleichung (\ref{para9-l-dash-ungleichung-nebenfall}) eine Tabelle für Fall 2 für die Intervalle
$$\lambda_1 \in [\lambda_{11},\lambda_{12}]=[0.34,0.38],[0.38,0.42],\ldots,[1.02,1.06],[1.06,1.099].$$
Aufgrund der besseren erzielten Werte, reicht es hier $\lambda_1$ in $0.04$-Schritten zu staffeln.\\
Wieder gilt, dass die rechte Seite von (\ref{para9-l-dash-ungleichung-nebenfall}) ohne die beiden Suprema monoton wachsend in $\lambda_1$ und $\lambda'$ ist.

Diesmal müssen wir zwei Suprema abschätzen, $A$ stehe für das erste und $B$ für das zweite Supremum aus (\ref{para9-l-dash-ungleichung-nebenfall}) (die Nummerierung laufe dabei immer von links nach rechts). Wir wählen die Parameter
$$\gamma=1.21- \frac{5\lambda_{12} }{12} \; \text{ und } \; k=0.77 + \frac{\lambda_{12}}{10} . $$
Für die Supremumanalyse aus §\ref{para-rechenkapitel} setzt man die folgenden Werte. Dabei hat man für $k_i$ bei den beiden Suprema $A$ und $B$ verschiedene Werte, jedoch für $s_i$ und $s_{ij}$ die gleichen:
\begin{eqnarray*}
 & & k_{1,A}=2k,\; k_{2,A}=0,\; k_{3,A}=2(k^2+\frac{3}{4}), \\
 & & k_{1,B}=\frac{1}{2},\; k_{2,B}=0,\; k_{3,B}=2k, \\
 & & s_1=\lambda_1, \; s_{11}:=\lambda_{11}, \;s_{12}:=\lambda_{12},\; (s_2=0,) \\
 & & \Delta s_1= 0.004, \; \Delta s_2=0, \; \Delta t= 0.004, \; x_1=15 .
\end{eqnarray*}

\noindent Bevor wir jetzt zu den beiden angekündigten Tabellen kommen, möchten wir noch folgende Erklärungen diesbezüglich notieren, welche auch in allen analogen späteren Tabellen Anwendung finden:
\begin{enumerate}
\item  Die folgenden Tabellen liefern Aussagen der Form $\lambda_1 \leq \lambda_{12} \Rightarrow \lambda' > c_2$. Jedoch haben wir eigentlich nur Aussagen der Form $\lambda_1 \in [\lambda_{11},\lambda_{12}] \Rightarrow \lambda' > c_2$ bewiesen. Da aber die Abschätzungen für kleinere $\lambda_1$ besser sind und da weiterhin der Fall $\lambda_1 \leq 0.34$ nach \cite[Theorem 1 (S.268)]{Hea92} nicht eintritt, ist alles korrekt.
\item Wir haben die Abschätzungen der Suprema $A\leq C$ (Fall 1) bzw. $A \leq C_1$ und $B \leq C_2$ (Fall 2) der jeweiligen Tabelle beigefügt. Man beachte dabei, dass wir in der Tabelle die aufgerundeten Werte der Abschätzungen $C$ bzw. $C_1$ und $C_2$ notieren. Die zweite Zeile in Tabelle 2  besagt z.B., dass unter den bekannten Voraussetzungen die Aussage gilt
      $$\lambda_1 \leq 0.38 \Rightarrow \lambda' > 1.96 .$$
    Zum Beweis dessen haben wir $\lambda^{\star}=0.887$ benutzt. Außerdem haben wir mittels der Analyse in §\ref{para-rechenkapitel} bewiesen, dass
    $$\sup_{t \in \R, \atop \lambda_1 \in [0.36,0.38]} Re  \left\{k F(-\lambda^{\star}+i t)- (k^2+\frac{3}{4})F(\lambda_1-\lambda^{\star} + i t) \right\} \leq C $$
    für eine gewisse konkrete Konstante $C \leq 0.0134$.
\end{enumerate}
\vspace*{0.25cm}\label{Tabelle2-3}
\begin{center}
\begin{tabular}[t]{|ll|ll|}
\multicolumn{4}{c}{\textbf{Tabelle 2. $\lambda'$-Abschätzungen}}\\
\multicolumn{4}{c}{($\chi_1$ oder $\rho_1$ komplex,  $ord \; \chi_1 \geq 5$)}\\\hline
$\lambda_1 \leq $ & $\lambda' >$ & $\lambda^{\star}$ &  $C \leq$ \\ \hline
0.36 & 2.06 & 0.903 & 0.0172\\
0.38 & 1.96 & 0.887 & 0.0134\\
0.40 & 1.86 & 0.871 & 0.0102\\
0.42 & 1.77 & 0.856 & 0.0074\\
0.44 & 1.69 & 0.842 & 0.0049\\
0.46 & 1.61 & 0.829 & 0.0032\\
0.48 & 1.53 & 0.816 & 0.0028\\
0.50 & 1.47 & 0.803 & 0.0025\\
0.52 & 1.40 & 0.791 & 0.0021\\
0.54 & 1.34 & 0.780 & 0.0018\\
0.56 & 1.28 & 0.769 & 0.0015\\
0.58 & 1.23 & 0.759 & 0.0012\\
0.60 & 1.18 & 0.749 & 0.0009\\
0.62 & 1.13 & 0.739 & 0.0008\\
0.64 & 1.09 & 0.730 & 0.0008\\
0.66 & 1.04 & 0.714 & 0.0007\\
0.68 & 1.00 & 0.712 & 0.0007\\
0.70 & 0.96 &  & 0.0012\\
0.72 & 0.93 &  & 0.0011\\
0.74 & 0.91 &  & 0.0010\\
0.76 & 0.89 &  & 0.0009\\
0.78 & 0.86 &  & 0.0008\\
0.80 & 0.84 &  & 0.0007\\
0.82 & 0.83 &  & 0.0006\\
0.827 & 0.827 &  & 0.0005\\
\hline
\end{tabular}
\hspace*{0.5cm}
\begin{tabular}[t]{|ll|ll|}
\multicolumn{4}{c}{\textbf{Tabelle 3. $\lambda'$-Abschätzungen}}\\
\multicolumn{4}{c}{($\chi_1$ oder $\rho_1$ komplex,  $ord \; \chi_1 \in \{2,3,4 \}$)}\\\hline
$\lambda_1 \leq $ & $\lambda' >$ & $  C_1 \leq$& $ C_2 \leq$\\ \hline
0.38 & 2.53 & 0.0060 & 0.0027\\
0.42 & 2.35 & 0.0051 & 0.0024\\
0.46 & 2.20 & 0.0043 & 0.0020\\
0.50 & 2.06 & 0.0035 & 0.0017\\
0.54 & 1.94 & 0.0028 & 0.0015\\
0.58 & 1.84 & 0.0021 & 0.0012\\
0.62 & 1.75 & 0.0015 & 0.0010\\
0.66 & 1.67 & 0.0010 & 0.0008\\
0.70 & 1.59 & 0.0006 & 0.0006\\
0.74 & 1.52 & 0.0006 & 0.0004\\
0.78 & 1.46 & 0.0006 & 0.0003\\
0.82 & 1.40 & 0.0006 & 0.0002\\
0.86 & 1.35 & 0.0006 & 0.0002\\
0.90 & 1.30 & 0.0006 & 0.0002\\
0.94 & 1.25 & 0.0006 & 0.0002\\
0.98 & 1.21 & 0.0006 & 0.0002\\
1.02 & 1.17 & 0.0006 & 0.0002\\
1.06 & 1.13 & 0.0006 & 0.0002\\
1.099 & 1.099 & 0.0006 & 0.0002\\
\hline
\end{tabular}
\end{center}
\vspace*{0.5cm}

\noindent Tabelle 3 liefert bessere Werte als Tabelle 2. Damit gilt Tabelle 2 für beide Fälle und ersetzt \cite[Table 8 (§9)]{Hea92}. Außerdem liefert nun Tabelle 2 eine geringfügige Verbesserung von \cite[Theorem 2a]{Hea92}, wobei die Konstante $0.696$ auf $0.702$ verbessert wird, was aus den folgenden drei Punkten folgt:
\begin{itemize}
 \item Ist $\chi_1$ oder $\rho_1$ komplex so liefert Tabelle 2, dass $\lambda' \geq 0.827$. Weiterhin liefert \cite[Lemma 9.4]{Hea92}, dass\footnote{Zwar steht in \cite[Lemma 9.4]{Hea92} die Aussage $\lambda_2 \geq 0.702$. Es wurde aber $\lambda_2>0.702$ bewiesen. In der Regel gilt für jegliche Abschätzungen aus \cite{Hea92}, in denen als bewiesene Abschätzung $\lambda \geq c$ präsentiert wird, dass in Wirklichkeit $\lambda > c$ gezeigt wurde. Im Übrigen folgt aus Stetigkeitsgründen sowieso, dass wenn mit der obigen Methode gezeigt wurde, dass $\lambda \geq c$, dann gilt auch $\lambda \geq c+\varepsilon$ für ein hinreichend kleines $\varepsilon>0$.} $\lambda_2 > 0.702$. Also kann wegen $L_1 \geq 1$ die Funktion
      $$\prod_{\chi \neq \chi_0} L(s,\chi)$$
      höchstens die beiden einfachen Nullstellen $\rho_1$ und $\overline{\rho_1}$ in
\begin{equation}\label{Theorem2a-new}\sigma \geq 1- \frac{0.702}{\mathscr{L}}, \; |t| \leq 1 \end{equation}
haben.
 \item Sind $\chi_1$ und $\rho_1$ beide reell, so liefern \cite[Lemma 8.4, Lemma 8.8]{Hea92} die hinreichenden Abschätzungen $\lambda' \geq 1.294$ und $\lambda_2 \geq 0.745$.
 \item Für die Funktion $L(s,\chi_0)$ gilt für hinreichend großes $q$ die Nullstellenfreiheit in (\ref{Theorem2a-new}) wegen der Nullstellenfreiheit der Riemannschen $\zeta$-Funktion auf $Re \{s\}=1$ und dem bekannten Sachverhalt, dass
     \begin{equation}\label{L-s-chi0-zeta-s}L(s,\chi_0)=\zeta(s) g(s,\chi_0)\end{equation}
     mit einer auf $Re  \{s\}>0$ nullstellenfreien Funktion $g(s,\chi_0)$.
\end{itemize}

\noindent Schlussendlich noch eine Bemerkung bezüglich zweier kleiner Verbesserungspotentiale:\\
Hätte man $x_1$ größer und $\Delta s_i,\; \Delta t$ kleiner gewählt, sowie die Parameter $\gamma$ und $k$ für jeden einzelnen Fall separat optimiert, so hätten die Abschätzungen teilweise um wenige "`Hundertstel"' verbessert werden können. Dieser Hinweis gilt für alle Fälle in dieser Arbeit, in der ein Supremum mittels §\ref{para-rechenkapitel} abgeschätzt wird.\\
Ein weiteres (und größeres) Verbesserungspotential für die Tabelle 2 liegt darin, die verbesserten Abschätzungen für $\lambda'$ (welche eben bewiesen wurden) und $\lambda_2$ (welche im nächsten Abschnitt bewiesen werden) zu nehmen. Dann bekommt man bessere Werte für $\lambda^{\star}$ und folglich bessere Werte in Tabelle 2.\\
Da wir jedoch (im Nachhinein) für $\lambda'$ nur knapp bessere Abschätzungen benötigen als für $\lambda_2$, und dies mit Tabelle 2 bereits hinreichend gegeben ist, berücksichtigen wir die beiden eben genannten Punkte nicht.

\section{$\lambda_2$-Abschätzungen}\label{para-l2-abschaetzungen}
Wie im letzten Abschnitt so gelte auch in diesem, dass mindestens eines der beiden Elemente $\rho_1$, $\chi_1$ komplex sei.\\
In \cite[Table 9 (§9)]{Hea92} benutzt Heath-Brown die in §\ref{hb-para-6} angesprochene Methode mit einer gewissen Anfangsungleichung und $\beta=\beta_1$. Er erhält daraus Abschätzungen für $\lambda_3$ bzw. auch für $\lambda_2$, wenn gewisse Sonderfälle ausgeschlossen werden. In Verbesserungspotential Nr. 5 \cite[S.334]{Hea92} wird darauf hingewiesen, dass man in der Herleitung auch $\beta=1-\mathscr{L}^{-1} \min \{\lambda_2,\lambda'\}$ hätte wählen können. Dies würde in Abschätzungen für $\lambda_4$ resultieren bzw. auch solchen für $\lambda_3$ und $\lambda_2$, wenn gewisse Sonderfälle ausgeschlossen werden würden. \\
Wir benutzen in diesem Abschnitt diesen Ansatz aus Verbesserungspotential Nr. 5. Weiterhin kann man mittels Verbesserungspotential Nr. 2 auch die eben angesprochenen Sonderfälle behandeln. Das führt dazu, dass wir letztendlich gute Verbesserungen für die $\lambda_2$-Abschätzungen bekommen. Wir sollten erwähnen, dass diese verbesserten $\lambda_2$-Abschätzungen den größten Beitrag für die Verbesserung der Linnikschen Konstante von $L=5.5$ auf $L=5.2$ liefern.

\subsection{Das wesentliche Lemma}
In diesem Abschnitt beweisen wir die folgende verfeinerte Version von \cite[Lemma 9.2 (S.306)]{Hea92}.

\begin{lemma}\label{para9-zentrales-neues-lemma} Sei $\chi_1$ oder $\rho_1$ komplex und $f$ eine Funktion, die Bedingung 1 und 2 erfüllt. Ferner seien $j \in \{2,3\}$ und $\lambda^{\star}$ eine Konstante mit  $0<\lambda^{\star} \leq \min \{\lambda',\lambda_2 \}$. Wenn $\rho'$ nicht existiert, so setzen wir nur $\lambda^{\star} \leq \lambda_2$ voraus. Schließlich seien $\varepsilon$ und $k$ positive Konstanten. Für $q \geq q_0(\varepsilon,f,k)$ gilt dann
\begin{equation}\label{l2-Hauptungleichung}
0 \leq (k^2+\frac{1}{2})\left(F(-\lambda^{\star})-F(\lambda_j-\lambda^{\star}) \right) - 2k F(\lambda_1 - \lambda^{\star})  + D + \varepsilon,
\end{equation}
wobei

$$D = \left\lbrace \begin{array}{ll}
\frac{f(0)}{6}(k^2+4k+\frac{3}{2})&\text{falls } \chi_1^2, \chi_1^3 \neq \chi_0,\chi_j,\overline{\chi_j}, \\
\\
A+\frac{f(0)}{6}(k^2+4k+\frac{5}{4})  &\chi_1^2 \in \{\chi_j,\overline{\chi_j}\} \text{ und } ord \;  \chi_1 \geq 6, \\
\\
2A+\frac{f(0)}{8}(k^2+4k+1) &\chi_1^2 \in \{\chi_j,\overline{\chi_j}\} \text{ und } ord \;  \chi_1=4, \\
\\
A+B+\frac{f(0)}{8}(k^2+4k+\frac{5}{4})&\chi_1^2 \in \{\chi_j,\overline{\chi_j}\} \text{ und } ord \;  \chi_1=5, \\
\\
B+\frac{f(0)}{6}(k^2+4k+\frac{3}{2}) &\chi_1^3 \in \{\chi_j,\overline{\chi_j}\} \text{ und } ord \;  \chi_1 \geq 7, \\
\\
2B+\frac{f(0)}{8}(k^2+4k+\frac{3}{2}) &\chi_1^3 \in \{\chi_j,\overline{\chi_j}\} \text{ und } ord \;  \chi_1=6, \\
\\
2A +\frac{f(0)}{6}(k^2+\frac{7}{2} k+1)& ord \;  \chi_1 =2,\\
\\
2B + \frac{f(0)}{6}(k^2+\frac{7}{2} k+\frac{11}{8})& ord \;  \chi_1 =3\\
  \end{array} \right.
$$
und
\begin{eqnarray*}
A&=&\sup_{t \in \R} Re  \left\{ \frac{1}{4}F(-\lambda^{\star}+it)-kF(\lambda_1-\lambda^{\star}+it) \right\},\\
B&=& \sup_{t \in \R} Re  \left\{-\frac{1}{4}F(\lambda_1-\lambda^{\star}+it) \right\}.
\end{eqnarray*}

\noindent Andere als die oben bei der Gleichung für $D$ genannten acht Fälle gibt es nicht.
\end{lemma}

\noindent Für den Beweis nummerieren wir die Fälle nacheinander von eins bis acht durch. Fall 3 ist beispielsweise der Fall $\chi_1^2 \in \{\chi_j,\overline{\chi_j}\}$ und $ord \;  \chi_1=4$. \\
Die Aussage im Lemma ist wohldefiniert, denn die acht Fälle überlappen sich nicht. Dazu erwähnen wir, dass Fall 2 sich nicht mit Fall 5 oder 6 überlappt. Wären wir nämlich gleichzeitig in Fall 2 und 5 oder in Fall 2 und 6 so folgte $\chi_1^3 \in \{\chi_1^2,\overline{\chi_1}^2\}$, also $\chi_1=\chi_0$ oder $ord \; \chi_1 \leq 5$, was beides ausgeschlossen ist. Die "`paarweise Disjunktheit"' der restlichen Fälle ist noch offensichtlicher.\\
Weiterhin wurden alle Fälle genannt, denn: Angenommen, wir würden uns in keinem der acht obigen Fälle befinden. Wegen Fall 1,7 und 8 bedeutet dies, dass $ord \;  \chi_1 \geq 4$ und
$\chi_1^2 \in \{\chi_j,\overline{\chi_j}\}$ oder $\chi_1^3 \in \{\chi_j,\overline{\chi_j}\}$. Da wir uns aber auch nicht in Fall 2-6 befinden, folgt  $\chi_1^3 \in \{\chi_j,\overline{\chi_j}\}$ und $ord \;  \chi_1 \in [4,5]$. Wäre $ord \;  \chi_1 =4$, dann folgt
$$\overline{\chi_1} =\chi_1^3 \in \{\chi_j,\overline{\chi_j}\},  $$
was wegen $j \neq 1$ ein Widerspruch zur Definition der $\chi_j$ ist. Wäre $ord \;  \chi_1 =5$, so folgt
$$\chi_1^2=\overline{\chi_1}^3 \in \{\overline{\chi_j},\chi_j\},$$
dann sind wir aber im Fall 4, Widerspruch und fertig mit der Begründung.\\

Wir brauchen noch eine weitere wichtige Bemerkung. Man betrachte dazu beispielsweise den Fall 2, also ist $\chi_1^2 \in \{\chi_j,\overline{\chi_j} \}$. Dann nehmen wir oBdA an, dass $\chi_1^2=\chi_j$, da $\chi_j$ und $\overline{\chi_j}$, als sie definiert wurden, austauschbar waren. Der Grund für diese "`Austauschbarkeit"' liegt in den einfachen Fakten
$$\rho \in R \text{ Nullstelle von } L(s,\chi_j) \Longleftrightarrow \overline{\rho} \in R \text{ Nullstelle von } L(s,\overline{\chi_j}),$$
$$Re \{\rho\}=Re \{\overline{\rho}\} $$
und für $i,j \in \N$ gilt
$$ \chi_j \neq \chi_i,\overline{\chi_i} \Longleftrightarrow  \overline{\chi_j} \neq \chi_i,\overline{\chi_i}.$$
Wir könnten also den Charakter $\overline{\chi_j}$ auch $\chi_{j,neu}$ nennen und die ganze Herleitung nochmal durchgehen. Dann würden wir etwas beweisen für den Fall $\chi_1^2=\chi_{j,neu}=\overline{\chi_j}$. Damit gilt die Ungleichung, die wir aus der Voraussetzung $\chi_1^2=\chi_j$ folgern, auch für den Fall $\chi_1^2=\overline{\chi_j}$.\\
Sollte man sich mit dieser Begründung nicht wohlfühlen, so kann man die folgende Analyse, die wir für den Fall $\chi_1^2=\chi_j$ machen, auch für den Fall $\chi_1^2=\overline{\chi_j}$ durchführen und sehen, dass die gleiche Ungleichung für das Lemma folgt.\\
Die gleiche Überlegung gilt für die Fälle 3,4,5 und 6. Also betrachten wir jeweils nur den Fall $\chi_1^2=\chi_j$ bzw. $\chi_1^3=\chi_j$ anstelle von $\chi_1^2 \in \{\chi_j,\overline{\chi_j} \}$ bzw. $\chi_1^3 \in \{\chi_j,\overline{\chi_j} \}$.\\

Nach diesen Vorbemerkungen möchten wir nun (\ref{l2-Hauptungleichung}) beweisen. Dazu betrachten wir die Ungleichung in \cite[S.306 2.Zeile]{Hea92}, welche nach der Standardmethode aus \cite[S.305 letzte Zeile]{Hea92} gefolgert wird (man benutze für diesen Schritt die Beziehung $Re\{x\} Re\{y\}=\frac{1}{2}(Re\{xy\}+Re\{x\overline{y}\})$). Also gilt diese Ungleichung auch, wenn man darin $\beta_1$ durch irgendeine andere komplexe Zahl $\beta \in \C$ ersetzt. Wie in Verbesserungspotential Nr. 5 vorgeschlagen wird, ersetzen wir $\beta_1$ durch $\beta^{\star}=1-\mathscr{L}^{-1} \lambda^{\star}$, wobei $\lambda^{\star}$ wie im Lemma definiert sei. Es folgt die Ungleichung
\begin{eqnarray}
 0 &\leq& (k^2+\frac{1}{2}) K(\beta^{\star},\chi_0) + (k^2+\frac{1}{2}) K(\beta^{\star}+i\gamma_j,\chi_j) +2k K(\beta^{\star}+i\gamma_1,\chi_1) \nonumber \\
 & &+ k K(\beta^{\star}+i(\gamma_1+\gamma_j),\chi_1 \chi_j)+ k K(\beta^{\star}+i(-\gamma_1+\gamma_j),\overline{\chi_1} \chi_j) \nonumber \\
 & &+ \frac{1}{4} K(\beta^{\star}+i(2\gamma_1+\gamma_j),\chi_1^2 \chi_j)+ \frac{1}{4} K(\beta^{\star}+i(-2\gamma_1+\gamma_j),\overline{\chi_1}^2 \chi_j) \nonumber \\
 & &+\frac{1}{2} K(\beta^{\star}+2i \gamma_1,\chi_1^2 ). \label{para9-l2-zentrale-ungleichung}
 \end{eqnarray}

\noindent Mit identischem Vorgehen wie beim Beweis von Lemma \ref{lemma-lstrich-ordchi-geq-5} und Lemma \ref{lemma-lstrich-ordchi-234} folgern wir nun aus (\ref{para9-l2-zentrale-ungleichung}) die Ungleichungen des Lemmas (für die acht verschiedenen Fälle). Wir fassen dir Vorgehensweise zusammen:\\
Man muss bei der Abschätzung eines Terms $K(\beta^{\star}+it,\chi)$ jeweils kontrollieren, ob der Charakter $\chi$ gleich $\chi_0,\; \chi_1$ oder $\overline{\chi_1}$ ist. Dann geht man zur Herleitung des Lemmas folgendermaßen vor:
\begin{itemize}
\item Ist $\chi=\chi_0$ so benutzt man Lemma \ref{hb-lemma-5-3}.
\item Ist $\chi=\chi_1$ so benutzt man Lemma \ref{hb-lemma-5-2-extended}. Dort muss man die Nullstelle $\rho_1$ berücksichtigen. Ist $\rho_1$ nicht bereits in der Menge $A_1$, so füge man sie in die Menge $A_2$ ein. Beim Fall $ord \; \chi_1=2$ müssen wir zusätzlich die Nullstelle $\overline{\rho_1} \neq \rho_1$ berücksichtigen.\\
    Wegen $\beta' \leq \beta^{\star}$ müssen keine weitere Nullstellen berücksichtigt werden. Wir wählen also insgesamt die Menge $A_2$ so, dass
    $$A_1 \cup A_2 = \left\lbrace \begin{array}{ll}
\{\rho_1, \overline{\rho_1}\} & \; \text{ für } ord \; \chi_1=2,\\
\{\rho_1\} & \; \text{ sonst.}
 \end{array} \right. $$
\item Ist $\chi=\overline{\chi_1}$, so berücksichtigt man die Nullstelle $\overline{\rho_1}$ (und ggf. auch $\rho_1$) und geht völlig analog zum Fall $\chi=\chi_1$ vor. Es folgt
    $$A_1 \cup A_2 = \left\lbrace \begin{array}{ll}
\{\overline{\rho_1}, \rho_1\} & \; \text{ für } ord \; \overline{\chi_1}=2,\\
\{\overline{\rho_1}\} & \; \text{ sonst.}
 \end{array} \right. $$
\item Sei $\chi \neq \chi_0,\; \chi_1,\; \overline{\chi_1}$. Wieder benutzt man Lemma \ref{hb-lemma-5-2-extended}. Diesmal gilt für alle Nullstellen $\rho \in R$ von $L(s,\chi)$ die Ungleichung $Re \{\rho\} \leq Re  \{\rho_2\} \leq \beta^{\star}$. Also ist $A_1=\emptyset$. Ist zusätzlich $\chi=\chi_j$, so wählen wir $A_2=\{\rho_j \}$ oder $A_2=\emptyset$. Für $\chi \neq \chi_j$ wählen wir immer $A_2=\emptyset$.
\item Für die Fälle in denen bekannt ist, dass $ord \; \chi_1 \leq 6$ benutzen wir in Lemma \ref{hb-lemma-5-2-extended} den Wert $\phi(\chi_1)=\phi(\chi_1^2)=\frac{1}{4}$. Gilt zusätzlich $\chi_j \in \{\chi_1^2, \chi_1^3 \}$ so benutzen wir für alle in (\ref{para9-l2-zentrale-ungleichung}) auftauchenden Charaktere $\chi \neq \chi_0$ den Wert $\phi(\chi)=\frac{1}{4}$. In den übrigen Fällen verwenden wir immer $\phi(\chi)\leq \frac{1}{3}$.
\end{itemize}

\noindent Wir schätzen nun die einzelnen "`$K$-Terme"' ab. Im Folgenden beachte man, dass nach Definition $\chi_j \neq \chi_0,\, \chi_1,\, \overline{\chi_1} $ gilt. Mit der Bezeichnung $K_0(s,\chi)=\mathscr{L}^{-1} K(s,\chi)-\varepsilon$ gilt

\begin{eqnarray*}
 K_0(\beta^{\star},\chi_0) &\leq& F(-\lambda^{\star}) \; \text{ für Fall 1-8,} \\
  \\
  K_0(\beta^{\star}+i \gamma_j,\chi_j) &\leq& \left\lbrace \begin{array}{ll}
-F(\lambda_j-\lambda^{\star})+ \frac{f(0)}{6}& \; \text{ Fall 1,2,5,7,8,}\\
-F(\lambda_j-\lambda^{\star})+ \frac{f(0)}{8}& \; \text{ Fall 3,4,6,}\\
 \end{array} \right.\\
 \\
  K_0(\beta^{\star}+i \gamma_1,\chi_1) &\leq& \left\lbrace \begin{array}{ll}
-F(\lambda_1-\lambda^{\star})+ \frac{f(0)}{6}& \; \text{ Fall 1,2,5,}\\
-F(\lambda_1-\lambda^{\star})+ \frac{f(0)}{8}& \; \text{ Fall 3,4,6,8,}\\
-F(\lambda_1-\lambda^{\star})-Re \{F(\lambda_1-\lambda^{\star}+2i\mu_1)\}+ \frac{f(0)}{8}& \; \text{ Fall 7.}
 \end{array} \right.
\end{eqnarray*}\\

\noindent Bei den restlichen Termen "`$K(s,\chi)$"' aus (\ref{para9-l2-zentrale-ungleichung}) gilt es zu untersuchen, wann der Charakter $\chi$ gleich $\chi_0$, $\chi_1$ oder $\overline{\chi_1}$ ist. Unter Berücksichtigung der obigen Zusatzbemerkung (also für $l \in \{2,3\}$ oBdA $\chi_1^l=\chi_j$ im Falle von $\chi_1^l \in \{\chi_j,\overline{\chi_j} \}$) hat man

$$
\begin{array}{ll}
\chi_1^2 \chi_j =\left\lbrace \begin{array}{ll}
\chi_1^2 \chi_j \neq \chi_0,\chi_1,\overline{\chi_1} & \; \text{ Fall 1},\\
\chi_1^4 \neq \chi_0,\chi_1,\overline{\chi_1} & \; \text{ Fall 2},\\
\chi_1^4 =\chi_0 & \; \text{ Fall 3},\\
\chi_1^4 =\overline{\chi_1} & \; \text{ Fall 4},\\
\chi_1^5 \neq \chi_0,\chi_1,\overline{\chi_1} & \; \text{ Fall 5},\\
\chi_1^5 =\overline{\chi_1} & \; \text{ Fall 6},\\
\chi_j \neq \chi_0,\chi_1,\overline{\chi_1} & \; \text{ Fall 7},\\
\overline{\chi_1} \chi_j \neq \chi_0,\chi_1,\overline{\chi_1} & \; \text{ Fall 8},
 \end{array} \right.
&
\overline{\chi_1}^2 \chi_j =\left\lbrace \begin{array}{ll}
\overline{\chi_1}^2 \chi_j \neq \chi_0,\chi_1,\overline{\chi_1} & \; \text{ Fall 1},\\
\chi_0 & \; \text{ Fall 2},\\
\chi_0 & \; \text{ Fall 3},\\
\chi_0 & \; \text{ Fall 4},\\
\chi_1 & \; \text{ Fall 5},\\
\chi_1 & \; \text{ Fall 6},\\
\chi_j \neq \chi_0,\chi_1,\overline{\chi_1} & \; \text{ Fall 7},\\
\chi_1 \chi_j \neq \chi_0,\chi_1,\overline{\chi_1} & \; \text{ Fall 8.}
 \end{array} \right.
 \end{array}
$$\\

\noindent Also folgt
\begin{eqnarray*}
K_0(\beta^{\star}+i(2\gamma_1+ \gamma_j),\chi_1^2 \chi_j) &\leq&  \left\lbrace \begin{array}{ll}
\frac{f(0)}{6}&  \; \text{ Fall 1,2,5,7,8,}\\
Re\{F(-\lambda^{\star}+i(2\mu_1+\mu_j))\}& \; \text{ Fall 3,}\\
-Re\{F(\lambda_1-\lambda^{\star}+i(3\mu_1+\mu_j))\}+ \frac{f(0)}{8}&  \; \text{ Fall 4,6,}
 \end{array} \right.\\
 \\
 K_0(\beta^{\star}+i(-2\gamma_1+ \gamma_j),\overline{\chi_1}^2 \chi_j) &\leq& \left\lbrace \begin{array}{ll}
\frac{f(0)}{6}& \; \text{Fall 1,7,8,}\\
Re\{F(-\lambda^{\star}+i(-2\mu_1+\mu_j))\}& \; \text{Fall 2,3,4,}\\
-Re\{F(\lambda_1-\lambda^{\star}+i(-3\mu_1+\mu_j))\}+ \frac{f(0)}{6}& \; \text{Fall 5,}\\
-Re\{F(\lambda_1-\lambda^{\star}+i(-3\mu_1+\mu_j))\}+ \frac{f(0)}{8}& \; \text{Fall 6.}
 \end{array} \right.\\
\end{eqnarray*}\\

\noindent Analog schließt man
\begin{eqnarray*}
 K_0(\beta^{\star}+i(\gamma_1+ \gamma_j),\chi_1 \chi_j) &\leq& \left\lbrace \begin{array}{ll}
\frac{f(0)}{6}& \; \text{ Fall 1,2,5,7,8,}\\
-Re\{F(\lambda_1-\lambda^{\star}+i(2\mu_1+\mu_j))\}+ \frac{f(0)}{8}&  \; \text{ Fall 3,}\\
\frac{f(0)}{8}&  \; \text{ Fall 4,6,}\\
 \end{array} \right.\\
 \\
 K_0(\beta^{\star}+i(-\gamma_1+ \gamma_j),\overline{\chi_1} \chi_j) &\leq& \left\lbrace \begin{array}{ll}
\frac{f(0)}{6}&  \; \text{ Fall 1,5,7,8,}\\
-Re\{F(\lambda_1-\lambda^{\star}+i(-2\mu_1+\mu_j))\}+ \frac{f(0)}{6}&  \; \text{ Fall 2,}\\
-Re\{F(\lambda_1-\lambda^{\star}+i(-2\mu_1+\mu_j))\}+ \frac{f(0)}{8}&  \; \text{ Fall 3,4,}\\
\frac{f(0)}{8}&  \; \text{ Fall 6,}\\
 \end{array} \right.\\
 \\
 K_0(\beta^{\star}+2i\gamma_1,\chi_1^2) &\leq& \left\lbrace \begin{array}{ll}
\frac{f(0)}{6}& \; \text{Fall 1,2,5,}\\
\frac{f(0)}{8}& \; \text{Fall 3,4,6,}\\
Re\{F(-\lambda^{\star}+2i\mu_1)\}& \; \text{Fall 7,}\\
-Re\{F(\lambda_1-\lambda^{\star}+3i\mu_1)\}+ \frac{f(0)}{8}& \; \text{Fall 8.}
 \end{array} \right.\\
\end{eqnarray*}

\noindent Trägt man für einen der acht Fälle die obigen acht Ungleichungen zusammen, passt das $\varepsilon$ an und beachtet man noch ggf. (\ref{Im-z-Vorzeichen-in-Re-z-unwichtig}), so folgt die entsprechende Ungleichung des Lemmas. Es folgt das Lemma.\\

\begin{bemerkung}Für den Beweis von Lemma \ref{para9-zentrales-neues-lemma} starteten wir mit der Ungleichung \cite[S.305, letzte Zeile]{Hea92}, also
  $$0 \leq \chi_0(n)\left(1+ Re  \left\{ \frac{\chi_j(n)}{n^{i\gamma_j}}\right\} \right) \left(k+ Re  \left\{ \frac{\chi_1(n)}{n^{i\gamma_1}}\right\} \right)^2. $$
  Wir haben getestet, was wir mit
  $$0 \leq \chi_0(n)\left(1+ Re  \left\{ \frac{\chi_j(n)}{n^{i\gamma_j}}\right\} \right) \left(k+ Re  \left\{ \frac{\chi_1(n)}{n^{i\gamma_1}}\right\} \right)^2 \left(1+ Re  \left\{ \frac{\chi_1(n)}{n^{i\gamma_1}}\right\}\right) $$
  bekämen. Das Ergebnis ist, dass man am Ende für die Fälle 5,6 und 8 bessere Abschätzungen für $\lambda_2$ im Bereich $\lambda_1 \leq 0.60$ bekommt. Diese haben wir letztendlich nicht nötig, also sparen wir uns dieses Vorgehen.   Man muss dabei beachten, dass durch die zusätzliche Klammer in der obigen Ungleichung die Anzahl der zu analysierenden Terme $K(s,\chi)$ zunimmt.
\end{bemerkung}
\vspace*{0.5cm}

\subsection{$\lambda_2$-Abschätzungen für die Fälle 1,2,3,4,6 und 8}

Mit Hilfe des eben bewiesenen Lemmas möchten wir zuerst Abschätzungen für $\lambda_2$ beweisen. Wir nehmen dafür an, dass $\lambda_2 \leq \lambda'$. Wäre $\lambda_2>\lambda'$, so hätten wir vermöge der $\lambda'$-Abschätzungen im letzten Abschnitt bereits sehr gute Abschätzungen. Nun benutzen wir Lemma \ref{para9-zentrales-neues-lemma} mit $j=2$ und $\lambda^{\star}=\lambda_2$. Man könnte jetzt die acht verschiedenen Ungleichungen aus dem Lemma hernehmen und für jeden Fall eine Tabelle aufstellen, analog zu Tabelle 2 und 3. Das Minimum über die jeweiligen Einträge der acht Tabellen wäre dann unsere allgemein gültige Abschätzung für $\lambda_2$.\\
Wir sparen uns jedoch etwas Arbeit, indem wir die Fälle 2,3,4,6 und 8 sozusagen gleichzeitig behandeln auf Kosten von geringfügig schlechteren Werten. Die Fälle 5 und 7 behandeln wir mit Hinblick auf §\ref{para9-l3bounds} und §\ref{Beweis-Theorem1-1-Fall2} separat. Wir beginnen mit

\subsubsection{Fall 1}
In diesem Fall gilt $\chi_1^2,\chi_1^3 \neq \chi_0,\chi_j,\overline{\chi_j}$. Nehme die dazugehörige erste Ungleichung aus Lemma \ref{para9-zentrales-neues-lemma}, das wäre
\begin{equation}\label{l2-ungleichung-fall1}0 \leq (k^2+\frac{1}{2})\left(F(-\lambda_2)-F(0) \right) - 2k F(\lambda_1 - \lambda_2) + \frac{f(0)}{6}(k^2+4k+\frac{3}{2})+ \varepsilon. \end{equation}
Die rechte Seite von (\ref{l2-ungleichung-fall1}) ist monoton wachsend in $\lambda_1$. Eine Monotonie in $\lambda_2$ lässt sich aber leider nicht so leicht ablesen. Aus diesem Grund muss man das Vorgehen im Vergleich zum Abschnitt mit den $\lambda'$-Abschätzungen ein wenig verändern (vergleiche \cite[S.307 oben]{Hea92}). Sei ein $\delta>0$ gegeben. Wir nehmen mal an, dass $$\lambda_1 \leq \lambda_{12}, \; \lambda_2 \in [\lambda_{2,alt},\lambda_{2,alt} +\delta], \; 2k-(k^2+\frac{1}{2}) \geq 0$$
 für gewisse konkrete Werte $\lambda_{11}$, $\lambda_{12}$ und $k$.
 Die Voraussetzung an das $k$ ist erfüllt, wenn \begin{equation}\label{l2-bed-an-k}k \in \left[1-\sqrt{\frac{1}{2}},1+\sqrt{\frac{1}{2}}\right]\supset [0.3,1.7].\end{equation} Es folgt aus (\ref{l2-ungleichung-fall1}), dass
\begin{eqnarray*}
  -\varepsilon &\leq& (k^2+\frac{1}{2})\left(F(-\lambda_2)-F(\lambda_{12}-\lambda_2) -F(0)\right)\\
   & &- (2k-(k^2+\frac{1}{2})) F(\lambda_{12} - \lambda_2) + \frac{f(0)}{6}(k^2+4k+\frac{3}{2}) \\
  &\leq & (k^2+\frac{1}{2})\left(F(-\lambda_{2,alt}-\delta)-F(\lambda_{12}-\lambda_{2,alt}-\delta) -F(0)\right) \\
  & &- (2k-(k^2+\frac{1}{2})) F(\lambda_{12} - \lambda_{2,alt}) + \frac{f(0)}{6}(k^2+4k+\frac{3}{2}).
\end{eqnarray*}
Den Term $F(\lambda_{12}-\lambda_2)$ haben wir dabei abgespalten, da dies für die Rechnungen vorteilhafter ist. Jetzt machen wir weiter wie immer: Bekommen wir für konkrete Werte von $\gamma,\; k,\; \lambda_{12},\; \lambda_{2,alt}$ und $\delta$ auf der rechten Seite der letzten Ungleichung etwas Negatives, so ist dies für hinreichend kleines $\varepsilon>0$ ein Widerspruch. Wissen wir zusätzlich, z.B. mittels \cite[Table 10 (§9)]{Hea92}, dass $\lambda_2 \geq \lambda_{2,alt}$, so haben wir bewiesen, dass
$$\lambda_1 \leq \lambda_{12} \Rightarrow \lambda_2 > \lambda_{2,alt}+\delta .$$
Das gleiche machen wir mit den Intervallen $\lambda_2 \in [\lambda_{2,alt}+j \delta,\lambda_{2,alt}+(j+1) \delta ]$ für $j \in \{1,2,\ldots,[(\lambda_{2,neu}-\lambda_{2,alt})/\delta] \}.$ Erhalten wir jedes Mal etwas Negatives, so haben wir insgesamt bewiesen, dass
$$\lambda_1 \leq \lambda_{12} \Rightarrow \lambda_2 > \lambda_{2,neu}.$$

\noindent Die soeben genannte Prozedur führen wir durch für die Intervalle
$$\lambda_1 \in [\lambda_{11},\lambda_{12}]=[0.34,0.36],[0.36,0.38],\ldots,[0.68,0.70].$$
Wir wählen \begin{equation}\label{Wahl-des-gamma-Fall1-l2}\gamma= 0.42+\lambda_{12}, \;k=0.59+ \frac{2}{5}\lambda_{12}, \; \delta=0.0001.\end{equation}
Dieses $k$ erfüllt immer (\ref{l2-bed-an-k}).
Die Werte für $\gamma$ und $k$ wurden nach Computerexperimenten mit Blick auf die zu involvierenden Fälle 2,3,4,6 und 8 gewählt. Als $\lambda_{2,alt}$ nehmen wir die Abschätzungen aus \cite[Table 10/Lemma 9.4 (§9)]{Hea92}. Der zweite Eintrag in "`Table 10"' besagt beispielsweise, dass für $\lambda_1 \leq 0.36$ gilt, dass $\lambda_2 \geq 0.903$. Also benutzen wir für den Fall $\lambda_1 \in [0.34,0.36]$ den Wert $\lambda_{2,alt}=0.903$.\\
Damit beweist man nun für den Fall 1 die entsprechenden Abschätzungen für $\lambda_2$, die in der nachfolgenden Tabelle 4 notiert sind. Diese gelten zuerst nur für den Fall, dass $\lambda_2 \leq \lambda'$, was eine Voraussetzung weiter oben war. Wenn jedoch letzteres nicht der Fall ist, so hat man vermöge Tabelle 2 mindestens genauso gute Abschätzungen. Damit gelten die Werte der Tabelle 4 immer im Fall 1. Bevor wir nun die Tabelle notieren, besprechen wir noch ein paar andere Fälle, für die wir die gleichen Abschätzungen beweisen möchten.\\

\subsubsection{Involvierung der Fälle 2,3,4,6 und 8 in die Tabelle für den Fall 1}
Für Fall 2 gehen wir völlig analog zu Fall 1 vor und wählen auch die gleichen Werte für $\gamma,\; k$ und $\delta$ (siehe (\ref{Wahl-des-gamma-Fall1-l2})). Diesmal brauchen wir aber zusätzlich eine Abschätzung für das Supremum $A$. Wir setzen dafür voraus, dass $\lambda_1 \in [\lambda_{11},\lambda_{12}]$. Außerdem nehmen wir mal an, dass \cite[Table 10 (§9)]{Hea92} die Abschätzung $\lambda_2 \geq \lambda_{2,alt}$ liefert und dass wir für Fall 1 die Abschätzung $\lambda_2 \geq \lambda_{2,neu}$ bewiesen haben. Wir möchten nun die gleiche Abschätzung auch für Fall 2 beweisen.\\ Also wählen wir für die Abschätzung des Supremums $A$ mittels §\ref{para-rechenkapitel} die folgenden Werte (beachte, dass $\gamma$ und $k$ bereits gewählt wurden):
\begin{eqnarray}
  & & k_{1,A}=\frac{1}{4}, \; k_{2,A}=k, \; k_{3,A}=0, \label{para9-A-Wahl-fuer-k}\\
  & &s_1=\lambda_2,\; s_{11}=\lambda_{2,alt}, \; s_{12}=\lambda_{2,neu}, \nonumber \\
  & &s_2=\lambda_1, \; s_{21}=\lambda_{11}, \; s_{22}=\lambda_{12}, \nonumber \\
  & &\Delta s_1=0.015,\; \Delta s_2=0.007, \; \Delta t=0.015,\; x_1=7. \nonumber
\end{eqnarray}
Wir erhalten dann eine Abschätzung $A \leq C_1$ und können damit (hoffentlich bzw. wenn wir alle Werte vorher richtig gewählt haben) im gleichen Stil wie für Fall 1 beweisen, dass im Fall 2 gilt: $\lambda_2 > \lambda_{2,neu}$.\\

Wir können nun die Fälle 3,4,6 und 8 gleichzeitig mit diesem Fall 2 behandeln. Dafür brauchen wir eine Abschätzung des Supremum $B$. Wir benutzen für $\gamma,\; k, \; s_{ij},\; \Delta s_i,\; \Delta t$ und $x_1$ die gleichen Werte wie bei der Abschätzung für $A$ und ersetzen lediglich (\ref{para9-A-Wahl-fuer-k}) durch
  $$k_{1,B}=0, \; k_{2,B}=\frac{1}{4}, \; k_{3,B}=0.$$
Damit erhalten wir eine Abschätzung $B \leq C_2$.\\
Wenn nun
$$D(Fall~3):=2C_1 + \frac{f(0)}{8}(k^2+4k+1) \leq C_1 + \frac{f(0)}{6}(k^2+4k+\frac{5}{4})=:D(Fall~2),$$
dann folgt aus der Ungleichung, die gemäß Lemma \ref{para9-zentrales-neues-lemma} für Fall 3 gilt, diejenige Ungleichung, die wir verwendet haben, um $\lambda_2$-Abschätzungen für Fall 2 zu beweisen. Alles was wir auf Basis dieser letzteren Ungleichung bewiesen haben, gilt dann also auch für Fall 3. Analog müssen wir für die Fälle $j=4,6$ und $8$ lediglich kontrollieren, ob
\begin{equation}\label{para9-D-Ungl}D(Fall~j)\leq D(Fall~2)\end{equation} ist. Letzteres haben wir gemacht und es hat immer gestimmt (dazu müssen wir sagen, dass wir die zu beweisenden $\lambda_2$-Abschätzungen extra so wählten, damit dies alles durchgeht). Insgesamt folgt
\label{Tabelle-4}
\begin{center}
\begin{tabular}[t]{|lc|cll|}
\multicolumn{5}{c}{\textbf{Tabelle 4. $\lambda_2$-Abschätzungen}}\\
\multicolumn{5}{c}{($\chi_1$ oder $\rho_1$ komplex, Fall 1,2,3,4,6,8)}\\\hline
$\lambda_1 \leq $\;\;\;\;& $\lambda_2 > \lambda_{2,neu}=$&  $\lambda_{2,alt}=$ & $C_1 \leq $ & $C_2 \leq $\\ \hline
0.36 & 1.69 & 0.903 & 0.0223 & 0.0152\\
0.38 & 1.69 & 0.887 & 0.0263 & 0.0181\\
0.40 & 1.69 & 0.871 & 0.0310 & 0.0214\\
0.42 & 1.69 & 0.856 & 0.0362 & 0.0252\\
0.44 & 1.67 & 0.842 & 0.0408 & 0.0287\\
0.46 & 1.59 & 0.829 & 0.0414 & 0.0297\\
0.48 & 1.52 & 0.816 & 0.0420 & 0.0307\\
0.50 & 1.45 & 0.803 & 0.0420 & 0.0315\\
0.52 & 1.39 & 0.791 & 0.0423 & 0.0324\\
0.54 & 1.31 & 0.780 & 0.0401 & 0.0317\\
0.56 & 1.23 & 0.769 & 0.0373 & 0.0305\\
0.58 & 1.13 & 0.759 & 0.0320 & 0.0274\\
0.60 & 1.04 & 0.749 & 0.0271 & 0.0245\\
0.62 & 0.96 & 0.739 & 0.0226 & 0.0216\\
0.64 & 0.88 & 0.730 & 0.0176 & 0.0182\\
0.66 & 0.82 & 0.721 & 0.0144 & 0.0156\\
0.68 & 0.76 & 0.712 & 0.0139 & 0.0126\\
0.70 & 0.71 & 0.704 & 0.0136 & 0.0099\\
\hline
\end{tabular}
\end{center}
\vspace*{0.5cm}

\noindent Beachte: Wir werden gleich zeigen, dass der Fall $\lambda_1 \leq 0.44$ gar nicht auftritt. Die Werte "`1.69"' wurden extra so gewählt - das ist kein Druckfehler.\\

\subsection{$\lambda_2$-Abschätzungen für die Fälle 5 und 7}
Fall 5 und 7 werden völlig analog zu Fall 2 behandelt, wobei wir diesmal etwas andere Parameter und andere zu beweisende $\lambda_2$-Abschätzungen wählen. Im Fall 5 wählen wir
\begin{eqnarray*}
  & &\gamma=0.76 +\frac{\lambda_{12}}{2}, \; k=0.84, \; \delta=0.0001 , \\
  & & k_{1,B}=0, \; k_{2,B}=\frac{1}{4}, \; k_{3,B}=0,  \\
  & &s_1=\lambda_2,\; s_{11}=\lambda_{2,alt}, \; s_{12}=\lambda_{2,neu}, \\
  & &s_2=\lambda_1, \; s_{21}=\lambda_{11}, \; s_{22}=\lambda_{12},  \\
  & &\Delta s_1=0.010,\; \Delta s_2=0.007, \; \Delta t=0.010,\; x_1=7.
\end{eqnarray*}
Wir bekommen in diesem Fall keine Verbesserung für $\lambda_1 \in [0.68,0.70]$ im Vergleich zum bereits bekannten Wert aus \cite[Table 10 (§9)]{Hea92}.\\

Im Fall 7 ist $\chi_1$ reell und $\rho_1$ komplex. Wir gehen genauso vor wie im Fall 2, beginnen wegen \cite[Lemma 9.5]{Hea92} unsere Werte bei $\lambda_1 \geq 0.50$ und gehen in $0.08$-Schritten vor. Wir wählen diesmal

\begin{eqnarray*}
  & &\gamma=0.61 +\frac{\lambda_{12}}{2}, \; k=0.81, \; \delta=0.0001, \\
  & & k_{1,A}=\frac{1}{4}, \; k_{2,A}=k, \; k_{3,A}=0,  \\
  & &s_1=\lambda_2,\; s_{11}=\lambda_{2,alt}, \; s_{12}=\lambda_{2,neu},  \\
  & &s_2=\lambda_1, \; s_{21}=\lambda_{11}, \; s_{22}=\lambda_{12},  \\
  & &\Delta s_1=\Delta s_2=\Delta t=0.015,\; x_1=7.
\end{eqnarray*}
Man beachte, dass im Fall $\lambda_1 \in [\lambda_{11},\lambda_{12}]$ mit $\lambda_{11}\geq 0.70$ wir aus \cite[Table 10 (§9)]{Hea92} keine Werte für $\lambda_{2,alt}$ ablesen können. In diesem Fall nehmen wir $\lambda_{2,alt}=\lambda_{11}$, was wegen $\lambda_1 \leq \lambda_2$ erlaubt ist. Es folgt

\vspace*{0.25cm}
\thinmuskip=0mu
\medmuskip=0mu
\thickmuskip=0mu
\begin{center}
\begin{tabular}[t]{|ll|ll|}
\multicolumn{4}{c}{\textbf{Tabelle 5. $\lambda_2$-Abschätzungen}}\\
\multicolumn{4}{c}{($\chi_1$ oder $\rho_1$ komplex, Fall 5)}\\\hline
$\lambda_1 \leq $ &  \begin{small}$\lambda_2 > \lambda_{2,neu}=$ \end{small}&  \begin{small}$\lambda_{2,alt}=$ \end{small} & $C_2 \leq$ \\ \hline
0.36 & 1.69 & 0.903 & 0.0664\\
0.38 & 1.69 & 0.887 & 0.0702\\
0.40 & 1.69 & 0.871 & 0.0742\\
0.42 & 1.69 & 0.856 & 0.0783\\
0.44 & 1.67 & 0.842 & 0.0799\\
0.46 & 1.56 & 0.829 & 0.0700\\
0.48 & 1.45 & 0.816 & 0.0606\\
0.50 & 1.36 & 0.803 & 0.0535\\
0.52 & 1.27 & 0.791 & 0.0465\\
0.54 & 1.19 & 0.780 & 0.0406\\
0.56 & 1.11 & 0.769 & 0.0348\\
0.58 & 1.04 & 0.759 & 0.0299\\
0.60 & 0.97 & 0.749 & 0.0249\\
0.62 & 0.91 & 0.739 & 0.0208\\
0.64 & 0.85 & 0.730 & 0.0167\\
0.66 & 0.79 & 0.721 & 0.0126\\
0.68 & 0.74 & 0.712 & 0.0092\\
\hline
\end{tabular}
\hspace*{0.2cm}
\begin{tabular}[t]{|ll|ll|}
\multicolumn{4}{c}{\textbf{Tabelle 6. $\lambda_2$-Abschätzungen}}\\
\multicolumn{4}{c}{($\chi_1$ reell und $\rho_1$ komplex, Fall 7)}\\\hline
$\lambda_1 \leq $ &   \begin{small}$\lambda_2 > \lambda_{2,neu}=$ \end{small}&  \begin{small}$\lambda_{2,alt}=$\end{small} & $C_1 \leq$ \\ \hline
0.54 & 1.43 & 0.780 & 0.0301\\
0.58 & 1.36 & 0.759 & 0.0276\\
0.62 & 1.28 & 0.739 & 0.0242\\
0.66 & 1.20 & 0.721 & 0.0206\\
0.70 & 1.11 & 0.704 & 0.0167\\
0.74 & 1.02 & 0.70 & 0.0128\\
0.78 & 0.93 & 0.74 & 0.0090\\
0.82 & 0.82 & 0.78 & 0.0070\\
\hline
\end{tabular}
\end{center}
\thinmuskip=3mu
\medmuskip=4mu plus 2mu minus 4mu
\thickmuskip=5mu plus 5mu
\vspace*{0.25cm}

Insgesamt erhalten wir die folgenden Abschätzungen für $\lambda_2$. Wähle dazu das Minimum der Einträge von Tabelle 4,5 und 6; dies ist letztendlich immer der Eintrag aus Tabelle 5. Zum Vergleich schreiben wir die alten Werte aus \cite[Table 10 (§9)]{Hea92} dazu.
\vspace*{0.25cm}

\begin{center}
\begin{tabular}[t]{|lc|c|}
\multicolumn{3}{c}{\textbf{Tabelle 7. $\lambda_2$-Abschätzungen}}\\
\multicolumn{3}{c}{($\chi_1$ oder $\rho_1$ komplex, alle Fälle)}\\\hline
& \begin{small}(neue Werte) \end{small}& \begin{small} (alte Werte)\end{small} \\
$\lambda_1 \leq $ &  $\lambda_2 >$ & $\lambda_2 >$ \\ \hline
0.36& 1.69& 0.903\\
0.38& 1.69& 0.887\\
0.40& 1.69& 0.871\\
0.42& 1.69& 0.856\\
0.44& 1.67& 0.842\\
0.46& 1.56& 0.829\\
0.48& 1.45& 0.816\\
0.50& 1.36& 0.803\\
0.52& 1.27& 0.791\\
0.54& 1.19& 0.780\\
0.56& 1.11& 0.769\\
0.58& 1.04& 0.759\\
0.60& 0.97& 0.749\\
0.62& 0.91& 0.739\\
0.64& 0.85& 0.730\\
0.66& 0.79& 0.721\\
0.68& 0.74& 0.712\\
0.70 & - & 0.704 \\
0.702 & - & 0.702\\
\hline
\end{tabular}
\end{center}
\vspace*{0.25cm}

\noindent Man bemerkt, dass die Güte der Verbesserung abnimmt, je mehr sich $\lambda_1$ dem Wert $0.70$ nähert. Dies stimmt mit der Beobachtung überein, dass in letzterem Fall die verwendete Ungleichung (\ref{l2-Hauptungleichung}) praktisch in diejenige aus \cite[Lemma 9.2]{Hea92} "`übergeht"'. Und aus letzterer wurden die Werte in \cite[Table 9 (§9)]{Hea92} gefolgert.

\section{$\lambda_3$-Abschätzungen}

\subsection{Für $\lambda_1 \in [0.52,0.62]$}\label{para9-l3bounds}
Sei $\chi_1$ oder $\rho_1$ komplex und $j \in \{2,3\}$. Ferner nehmen wir, dass $\lambda_1 \in [\lambda_{11},\lambda_{12}]$ für gewisse konkrete Konstanten $\lambda_{ij}$. Weiterhin sei $\lambda^{\star}$ das Minimum der $\lambda'$- und $\lambda_2$-Abschätzung aus Tabelle 2 und 7 unter der Voraussetzung $\lambda_1 \leq \lambda_{12}$. Wir leiten nun folgendermaßen Abschätzungen für $\lambda_3$ her:\\
Sind wir in Fall 7, so nehmen wir für $\lambda_j$ die $\lambda_2$-Abschätzung für Fall 7 aus Tabelle 6. \\ Für den Fall 1, sowie Fall 2,3,4,8, beweisen wir $\lambda_j$-Abschätzungen mittels Lemma \ref{para9-zentrales-neues-lemma}. Das Minimum der Abschätzungen für Fall 7, Fall 1 und Fall 2,3,4,8 gilt dann für $\lambda_j$, falls wir uns nicht in Fall 5 oder 6 befinden, also $\chi_1^3 \in \{\chi_j, \overline{\chi_j}\}, \; ord \; \chi_1\geq 6$. Letzteres kann aber nicht gleichzeitig für $j=2$ und $j=3$ erfüllt sein, da sonst $\chi_2 \in \{\chi_3,\overline{\chi_3}\}$, was ausgeschlossen ist. Also gilt das angesprochene Minimum für $\lambda_2$ oder $\lambda_3$ und damit wegen $\lambda_3 \geq \lambda_2$ für $\lambda_3$ allgemein.\\

Für den Beweis der Abschätzungen in Fall 1 und Fall 2 benutzen wir die gleichen Parameter $\gamma$ und $k$ wie oben bei der Herleitung von Tabelle 4, also (\ref{Wahl-des-gamma-Fall1-l2}). Außerdem benutzen wir auch die für Tabelle 4 berechneten Abschätzungen von $A$ und $B$. Die Abschätzungen der beiden Suprema $A$ und $B$ wurden dort sogar für einen größeren $\lambda^{\star}$-Bereich bewiesen, als wir es jetzt nötig haben (jetzt ist $\lambda^{\star}$ eine feste Konstante $c$, während bei den genannten Berechnungen $\lambda^{\star}=\lambda_2$ in einem speziellen Intervall liegen durfte, welches die Konstante $c$ enthält).\\
Diesmal sind die rechten Seiten der Ungleichungen im Lemma \ref{para9-zentrales-neues-lemma} (ohne die Suprema) offensichtlich monoton wachsend in $\lambda_1$ und $\lambda_j$. Wir brauchen also nicht analog zu oben vorzugehen mit einem $\delta>0$. Stattdessen erhält man sofort Abschätzungen für die beiden Fälle 1 und 2 analog zur Herleitung der $\lambda'$-Abschätzungen aus Tabelle 2. Schlussendlich benutzen wir dann die für Tabelle 4 bewiesene Tatsache, dass aus der jeweiligen Ungleichung für Fall 3,4 und 8 diejenige für Fall 2 folgt. Also gelten die in Fall 2 bewiesenen Abschätzungen auch für die Fälle 3,4 und 8.

Man bekommt die folgende Tabelle. Wir haben das Minimum der drei verschiedenen Abschätzungen - welches also die allgemein gültige Abschätzung für $\lambda_3$ darstellt - in eine extra Spalte geschrieben.

\begin{center}
 \begin{tabular}[t]{|lc|ccc|c|}
\multicolumn{6}{c}{\textbf{Tabelle 8. $\lambda_3$-Abschätzungen}}\\
\multicolumn{6}{c}{($\chi_1$ oder $\rho_1$ komplex)}\\\hline
 & (alle Fälle) & (Fall 1) & (Fall 2,3,4,8) & (Fall 7) & \\
$\lambda_1 \leq $\;\;\;\;\;\; & $\lambda_3 > $ & $\lambda_3 > $& $\lambda_3 >$  & $\lambda_3 >$ & $\lambda^{\star}$ \\ \hline
0.52& 1.320 & 1.352 & 1.320& 1.43 & 1.27\\
0.54& 1.243 & 1.253 & 1.243& 1.43 & 1.19\\
0.56& 1.160 & 1.160 & 1.167& 1.36 & 1.11\\
0.58& 1.079 & 1.079 & 1.103& 1.36 & 1.04\\
0.60& 1.001 & 1.001 & 1.038& 1.28 & 0.97\\
0.62& 0.933 & 0.933 & 0.979& 1.28 & 0.91\\
\hline
\end{tabular}
\end{center}
\vspace*{0.25cm}

\noindent Man beachte, dass die erste Zeile in Tabelle 8 ($\lambda_3 \geq 1.320$) zuerst nur unter der Voraussetzung $\lambda_1 \in [0.50,0.52]$ bewiesen wurde (dies wird vorausgesetzt bei der Abschätzung der Suprema $A$ und $B$). Tabelle 7 liefert schließlich im Fall $\lambda_1 \leq 0.50$ die Aussage $\lambda_3 \geq \lambda_2 \geq 1.36 \geq 1.320$.\\
Die Zeile $x$ in Tabelle 8 ($x \in \{3,4,5,6,7\}$) wurde analog zuerst nur für $\lambda_1 \in [0.46+0.02x,0.48+0.02x]$ bewiesen. Indem man die Abschätzung aus Zeile $(x-1)$ verwendet, erhält man dann die Aussage für $\lambda_1 \leq 0.48+0.02x$.

\subsection{Für $\lambda_1 \in [0.62,0.72]$ oder $\chi_1$ und $\rho_1$ beide reell}\label{para10analysis}
In \cite[§10]{Hea92} beweist Heath-Brown
\begin{lemma}[Lemma 10.3 aus \cite{Hea92}]\label{hb-lemma-10-3}
  Sei $\varepsilon>0$. Dann gibt es eine absolute Konstante $q_0(\varepsilon)$, so dass für $q \geq q_0(\varepsilon)$ gilt $$\lambda_3 \geq \frac{6}{7}-\varepsilon.$$
\end{lemma}
\noindent Durch analoges Vorgehen beweisen wir
\begin{lemma}\label{lemma-para10-analysis} Es gibt eine absolute Konstante $q_0$, so dass für $q \geq q_0$ in den Fällen $\chi_1$ komplex bzw. $\chi_1$ und $\rho_1$ beide reell die folgenden zwei Tabellen gelten.\\
\end{lemma}
 \begin{tabular}[t]{|lll|}
\multicolumn{3}{c}{\textbf{Tabelle 9. $\lambda_3$-Abschätzungen}}\\
\multicolumn{3}{c}{($\chi_1$ komplex)}\\\hline
$\lambda_1 \in $ & Zusatzbedingung & $\lambda_3 >$ \\ \hline
$[0.62,0.64]$ & - & 0.902\\
$[0.64,0.66]$ & - & 0.898\\
$[0.64,0.66]$ & $\lambda_2 \leq 0.86$ & 0.938\\
$[0.66,0.68]$ & - & 0.893\\
$[0.66,0.68]$ & $\lambda_2 \leq 0.83$ & 0.960\\
$[0.68,0.72]$ & - & 0.883\\
$[0.68,0.72]$ & $\lambda_2 \leq 0.81$ & 0.962\\
\hline
\end{tabular}
\hspace*{0.5cm}
\begin{tabular}[t]{|ll|}
\multicolumn{2}{c}{\textbf{Tabelle 10. $\lambda_3$-Abschätzungen}}\\
\multicolumn{2}{c}{($\chi_1$ und $\rho_1$ reell)}\\\hline
$\lambda_1 \in $ & $\lambda_3 >$ \\ \hline
$[0.44,0.60]$ & 1.175\\
$[0.60,0.68]$ & 1.078\\
$[0.68,0.78]$ & 0.971\\
\hline
\end{tabular}
\vspace*{0.25cm}\label{Tabelle-10}

\noindent Die zweite Zeile in Tabelle 9 besagt beispielsweise  ($\lambda_1 \in [0.64,0.66] \Rightarrow \lambda_3>0.898$), während die dritte Zeile aus der gleichen Tabelle besagt, dass ($\lambda_1 \in [0.64,0.66] \text{ und } \lambda_2 \leq 0.86 \Rightarrow \lambda_3 >0.938$).\\
Der Fall $\chi_1$ reell und $\rho_1$ komplex wird im Lemma nicht behandelt. Man hätte jedoch auf die gleiche Art und Weise die gleichen Werte beweisen können wie für den Fall $\chi_1$ komplex. Dies ist für uns nicht notwendig, deswegen lassen wir es aus.

Wir kommen jetzt zum Beweis des Lemmas. Der Beweis läuft im Wesentlichen genauso ab wie der Beweis der allgemeinen Abschätzung $\lambda_3 \geq \frac{6}{7} - \varepsilon$, die ohne Einschränkungen an $\lambda_1$ gilt. Wir führen den Beweis trotzdem vollständig durch und erwähnen auch ein paar Details, die in \cite{Hea92} nicht explizit aufgeschrieben werden.\\
Wir benutzen dafür die Definitionen für \label{Definition-Sigma-j} $\Sigma_2$ und $\Sigma_3$ und die Ungleichungen (10.2) und (10.4) aus \cite[S.310]{Hea92}.

\subsubsection{Beweis für $\chi_1$ komplex}

\textbf{Fall 1: }Es sei keiner der Charaktere in $\Sigma_3$ gleich $\chi_0$, dies gilt genau dann wenn $$\chi_3,\overline{\chi_3} \neq \chi_1 \chi_2, \chi_1 \overline{\chi_2}  .$$
Mit den beiden Lemmata \ref{hb-lemma-5-3} und \ref{hb-lemma-5-2-extended} folgt die Ungleichung \cite[(10.5)]{Hea92} und schließlich dann \cite[(10.6)]{Hea92}. Wir notieren letztere Ungleichung hier:
\begin{equation}\label{hb-ungleichung-10-6}0 \leq F(-\lambda_1)-F(\lambda_3-\lambda_1)-F(\lambda_2-\lambda_1)-F(0)+\frac{7}{6}f(0)+\varepsilon .
\end{equation}

\noindent Wir möchten nun beweisen, dass (\ref{hb-ungleichung-10-6}) immer gilt, wenn $\chi_1$ komplex und sagen wir $\lambda_1 \in [0.44,0.85]$ ist. Also müssen wir für die restlichen Fälle untersuchen, wie sich die Ungleichung verändert, wenn einer (oder mehrere) der Charaktere in $\Sigma_3$ doch gleich $\chi_0$ ist.\\
Der Fall $\chi_3 \in \{\chi_1 \chi_2, \chi_1 \overline{\chi_2} \}$ folgt nach Umbenennung aus dem Fall $\overline{\chi_3} \in \{\chi_1 \chi_2, \chi_1 \overline{\chi_2} \}$; vergleiche dafür die Bemerkung im Beweis von Lemma \ref{para9-zentrales-neues-lemma}. Der Fall $\overline{\chi_3}=\chi_1 \overline{\chi_2}$ folgt wiederum nach Umbenennung aus dem Fall $\overline{\chi_3}=\chi_1 \chi_2$, womit wir nur noch diesen letzten Fall betrachten müssen. Sei also $\chi_1 \chi_2 \chi_3 =\chi_0$. Wir unterscheiden die drei Fälle $\chi_2$ und $\chi_3$ komplex, $\chi_2$ komplex und $\chi_3$ reell, $\chi_2$ reell und $\chi_3$ komplex.
Da $\chi_1$ komplex ist, tritt der Fall, dass $\chi_2$ und $\chi_3$ beide reell sind, nicht auf, da sonst $\chi_0 (n)=\chi_1 \chi_2 \chi_3 (n)$ komplexe Werte annehmen würde. Wir weisen nochmal darauf hin, dass bei den nächsten drei Fällen stets $\chi_1 \chi_2 \chi_3=\chi_0$ angenommen wird. Wir beginnen mit\\

\noindent \textbf{Fall 2}: Seien $\chi_2$ und $\chi_3$ komplex.\\
Es ist $$\chi_1 \chi_2 \overline{\chi_3}= \chi_1 \chi_2 \chi_3 \overline{\chi_3}^2=\overline{\chi_3}^2.$$
Da nach Voraussetzung $\chi_3$ komplex ist, ist $\overline{\chi_3}^2\neq \chi_0$. Analog sind $\chi_1 \overline{\chi_2} \chi_3=\overline{\chi_2}^2$ und $\chi_1 \overline{\chi_2} \overline{\chi_3} = \chi_1^2$ ungleich $\chi_0$. Also muss man nur den Term $K(\beta_1 +i(\gamma_1+\gamma_2+\gamma_3),\chi_1\chi_2\chi_3)$ mit Lemma \ref{hb-lemma-5-3} behandeln anstelle von Lemma \ref{hb-lemma-5-2-extended}. Dann bekommt man anstatt von \cite[(10.5)]{Hea92} die Ungleichung
$$\Sigma_3 \leq \mathscr{L} Re  \{F(-\lambda_1 + i \mu)\} +  \frac{3}{2}f(0)\phi \mathscr{L} + \varepsilon \mathscr{L} ,$$
wobei $\mu=\mu_1+\mu_2+\mu_3$.\\
Wegen $\chi_2 \chi_3 =\overline{\chi_1}$ können wir beim Term $K(\beta_1+i(\gamma_2+\gamma_3),\chi_2 \chi_3)$ in $\Sigma_2$ die Nullstelle $\overline{\rho_1}$ verwenden. Anstelle von \cite[(10.4)]{Hea92} haben wir dann
$$\Sigma_2 \leq -\mathscr{L} Re  \{F(i \mu)\} + \frac{6}{2}f(0)\phi \mathscr{L} + \varepsilon \mathscr{L}. $$
Insgesamt erhalten wir aus \cite[(10.2)]{Hea92} wieder (\ref{hb-ungleichung-10-6}), aber mit dem zusätzlichen Term (benutze $\phi \leq \frac{1}{3}$)
$$+ Re  \left\{\frac{1}{4}F(-\lambda_1 + i \mu) - \frac{1}{2} F(i \mu)\right\} - \frac{1}{24}f(0)$$
auf der rechten Seite. Damit gilt (\ref{hb-ungleichung-10-6}), vorausgesetzt wir haben
\begin{equation}\label{para10-analysis-bedingung-1}\sup_{t \in \R} Re  \{F(-\lambda_1 + i t) - 2 F(i t)\}\leq \frac{1}{6}f(0).\end{equation}

\noindent \textbf{Fall 3 und 4:} Sei ein $\chi_j$ reell ($j \in \{2,3\}$) und der entsprechend andere Charakter komplex. Dann erhalten wir mit $\tilde{\mu_j}:=\mu -2 \mu_j$
$$\Sigma_3 \leq \mathscr{L} \Big{(}Re  \{F(-\lambda_1 + i \mu)\}+Re  \{F(-\lambda_1 + i \tilde{\mu_j})\} \Big{)} + \frac{2}{2}f(0)\phi\mathscr{L} + \varepsilon \mathscr{L}. $$ Wenn $j=3$ so gilt $\chi_2 \overline{\chi_3} =\chi_2 \chi_3= \overline{\chi_1}$, also
\begin{equation}\label{para10-analysis-Sigma-2}\Sigma_2 \leq -\mathscr{L} \Big{(}Re \{F(i\mu)\}+Re \{F(i \tilde{\mu_j})\} \Big{)} + \frac{6}{2}f(0)\phi\mathscr{L} + \varepsilon \mathscr{L}. \end{equation} Für $j=2$ haben wir $\chi_2 \overline{\chi_3}=\overline{\chi_2 \chi_3} =\chi_1$, also folgt mit (\ref{Im-z-Vorzeichen-in-Re-z-unwichtig}) wieder (\ref{para10-analysis-Sigma-2}). In Fall 3 und 4 bekommen wir damit wieder die Ungleichung (\ref{hb-ungleichung-10-6}), diesmal aber mit dem zusätzlichen Term
$$+ Re  \left\{\frac{1}{4}F(-\lambda_1 + i \mu) - \frac{1}{2} F(i \mu)\right\}+ Re  \left\{\frac{1}{4}F(-\lambda_1 + i \tilde{\mu_j}) - \frac{1}{2} F(i \tilde{\mu_j})\right\} - \frac{2}{24}f(0)$$
auf der rechten Seite. Es reicht wieder die Bestätigung von (\ref{para10-analysis-bedingung-1}), damit (\ref{hb-ungleichung-10-6}) in diesem Fall gilt.
Also beweisen wir (\ref{para10-analysis-bedingung-1}) und zwar für die Funktion $f(t)$ mit $\gamma=1.25$. Wir benutzen in §\ref{para-rechenkapitel}
\begin{eqnarray*}
  & &k_1=1, \;  k_2=0, \;  k_3=2,\\
  & &s_1=\lambda_1, \; s_{11}=0.44, \; s_{12}=0.85, \; (s_2=0,) \; \\
  & &\Delta s_1=0.03, \; \Delta s_2 =0, \; \Delta t=0.03 ,\; x_1=6
\end{eqnarray*}
und bekommen
$$\sup_{t \in \R} Re  \{F(-\lambda_1 + i t) - 2 F(i t)\} < 0.18 < 0.54\ldots=\frac{f(0)}{6} .$$\\
Wir folgern nun aus (\ref{hb-ungleichung-10-6}) die in Tabelle 9 behaupteten Werte. Für $f$ wurde bereits die Funktion mit $\gamma=1.25$ gewählt. Die rechte Seite von (\ref{hb-ungleichung-10-6}) ist monoton wachsend in $\lambda_2$ und $\lambda_3$. Eine Monotonie in $\lambda_1$ ist jedoch nicht so leicht festzustellen. Wir benutzen, dass $F(-\lambda_1)-F(\lambda_3-\lambda_1)$ monoton wachsend in $\lambda_1$ ist. Wenn also $\lambda_1 \in [\lambda_{11}, \lambda_{12}],$ $\lambda_2 \leq \lambda_{22}$ und $\lambda_3 \leq \lambda_{32}$, so folgt
$$-\varepsilon \leq F(-\lambda_{12})-F(\lambda_{32}-\lambda_{12})-F(\lambda_{22}-\lambda_{11})-F(0)+\frac{7}{6}f(0). $$
Ist $\lambda_1 \in [0.62,0.64]$, so berechnen wir die rechte Seite der letzten Ungleichung für $(\lambda_{11},\lambda_{12})=(0.62+j\delta, 0.62 + (j+1)\delta)$ mit $\delta=0.0001$ und $j=0,..,[(0.64-0.62)/\delta]$. Dabei nehmen wir $\lambda_{22}=\lambda_{32}=0.902$. Man beachte, dass wenn $\lambda_2>0.902$ wäre, wir wegen $\lambda_3 \geq \lambda_2$ sofort fertig wären. Eine Rechnung zeigt nun, dass für jegliche genannten $j$ wir etwas Negatives bekommen, womit die Aussage des Lemmas für $\lambda_1 \in [0.62,0.64]$ bewiesen wäre.\\ Analog und mit gleichem $\gamma,\delta$ beweisen wir alle Werte im Fall $\chi_1$ komplex. Die Zusatzbedingung $\lambda_2 \leq c$ baut man dabei auf die offensichtliche Weise $\lambda_{22}=c$ ein .\\

\subsubsection{Beweis für $\chi_1$ und $\rho_1$ beide reell}
Ausgangspunkt ist diesmal die Ungleichung \cite[(10.2)]{Hea92}, wobei wir dort $\beta_1$ durch $\beta_2$ ersetzen. Ziel ist es durch das übliche Vorgehen eine Ungleichung vom Typ (\ref{hb-ungleichung-10-6}) zu folgern. Dazu müssen wir alle Terme $K(\beta_2+it,\chi)$ in $\Sigma_2$ und $\Sigma_3$ mit $\chi \in \{\chi_0,\chi_1 \}$ gesondert behandeln.\\
Zuerst bemerkt man, dass alle Charaktere, die in $\Sigma_3$ auftauchen, ungleich $\chi_1$ sind. Auch sieht man schnell, dass alle Charaktere in $\Sigma_2$ ungleich $\chi_0$ und alle in $\Sigma_2$ außer $\chi_2 \chi_3$ und $\chi_2 \overline{\chi_3}$ ungleich $\chi_1$ sind.\\
Wir können außerdem annehmen, dass $\lambda_2 \leq 1.294$. Dann gilt nach \cite[Lemma 8.4]{Hea92}, dass $\lambda' \geq \lambda_2$. Damit kann man bei der Behandlung eines Terms $K(\beta_2+it,\chi_1)$ mittels Lemma \ref{hb-lemma-5-2-extended} jeweils $A_1 \cup A_2 =\{\rho_1\}$ setzen (da $\rho_1$ reell, fällt $\overline{\rho_1}$ weg).  Wir unterscheiden die folgenden Fälle:\\

\noindent \textbf{Fall 1: }Seien alle Charaktere in $\Sigma_3$ ungleich $\chi_0$.\\Dann sind aber auch $\chi_2 \chi_3$ und $\chi_2 \overline{\chi_3}$ ungleich $\chi_1$. Also sind in $\Sigma_2$ und $\Sigma_3$ alle Charaktere $\neq \chi_0,\chi_1$. Mit $\phi(\chi_1)=\frac{1}{4}$ und sonst $\phi(\chi)\leq \frac{1}{3}$ folgt
\begin{equation}\label{para10-ungleichung-fuer-chi1reell}
0 \leq F(-\lambda_2)-F(\lambda_3-\lambda_2)-F(0)-F(\lambda_1-\lambda_2) +\frac{9}{8}f(0)+\varepsilon .
\end{equation}

\noindent \textbf{Fall 2: }Nehme an, dass mindestens ein Charakter in $\Sigma_3$ gleich $\chi_0$ ist, nach Umbenennung können wir $\chi_1 \chi_2 \chi_3 =\chi_0$ annehmen. Es folgt $\chi_1 \overline{\chi_2} \overline{\chi_3} =\chi_0$, $\chi_2 \chi_3 =\chi_1$, $\chi_1 \chi_3=\overline{\chi_2}$ und $\chi_1 \overline{\chi_3} = \chi_2$.  Wir müssen zwei Fälle unterscheiden.\\
\noindent \textbf{Fall 2.1: }Seien $\chi_2$ und $\chi_3$ beide komplex. \\
Dann ist $\chi_1 \overline{\chi_2} \chi_3= \overline{\chi_2}^2 \neq \chi_0$ und analog $\chi_1 \chi_2 \overline{\chi_3}  \neq \chi_0$. Außerdem ist $\chi_2 \overline{\chi_3} \neq \chi_1$, da sonst aus $\chi_1 \chi_2 \chi_3=\chi_0$ ein Widerspruch zu $\chi_3$ komplex folgen würde. Weiter benutzt man die Nullstelle $\rho_2$ bzw. $\overline{\rho_2}$ von $L(s,\chi_2)$ bzw. $L(s,\overline{\chi_2})$. Zusammen mit den allgemeinen Bemerkungen für Fall 2, sowie der Voraussetzung $\mu_1=0$, erhalten wir dann (\ref{para10-ungleichung-fuer-chi1reell}) mit dem folgenden zusätzlichen Term auf der rechten Seite:
\begin{eqnarray}
  & &\sup_{t \in \R} \frac{1}{2} Re  \Big{\{}F(-\lambda_2 + it)-F(\lambda_1-\lambda_2 +it)-2F(it) \Big{\}} - (\frac{9}{8}-\frac{49}{48}) f(0) \label{para10-zusaetzlicher-Term-chi1reell}
\end{eqnarray}\\
\noindent \textbf{Fall 2.2: }Seien $\chi_2$ und $\chi_3$ beide reell.\\
Dann sind alle vier Charaktere in $\Sigma_3$ gleich $\chi_0$ und  $\chi_2 \overline{\chi_3}=\chi_1$. Diesmal verwenden wir für die beiden Charaktere $\chi_1 \chi_3, \;\chi_1 \overline{\chi_3}$, die gleich $\chi_2$ sind, beides Mal die Nullstelle $\rho_2$. Außerdem ist $\phi(\chi)=\frac{1}{4}$ für alle Charaktere $\chi \neq \chi_0$, die in \cite[(10.2)]{Hea92} auftauchen. Man bekommt (\ref{para10-ungleichung-fuer-chi1reell}) mit dem folgenden zusätzlichen Term auf der rechten Seite:
\begin{equation}\label{para10-zusaetzlicher-Term-2-chi1reell}\sup_{t \in \R}  Re  \Big{\{}F(-\lambda_2 + it)-F(\lambda_1-\lambda_2 +it)-F(it) \Big{\}} - (\frac{9}{8}-\frac{6}{8}) f(0).\end{equation}\\

\noindent Nun folgt mit (\ref{rechenkapitel-A-grosser-0}) und der Bedingung 2, die die Funktion $f$ erfüllt, dass (\ref{para10-zusaetzlicher-Term-chi1reell}) und (\ref{para10-zusaetzlicher-Term-2-chi1reell}) jeweils kleiner oder gleich
\begin{equation}\label{para10-zusaetzlicher-Term-3-chi1reell}
  \sup_{t \in \R}  Re  \Big{\{}F(-\lambda_2 + it)-F(\lambda_1-\lambda_2 +it)-F(it) \Big{\}} - \frac{5}{48}f(0)
\end{equation}
sind. Wenn wir also zeigen, dass der Wert von (\ref{para10-zusaetzlicher-Term-3-chi1reell}) $\leq 0$ ist, so gilt auch für Fall 2 die Ungleichung (\ref{para10-ungleichung-fuer-chi1reell}). Mit $\gamma=1.04$ und den Werten
\begin{eqnarray*}
  & &k_1=1, \; k_2=1, \; k_3=1,\\
  & &s_1=\lambda_2,\; s_{11}=0.44,\; s_{12}=1.175, \\
  & &s_2=\lambda_1, \; s_{21}=0.44,\; s_{22}=0.80, \\
  & &\Delta s_1=0.03, \; \Delta s_2 =0.03, \; \Delta t=0.03 , \; x_1=6
\end{eqnarray*}
folgt
 $$\sup_{t \in \R}  Re  \Big{\{}F(-\lambda_2 + it)-F(\lambda_1-\lambda_2 +it)-F(it) \Big{\}} < 0.10 < 0.13\ldots= \frac{5}{48}f(0) .$$

\noindent Damit gilt (\ref{para10-ungleichung-fuer-chi1reell}) immer im Fall $\chi_1$ und $\rho_1$ beide reell (für $\gamma=1.04$ wohlgemerkt). Den Rest beweisen wir jetzt völlig analog wie im Fall $\chi_1$ komplex. Zuerst erhalten wir aus
(\ref{para10-ungleichung-fuer-chi1reell})
$$-\varepsilon \leq F(-\lambda_{22})-F(\lambda_{32}-\lambda_{22})-F(0)-F(\lambda_{12}-\lambda_{21}) +\frac{9}{8}f(0) .$$
Diesmal müssen wir $\lambda_2$ anstelle von $\lambda_1$ staffeln. Für den ersten Eintrag in Tabelle 10 nimmt man die Werte
$\lambda_{11}=0.44, \; \lambda_{12}=0.60, \; \lambda_{32}=1.175$ und
$$[\lambda_{21},\lambda_{22}]=\big{[}0.44 + j \delta,0.44+ (j+1) \delta \big{]} \;\;\; (j \in \{0,1,\ldots,[\frac{\lambda_{32}-\lambda_{11}}{\delta}] \}, \; \delta=0.001).$$
Für alle $j$ erhält man dann etwas Negatives auf der rechten Seite der letzten Ungleichung, damit einen Widerspruch, und damit die Behauptung $\lambda_3 > \lambda_{32}=1.175$. Analog für die restlichen Einträge aus Tabelle 10.

\section{$\lambda_1$-Abschätzungen, Beweis von Theorem \ref{theorem-1-von-HB}}
Generalvoraussetzung dieses Abschnitts ist wieder, dass $\chi_1$ oder $\rho_1$ komplex ist. Unter dieser Voraussetzung beweist Heath-Brown (\cite[Lemma 9.5]{Hea92}) für hinreichend großes $q$
$$\lambda_1 \geq \left\lbrace \begin{array}{ll}
0.364& ord \; \chi_1 \geq 6, \\
0.397& ord \; \chi_1 = 5, \\
0.348& ord \; \chi_1 = 4, \\
0.389& ord \; \chi_1 = 3, \\
0.518& ord \; \chi_1 = 2. \\
\end{array} \right. $$

\noindent Wir verbessern diese Werte, indem wir einerseits die neuen $\lambda'$- und $\lambda_2$-Abschätzungen benutzen und andererseits das Verbesserungspotential Nr. 2 für die Fälle $ord \; \chi_1 \leq 5$.\\
Wir beginnen mit derselben Ungleichung wie Heath-Brown \cite[(9.16), S.308]{Hea92}:
\begin{eqnarray}
  0 &\leq& 14379 K(\beta,\chi_0) + 24480 K(\beta+i\gamma_1,\chi_1) + 14900 K(\beta+2i\gamma_1,\chi_1^2) \nonumber \\
  & &+ 6000 K(\beta+3i\gamma_1,\chi_1^3) + 1250 K(\beta+4i\gamma_1,\chi_1^4) . \label{l1-bounds-zentrale-ungleichung}
\end{eqnarray}
\vspace*{0.1cm}

\noindent Wir wählen $\beta=\beta^{\star}=1-\mathscr{L}^{-1} \lambda^{\star}$ mit einem $\lambda^{\star} \leq \min \{\lambda_2,\lambda'\}$ und benutzen die übliche Vorgehensweise mit den beiden Lemmata \ref{hb-lemma-5-3} und \ref{hb-lemma-5-2-extended}. Im Fall $ord \; \chi_1=2$ muss man dabei die beiden Nullstellen $\rho_1$ und $\overline{\rho_1}$ berücksichtigen, ansonsten nur die Nullstelle $\rho_1$. Es folgt

\begin{equation}\label{l1-bounds-endungleichung}0 \leq 14379 F(-\lambda^{\star}) - 24480 F(\lambda_1-\lambda^{\star}) + D + \varepsilon,\end{equation}
wobei
$$D=\left\lbrace \begin{array}{ll}
\frac{46630}{6}f(0)& ord \; \chi_1 \geq 6, \\
\\
\frac{46630}{8}f(0) + \sup\limits_{t \in \R} Re  \big{\{}-1250 F(\lambda_1-\lambda^{\star}+it) \big{\}} & ord \; \chi_1 = 5, \\
\\
\frac{45380}{8}f(0) + \sup\limits_{t \in \R} Re  \big{\{} 1250 F(-\lambda^{\star}+it)-6000 F(\lambda_1-\lambda^{\star}+it) \big{\}} &  ord \; \chi_1 =4, \\
\\
\frac{40630}{8}f(0) + \sup\limits_{t \in \R} Re  \big{\{} 6000 F(-\lambda^{\star}+it)-16150 F(\lambda_1-\lambda^{\star}+it) \big{\}} &  ord \; \chi_1 =3, \\
\\
\frac{30480}{8}f(0) + \sup\limits_{t \in \R} Re  \big{\{} 14900 F(-\lambda^{\star}+it)-30480 F(\lambda_1-\lambda^{\star}+it) \big{\}} & \\
 \;\;\;\;\;+ \sup\limits_{t \in \R} Re  \big{\{} 1250 F(-\lambda^{\star}+it)-6000 F(\lambda_1-\lambda^{\star}+it) \big{\}} &  ord \; \chi_1 =2 . \\
\end{array} \right.
$$

\noindent Die verschiedenen Suprema schätzen wir wie immer mit §\ref{para-rechenkapitel} ab. Dazu wählt man
$k_3=0$ und $k_1$, $k_2$ passend zum jeweiligen Supremum, sowie
\begin{eqnarray*}
  & &s_1=\lambda^{\star}=s_{11}=s_{12},\\
  & &s_2=\lambda_1,\; s_{21}=\lambda_{1,alt}, \; s_{22}=\lambda_{1,Ann}, \\
  & &\Delta s_1=0, \; \Delta s_2= 0.005, \; \Delta t=0.005, \; x_1=12 .
\end{eqnarray*}
Dabei ist $\lambda_{1,alt}$ die jeweils alte $\lambda_1$-Abschätzung aus \cite[Lemma 9.5]{Hea92}. Der Wert $\lambda_{1,Ann}$ ist irgendein konkreter Wert, der am Ende leicht oberhalb der bewiesenen $\lambda_1$-Abschätzung liegen wird. Wir nehmen nun an, dass $\lambda_1 \leq \lambda_{1,Ann}$ und wählen ein dazugehöriges $\lambda^{\star}$, das wir aus den Abschätzungen für $\lambda'$ und $\lambda_2$ in den Tabellen 2 und 7 ablesen. Für den Fall $ord \; \chi_1=2$ benutzen wir dabei die besseren Abschätzungen aus den Tabellen 3 und 6.\\
Wegen der Monotonie von (\ref{l1-bounds-endungleichung}) (ohne die Suprema) in $\lambda_1$ erhalten wir dann Abschätzungen $\lambda_1 > \lambda_{1,neu}$ auf die übliche Art und Weise. Außerdem müssen wir noch den Parameter $\gamma$ festlegen, den wir für den jeweiligen Fall benutzen möchten. In der folgenden Tabelle fassen wir für jeden der fünf Fälle alle relevanten Daten zusammen. Das wären einerseits $\gamma, \, \lambda^{\star},\, \lambda_{1,alt},\, \lambda_{1,Ann}$ und die Abschätzung $C$ des Supremums (bzw. der beiden Suprema zusammen im Fall $ord \; \chi_1=2$) und andererseits der bewiesene Wert $\lambda_{1,neu}$, der nach Annahme $\leq \lambda_{1,Ann}$ sein muss.

\vspace*{0.25cm}

\begin{center}

\begin{tabularx}{\textwidth}{|Xp{2.3cm}|p{1.2cm}p{1.2cm}p{1.2cm}p{1.2cm}p{1.2cm}|}
\multicolumn{7}{c}{\textbf{Tabelle 11. $\lambda_1$-Abschätzungen}}\\
\multicolumn{7}{c}{($\chi_1$ oder $\rho_1$ komplex)}\\
\hline
 $ord \text{ } \chi_1$ &
 $\lambda_1>\lambda_{1,neu}=$&
 $\lambda_{1,alt}$&
 $\lambda_{1,Ann}$&
 $\lambda^{\star}$ \newline&
 $\gamma$&
 C $\leq$\\
\hline
 $\geq 6$ & 0.440 & 0.364 & 0.44 & 1.67 & 1.00 &-\\ \hline
$= 5$ & 0.493 & 0.397 & 0.50 & 1.36 & 0.90 & 120\\ \hline
$= 4$ & 0.478 & 0.348 & 0.48 & 1.45 & 0.82 & 235\\ \hline
$= 3$ & 0.498 & 0.389 & 0.50 & 1.36 & 0.82 & 290\\ \hline
$= 2$ & 0.628 & 0.518 & 0.66 & 1.20 & 0.70 & 58\\ \hline
\end{tabularx}
\end{center}
\vspace*{0.25cm}
\noindent Die bewiesenen Werte halten wir im folgenden Lemma fest, welches nun \cite[Lemma 9.5]{Hea92} ersetzt:

\begin{lemma}\label{lemma-l1-abschaetzungen-werte} Sei $\chi_1$ oder $\rho_1$ komplex. Dann gibt es eine absolute Konstante $q_0$, so dass für $q\geq q_0$ gilt\\

$$\lambda_1 > \left\lbrace \begin{array}{ll}
0.440& ord \; \chi_1 \geq 6, \\
0.493& ord \; \chi_1 = 5, \\
0.478& ord \; \chi_1 = 4, \\
0.498& ord \; \chi_1 = 3, \\
0.628& ord \; \chi_1 = 2. \\
\end{array} \right. $$
\end{lemma}

\vspace*{0.25cm}
\begin{bemerkung}
Liu und Wang (1998, \cite[S.345-346]{LiuWan98}) benutzen die Ungleichung $$0 \leq (1+\cos \theta)(1+2 \cos \theta)^2$$ anstelle von derjenigen Ungleichung, aus der man (\ref{l1-bounds-zentrale-ungleichung}) schließt, nämlich \begin{equation}\label{l1-abs-anfangs-trig-ungleichung}0 \leq (3+10\cos \theta)^2(9+10\cos \theta)^2.\end{equation} Damit beweisen sie für den Fall $ord \; \chi_1=4$ den Wert $0.3711$ anstelle von $0.348$ (Heath-Brown) bzw. $0.478$ (Lemma \ref{lemma-l1-abschaetzungen-werte}).\\
\end{bemerkung}
\vspace*{0.25cm}
\begin{bemerkung}
Durch eine Verfeinerung der Rechnungen hätte man folgende minimale Verbesserungen zu dem eben Bewiesenen erzielen können:\\

\noindent Für alle Fälle: Eine Variation der Konstanten $3$ und $9$ in der Ungleichung (\ref{l1-abs-anfangs-trig-ungleichung}) bringt einen Gewinn von ca. $0.000-0.001$ bei den neuen $\lambda_1$-Abschätzungen.\\

\noindent Für den Fall $ord \; \chi_1 \geq 6$: Hätte man bessere $\lambda'$- und $\lambda_2$-Abschätzungen im Bereich $\lambda_1 \leq 0.46$ bewiesen (indem man z.B. für $\lambda_2$ acht verschiedene optimale Tabellen aufgestellt hätte), so hätte man eine Verbesserung von mindestens ca. $0.01$.\\

\noindent Für die Fälle $ord \; \chi_1 \leq 5$: Einerseits hätte man spezielle, bessere $\lambda'$- und $\lambda_2$-Abschätzungen beweisen können (wegen $\phi(\chi_1)=\frac{1}{4}$). Andererseits hätte man die bewiesenen Abschätzungen, z.B. $\lambda_1 > 0.493$ bei $ord \; \chi_1=5$, dazu benutzen können das Supremum nochmal abzuschätzen, diesmal unter der Voraussetzung $\lambda_1 \geq \lambda_{1,alt}= 0.493$. Dies kann einen kleineren Wert für die Abschätzung des Supremums liefern. Beide eben erwähnten Verfeinerungen liefern Verbesserungen von zusammen wenigen Hundertstel.\\

\noindent Man beachte, dass in unserer Arbeit keine der eben genannten Verfeinerungen auch nur einen minimalen Gewinn für die letztendlich zulässige Linniksche Konstante $L$ liefern würde.
\end{bemerkung}
\vspace*{0.25cm}
\begin{bemerkung}
Aus (\ref{L-s-chi0-zeta-s}) folgt, dass $L(s,\chi_0)$ keine Nullstellen in der Region
\begin{equation}\label{nullstellenfreie-region}
\sigma \geq 1- \frac{0.440}{\mathscr{L}}, \; |t|\leq 1
\end{equation}
hat. Lemma \ref{lemma-l1-abschaetzungen-werte} liefert also das Theorem \ref{theorem-1-von-HB} im Fall $\chi_1$ oder $\rho_1$ komplex. Für den Fall $\chi_1$ und $\rho_1$ beide reell folgt das Theorem aus \cite[Lemma 8.4, Lemma 8.8]{Hea92}. Damit ist Theorem \ref{theorem-1-von-HB} bewiesen.
\end{bemerkung}

\section{Abschätzungen der Nullstellendichte}\label{para11-neu}
Seien $q \in \N$ und $\lambda>0$ fest. Dann definiere
$$N(\lambda)=\#\{\chi \; mod \; q \, |\, \chi \neq \chi_0, \, L(s,\chi) \text{ hat eine Nullstelle in } \sigma \geq 1-\mathscr{L}^{-1}\lambda , \; |t| \leq 1 \}. $$
 So wie in \cite[S.316]{Hea92} nummerieren wir die $N(\lambda)$ verschiedenen Charaktere $\chi \neq \chi_0$, die im betrachteten Bereich eine Nullstelle haben, durch, nämlich $\chi^{(1)},\chi^{(2)},\ldots,\chi^{(N(\lambda))}$. Zu jedem Charakter wählen wir eine zugehörige Nullstelle $\chi^{(k)}$ mit maximalem Realteil, also \label{Definition-rho-j-oben}
$$Re \{\rho^{(k)}\} = \max \{Re \{\rho\} \,|\, L(\rho,\chi^{(k)})=0, \; Re \{\rho\} \geq 1-\mathscr{L}^{-1}\lambda , \; |Im\{\rho\}| \leq 1\} .$$
Setze
$$\rho^{(k)}=\beta^{(k)}+i\gamma^{(k)}, \; \beta^{(k)}=1-\mathscr{L}^{-1}\lambda^{(k)}, \; \gamma^{(k)}=\mathscr{L}^{-1}\mu^{(k)}. $$

\noindent Heath-Brown beweist nun für $ 0 < \lambda \leq \frac{1}{3}\log \log \mathscr{L}, \; \varepsilon>0$ und $q \geq q_0(\varepsilon)$, dass (\cite[(11.5)]{Hea92})
\begin{equation}\label{para11-abschaetzung-1-N-l}N(\lambda)\leq (1+\varepsilon)\frac{67}{6\lambda}(e^{\frac{73\lambda}{30}}-e^{\frac{16\lambda}{15}}). \end{equation}
Im Grunde genommen zeigt er noch mehr, nämlich für $\lambda_0=\frac{1}{3} \log \log \mathscr{L}$, $\varepsilon$, $c_1$, $c_2>0$, $\phi=\max\limits_{\chi \; mod \; q}\phi(\chi)$ und $q \geq q_0(\varepsilon,c_1,c_2)$, dass (\cite[(11.3)]{Hea92})
\begin{equation}\label{para11-abschaetzung-2}
  \sum_{1\leq k\leq N(\lambda_0)} \frac{\lambda^{(k)}}{e^{(4\phi+6c_1+2c_2)\lambda^{(k)}}-e^{(2\phi + 4c_1)\lambda^{(k)}}}\leq \frac{2\phi + 2c_1 +c_2}{4c_1c_2}+\varepsilon .
\end{equation}
Die Ungleichung (\ref{para11-abschaetzung-1-N-l}) folgt aus (\ref{para11-abschaetzung-2}) mit $\phi \leq \frac{1}{3}$, $c_1=\frac{1}{10},$ und $c_2=\frac{1}{4}$. Dafür muss man beachten, dass $\lambda/(e^{A\lambda}-e^{B\lambda})$ monoton fallend ist\footnote{Das zeigt man wohl am schnellsten mittels $(e^{Ax}-e^{Bx})/x=\int_B^A e^{t x} \,dt$.} für $A>B>0$.\\

\noindent In Verbesserungspotential Nr. 7 \cite[S.336-337]{Hea92} beschreibt Heath-Brown eine Variation in der Herleitung der letzten Ungleichung, welche eine Verbesserung von (\ref{para11-abschaetzung-1-N-l}) liefert. Für den Beweis von Theorem \ref{haupttheorem} ist es jedoch von großem Vorteil, wenn man eine Verbesserung von (\ref{para11-abschaetzung-2}) hätte. \\
Eine sehr simple Methode, um eine kleine Verbesserung von (\ref{para11-abschaetzung-2}) zu erzielen, ist so vorzugehen wie Heath-Brown in \cite[S.336-337]{Hea92}, jedoch mit dem einzigen Unterschied, dass man einen gewissen freien Parameter $w_{\chi}$ anders wählt. Dies führt zum folgenden Lemma. Der Beweis dieses Lemmas ist Inhalt des gesamten restlichen Paragraphen \ref{para11-neu}.

\begin{lemma}\label{para11-lemma1}
   Seien $\delta, c_1, c_2>0$ und $\lambda_0=\frac{1}{3} \log \log \mathscr{L}$. Weiterhin seien $u=\frac{1}{3}+2c_1$, $v=\frac{1}{3}+2c_1+c_2$ und $x=\frac{2}{3} +3c_1+c_2$. Schließlich sei $w_1:\R_{\geq 0} \rightarrow \R$ eine stetige Funktion, die auf $[0,v)$ und $(v,\infty)$ stetig differenzierbar ist. Für $t \geq 0$ gelte dabei  $1\ll w_1(t) \ll 1$  und $w_1'(t) \ll 1 $ mit absoluten impliziten Konstanten. Dann gilt für $q \geq q_0(\delta,c_1,c_2,w_1)$
  \begin{equation}\label{para11-zentrale-ungleichung}\sum_{1\leq k\leq N(\lambda_0)} \left(\int_u^x w_1(t)^2 e^{2 \lambda^{(k)} t} \,dt \right)^{-1} \leq \frac{1+ \delta}{c_1 c_2^2} \int_u^x w_1(t)^{-2} \min \{t-u,v-u\}\, dt .\end{equation}
\end{lemma}
\noindent Man bemerke, dass das Aussehen der Funktion $w_1(t)$ außerhalb von $t\in [u,x]$ keine Rolle für die Aussage des Lemmas spielt. Setzt man $w_1(t) \equiv 1$, so folgt aus diesem Lemma die Aussage (\ref{para11-abschaetzung-2}). Als Funktion $w_1(t)$ wird später die Funktion in (\ref{para15-definition-w-1}) recht optimal sein.

Wir kommen jetzt zum Beweis des Lemmas. Dieser läuft wortwörtlich so ab, wie der Beweis von (\ref{para11-abschaetzung-2}), mit dem einzigen Unterschied, dass an einer Stelle der frei wählbare Parameter $w_1(t)$ eingefügt wird und die folgenden Rechnungen dann entsprechend anzupassen sind. Wir werden den Beweis ab der Stelle, wo der neue Parameter eingefügt wird, beginnen. Im Folgenden sei dazu immer $q\geq q_0$, wobei $q_0$ eine absolute Konstante ist, die von allen im Lemma benutzten Parametern abhängen darf. Weiterhin sei $\lambda_0$ und $0<u<v<x$ so definiert wie im Lemma. Definiere zusätzlich $w:=c_1-\delta$ und setze $U=q^u,\, V=q^v,\, W=q^w, \, X=q^x$. Die gewählten Werte für $w$, $v$ und $x$ erfüllen die Bedingung \cite[(11.10)]{Hea92}. Damit gilt die Ungleichung \cite[Seite 318 letzte Zeile]{Hea92}:
\begin{equation}\label{para11-ungleichung-am-anfang}
  \sum_{\chi} w_{\chi} \leq (1+O(\mathscr{L}^{-1})) \sum_{\chi} \Big{|}\sum_{n=1}^{\infty} w_{\chi} (\sum_{d|n} \psi_d )(\sum_{d|n} \theta_d ) \chi(n) n^{-\rho(\chi)} (e^{-n/X}-e^{-n\mathscr{L}^2/U} )\Big{|}^2 .
\end{equation}
Die Summation läuft dabei über die $N(\lambda_0)$ verschiedenen Charaktere $\chi=\chi^{(k)}$, die weiter oben definiert wurden. Außerdem ist $\rho(\chi)=\rho^{(k)}$ für $\chi=\chi^{(k)}$ und $w_{\chi} \geq 0$ sind später zu wählende Gewichte. Wir müssen auch die Definitionen von $\theta_d$ und $\psi_d$ nachreichen, nämlich (\cite[(11.6)]{Hea92} und \cite[(11.7)]{Hea92})

\begin{equation}\label{para11-definition-psi-d}
  \psi_d = \left\lbrace \begin{array}{ll}
\mu(d) &  1 \leq d \leq U, \\
\mu(d) \frac{\log (V/d)}{\log (V/U)} & U \leq d \leq V, \\
0 &  V \leq d , \\
\end{array} \right.
\end{equation}
\begin{equation}\label{para11-definition-theta-d}
  \theta_d = \left\lbrace \begin{array}{ll}
\mu(d) \frac{\log (W/d)}{\log W} &  1 \leq d \leq W, \\
0 & W \leq d .\\
\end{array} \right.
\end{equation}
Dabei ist $\mu(d)$ die Möbiussche $\mu$-Funktion. Wir setzen nun
\begin{equation}\label{vb7-definition-a-n}a_{n \chi}=w_1(\mathscr{L}^{-1} \log n)w_{\chi}^{\frac{1}{2}} (\sum_{d|n} \theta_d ) \chi(n) n^{\frac{1}{2}-\rho(\chi)} (e^{-n/X}-e^{-n\mathscr{L}^2/U})^{\frac{1}{2}}\end{equation}
und
\begin{equation}\label{vb7-definition-b-n}b_n=w_1(\mathscr{L}^{-1} \log n)^{-1} (\sum_{d|n} \psi_d )n^{-\frac{1}{2}}(e^{-n/X}-e^{-n\mathscr{L}^2/U})^{\frac{1}{2}}, \end{equation}
haben also $$\sum_{\chi} w_{\chi} \leq (1+O(\mathscr{L}^{-1})) \sum_{\chi} \Big{|}\sum_{n=1}^{\infty} a_{n\chi} b_n \Big{|}^2 .$$
Wir weisen daraufhin, dass in \cite[§11]{Hea92} die $a_{n\chi}$ und $b_n$ ohne das $w_1(t)$ definiert wurden. In Verbesserungspotential Nr. 7 wird die Einführung dieses neuen Parameters vorgeschlagen.\\
Zuerst haben wir für alle komplexe Folgen $(C_{\chi})=(C_{\chi^{(k)}})_k$, dass
\begin{eqnarray}
 \sum_{n=1}^{\infty} \left|\sum_{\chi} a_{n\chi}C_{\chi} \right|^2&=&\sum_{n=1}^{\infty} (\sum_{\chi} a_{n\chi}C_{\chi})\overline{(\sum_{\chi} a_{n\chi}C_{\chi})} \nonumber \\
 &=&\sum_{\chi,\chi'} \sum_{n=1}^{\infty} a_{n\chi} \overline{a_{n\chi'}} C_{\chi}\overline{C_{\chi'}} \nonumber \\
 &=&\sum_{\chi \neq \chi'} C_{\chi}\overline{C_{\chi'}} \sum_{n=1}^{\infty} a_{n\chi} \overline{a_{n\chi'}}+\sum_{\chi} |C_{\chi}|^2 \sum_{n=1}^{\infty} |a_{n\chi}|^2 . \label{para11-zwischensumme}
\end{eqnarray}
Nun gilt für $\chi \neq \chi'$ (dies beweisen wir später)
\begin{equation}\label{para11-nichtdiagonalen-Terme-klein}
  \sum_{n=1}^{\infty} a_{n\chi} \overline{a_{n\chi'}}\ll w_{\chi}^{\frac{1}{2}} w_{\chi'}^{\frac{1}{2}} \mathscr{L}^{-1}.
\end{equation}
Also folgt
\begin{eqnarray*}
\sum_{\chi \neq \chi'} C_{\chi}\overline{C_{\chi'}} \sum_{n=1}^{\infty} a_{n\chi} \overline{a_{n\chi'}} &\ll& \mathscr{L}^{-1} \sum_{\chi \neq \chi'} C_{\chi}\overline{C_{\chi'}} w_{\chi}^{\frac{1}{2}} w_{\chi'}^{\frac{1}{2}} \\
&\ll& \mathscr{L}^{-1} \left(\sum_{\chi} |C_{\chi}| w_{\chi}^{\frac{1}{2}}\right)^2\\
&\ll& \mathscr{L}^{-1} (\sum_{\chi} |C_{\chi}|^2)(\sum_{\chi} w_{\chi}) ,
\end{eqnarray*}
wobei wir zuletzt die Cauchy Schwarz Ungleichung verwendeten. Sei $\varepsilon>0$ und setze $Re \{\rho(\chi) \}=1-\mathscr{L}^{-1}\lambda(\chi)$. Für $q \geq q_0=q_0(\varepsilon)$ gilt dann (das wird auch später bewiesen)
\begin{equation}\label{para11-zweite-ungleichung}
  \sum_{n=1}^{\infty} |a_{n\chi}|^2 \leq w_{\chi} \frac{1+\varepsilon}{w} \int_u^x w_1(t)^2 e^{2 \lambda(\chi)t} \, dt.
\end{equation}
Wählen wir also \begin{equation}\label{para11-definition-von-w-chi}w_{\chi}= \left(w^{-1} \int_u^x w_1(t)^2 e^{2 \lambda(\chi)t} \,dt \right)^{-1}\end{equation}
so ist die zweite Summe in (\ref{para11-zwischensumme})
$$\leq (1+\varepsilon)\sum_{\chi} |C_{\chi}|^2 .$$
Setzen wir nun $\mathscr{B} = 1+ \varepsilon + \varepsilon \sum\limits_{\chi} w_{\chi},$
so folgt für jede Wahl von $(C_{\chi})$, dass
\begin{equation}\label{para11-dualitaetsprinzip-anfang}\sum_{n=1}^{\infty} \left|\sum_{\chi} a_{n\chi}C_{\chi} \right|^2 \leq \mathscr{B} \sum_{\chi} |C_{\chi}|^2 . \end{equation}

\noindent Nun zitieren wir die dritte und letzte Ungleichung, deren Beweis wir auf das Ende des Abschnitts verschieben, nämlich
\begin{equation}\label{para11-erste-ungleichung}
  \sum_{n=1}^{\infty} |b_n|^2 \leq \frac{1+\varepsilon}{(v-u)^2} \int_u^x w_1(t)^{-2} \min\{t-u,v-u\} \,dt.
\end{equation}
Nach dem Dualitätsprinzip, welches ein schnelles Korollar der Cauchy Schwarz Ungleichung ist, siehe z.B. \cite[S.170]{IwaKow04}, folgt insgesamt für hinreichend kleines $\varepsilon>0$, dass
\begin{eqnarray*}
  \sum_{\chi} w_{\chi} &\leq& (1+\varepsilon) \sum_{\chi} \Big{|}\sum_{n=1}^{\infty} a_{n\chi} b_n \Big{|}^2 \\
  &\leq& (1+\varepsilon)\mathscr{B} \sum_{n=1}^{\infty} |b_n|^2 \\
  &\leq& (1+ 3\varepsilon+ 2\varepsilon \sum_{\chi} w_{\chi}) \frac{1}{(v-u)^2} \int_u^x w_1(t)^{-2} \min\{t-u,v-u\} \,dt .
\end{eqnarray*}
Nach Anpassung des $\varepsilon$ folgt die Aussage des Lemmas.\\

\noindent Damit verbleibt noch der Beweis von (\ref{para11-nichtdiagonalen-Terme-klein}), (\ref{para11-zweite-ungleichung}) und (\ref{para11-erste-ungleichung}). Heath-Brown erwähnt deren Richtigkeit in \cite[S.335 unten]{Hea92}, gibt jedoch keine Beweise an. Dies ist aber auch nicht verwunderlich. Die Ungleichungen folgen auf völlig analoge Weise wie im Fall $w_1(t) \equiv 1$, der ja in \cite[§11]{Hea92} behandelt wird. (\ref{para11-zweite-ungleichung}) und (\ref{para11-erste-ungleichung}) folgen dabei "`straight-forward"' (wenngleich ein wenig mühselig) durch partielle Summation und (\ref{para11-nichtdiagonalen-Terme-klein}) wird mit Hilfe des Cauchyschen Integralsatzes bewiesen. Wir geben die drei Beweise in den folgenden Abschnitten an. Wir weisen noch darauf hin, dass jegliche implizite Konstanten in den folgenden drei Unterabschnitten von allen im Lemma benutzten festen Parametern abhängen dürfen, aber (natürlich) nicht von $q$.\\

\subsection{Beweis von (\ref{para11-nichtdiagonalen-Terme-klein})}
Um einen rigorosen Beweis durchzuführen, gehen wir zurück zu \cite[§2]{Hea92}. Dort benutzt Heath-Brown Abschätzungen der Charaktersummen durch Burgess, um eine Abschätzung für $$L(s,\chi)=\sum_{n=1}^{\infty} \chi (n) n^{-s}$$ zu beweisen. Wir brauchen hier eine analoge Abschätzung für $$\sum_{1\leq n\leq T} \chi (n) n^{-s}$$ für beliebiges $T\geq 1$. Die Vorgehensweise dafür ist wortwörtlich die gleiche. Wir benutzen das in \cite[Lemma 2.1]{Hea92} zitierte Resultat von Burgess für $k=3$, also

\begin{lemma}[Burgess, 1986, \cite{Bur86}]\label{lemma-burgess}
 Seien $q\geq 1$ und $\chi \neq \chi_0$ ein Charakter $mod \; q$. Ferner seien $M\geq 1$ und $1 \leq H \leq q$. Für alle $\varepsilon>0$ gilt dann
 $$\sum_{M <n \leq M+H} \chi(n) \ll_{\varepsilon} q^{\frac{1}{9}+\varepsilon} H^{\frac{2}{3}} . $$
\end{lemma}

\noindent Aus diesem Lemma folgern wir
\begin{lemma}[vergleiche Lemma 2.5 aus \cite{Hea92}] \label{L-s-chi-abschaetzung}
Sei $\chi \neq \chi_0$ ein Charakter $mod \; q$. Ferner sei $T\geq 1$ oder $T=\infty$ und $k \geq 3$. Außerdem sei $s=\sigma+it \in \C$ mit $Re  \{s\} = 1-\frac{1}{k}$. Dann gilt
  $$\sum_{1\leq n \leq T} \frac{\chi(n)}{n^{s}} \ll_{\varepsilon} (1+|t|) q^{\frac{1}{3k}+\varepsilon}. $$
\end{lemma}
\begin{beweis}
Der Beweis läuft praktisch identisch zu dem von \cite[Lemma 2.5]{Hea92} ab. Sei zuerst $T\leq q^{\frac{1}{3}}$. Es folgt dann
$$\sum_{1\leq n \leq T} \frac{\chi(n)}{n^s}\leq \sum_{1\leq n \leq q^{\frac{1}{3}}} \frac{1}{n^{1-\frac{1}{k}}}
\leq q^{\frac{1}{3k}}  \sum_{1\leq n \leq q^{\frac{1}{3}}} \frac{1}{n} \ll_{\varepsilon} q^{\frac{1}{3k}+\varepsilon}. $$
Ist $q^{\frac{1}{3}} < T \leq q$, so haben wir mit partieller Summation und Lemma \ref{lemma-burgess} folgenden zusätzlichen Term:
\begin{eqnarray*}
  \sum_{q^{\frac{1}{3}} < n \leq T} \frac{\chi(n)}{n^s} &\ll& \Big{|}\sum_{q^{\frac{1}{3}} < n \leq T} \chi(n) \Big{|} T^{-1+\frac{1}{k}} + (1+|t|) \int_{q^{\frac{1}{3}}}^T x^{-2+\frac{1}{k}} \Big{|} \sum_{q^{\frac{1}{3}}<n\leq x} \chi(n)\Big{|} \,dx \\
  &\ll_{\varepsilon}& q^{\frac{1}{9}+\varepsilon}T^{\frac{2}{3}} T^{-1+\frac{1}{k}} + (1+|t|)\int_{q^{\frac{1}{3}}}^T x^{-2+\frac{1}{k}+\frac{2}{3}} q^{\frac{1}{9}+\varepsilon} \, dx \\
  &\ll& (1+|t|)q^{\frac{1}{3k}+2 \varepsilon}.
\end{eqnarray*}
Ist schließlich $T>q$, so schätzen wir den zusätzlichen Term durch partielle Summation mit der Pólya-Vinogradov Ungleichung ($\sum\limits_{y\leq n \leq x}\chi(n) \ll_{\varepsilon} q^{\frac{1}{2}+\varepsilon}$) ab zu
$$ \sum_{q<n \leq T} \frac{\chi(n)}{n^s} \ll_{\varepsilon} (1+|t|) q^{-\frac{1}{2}+\frac{1}{k}+\varepsilon} \ll (1+|t|) .$$
\end{beweis}

Wir führen nun den Beweis der Ungleichung (\ref{para11-nichtdiagonalen-Terme-klein}). Dieser läuft analog zum Fall $w_1(t) \equiv 1$, welcher in
\cite[S.318, dritte Zeile von unten]{Hea92} "`geführt"' wird. Geführt in Anführungsstrichen, da für den Beweis gesagt wird, dass dieser analog zur Behandlung der Gleichung (11.8) in \cite{Hea92} geschieht, was auch absolut zutreffend ist.\\
Seien also $\chi \neq \chi'$ zwei Charaktere aus der betrachteten Menge mit zugehöriger Nullstelle $\rho$ bzw. $\rho'$, also sind $$\rho,\rho' \in \{s \in \C\,| \, \sigma \geq 1- \mathscr{L}\lambda_0 , \; |t| \leq 1 \}.$$

\noindent Wir ändern die Funktion $w_1(t)$ so ab, dass $w_1'(t)=0$ für $t \geq 2x$. Dabei soll $w_1(t)$ weiterhin die Bedingungen aus dem Lemma erfüllen. Dies ist möglich und hat keinen Einfluss auf die Aussage des Lemmas. Betrachte die folgende offensichtlich konvergente Summe:
\begin{eqnarray*}
  S&:=&w_{\chi}^{-\frac{1}{2}} w_{\chi'}^{-\frac{1}{2}} \sum_{n=1}^{\infty} a_{n\chi} \overline{a_{n\chi'}} \nonumber \\
  &=&\sum_{n=1}^{\infty} (\sum_{d|n} \theta_d)^2 \chi \overline{\chi'} (n) n^{1-\rho-\overline{\rho'}}(e^{-n/X}-e^{-n \mathscr{L}^2 /U}) w_1(\mathscr{L}^{-1} \log n)^2 . \label{para11-nondiagonal-term}
\end{eqnarray*}

\noindent Dabei sei $\Gamma(s)$ die Gamma-Funktion. Für $y >0$ gilt (siehe z.B. \cite[S.220]{Bru95}) $$e^{-\frac{1}{y}}=\frac{1}{2\pi i}\int_{1-i\infty}^{1+i\infty} \Gamma(s) y^s \,ds .$$
Wir vertauschen nun Summe und Integral, was mit dem Satz der majorisierten Konvergenz erlaubt ist. Als Majorante kann wegen $|\theta_d|,|w_1(t)|,|\chi(n)| \ll 1$ und $\sum\limits_{d|n} 1 \ll n^{\varepsilon}$ die Funktion $C X |\Gamma(s)|$ genommen werden für ein $C>0$. Es folgt
\begin{equation}\label{para11-nondiagonal-term-integral}S=\frac{1}{2 \pi i} \int_{1-i\infty}^{1+i\infty} \Gamma(s) \left(X^s-\left(\frac{U}{\mathscr{L}^2}\right)^s\right) G(s) \, ds, \end{equation}
wobei $G(s)$ wie folgt definiert ist ($s':=s+\rho+\overline{\rho'}-1$):
\begin{eqnarray}
  G(s)&=&\sum_{n=1}^{\infty} \frac{ (\sum_{d|n} \theta_d)^2 \chi \overline{\chi'} (n) w_1(\mathscr{L}^{-1} \log n)^2}{n^{s'}} \nonumber \\
  &=&\sum_{d_1 \leq W, d_2 \leq W} \theta_{d_1} \theta_{d_2} \sum_{n=1, \atop d_1 |n, d_2|n}^{\infty} \frac{w_1(\mathscr{L}^{-1} \log n)^2 \chi \overline{\chi'} (n)}{n^{s'}} \nonumber \\
  &=&\sum_{d_1 \leq W, d_2 \leq W} \frac{\theta_{d_1} \theta_{d_2} \chi \overline{\chi'} (kgV(d_1,d_2))}{kgV(d_1,d_2)^{s'}} \sum_{j=1}^{\infty} \frac{w_1(\mathscr{L}^{-1} \log \{j \; kgV(d_1,d_2)\})^2 \chi \overline{\chi'} (j)}{j^{s'}} . \nonumber \\
  & & \label{para11-H(s)}
\end{eqnarray}
Die Umordnung der Summe war wegen absoluter Konvergenz für jedes einzelne $s$ erlaubt. Außerdem haben wir $\theta_d=0$ für $d>W$ und die Äquivalenz
$$d_1|n \text{ und } d_2|n \Longleftrightarrow kgV(d_1,d_2)|n$$
benutzt.\\
Da $w_1(t)$ ab einer Stelle konstant ist und $\sum\limits_{n=1}^{\infty} \chi \overline{\chi'}(n) n^{-s}$ holomorph in $Re  \{s\}>0$, folgt, dass auch $G(s)$ in $Re \{s'\}>0$ bzw. $Re \{s\}>1-Re \{\rho\}-Re \{\overline{\rho'}\}$ holomorph ist.\\

\noindent Sei $k \geq 3$ fest gewählt. Wir verschieben die Integrationslinie in (\ref{para11-nondiagonal-term-integral}) auf $Re  \{s\} =2-\frac{1}{k}-Re \{\rho\}-Re \{\overline{\rho'}\}$. Die horizontalen Integrale verschwinden (beachte, dass $q$ fest ist), da $\Gamma(s)$ für großen Imaginärteil gleichmäßig exponentiell gegen $0$ geht, während $L(s,\chi)$ nur polynomiell gegen $\infty$ geht
(\cite[Lemma 2.6.2]{Bru95}). Da $$X^s-\left(\frac{U}{\mathscr{L}^2}\right)^s$$ die Nullstelle $s=0$ hat und $\Gamma(s)$ bei $s=0$ einen einfachen Pol, ist der Integrand im betrachteten Bereich holomorph und wir bekommen mit dem Cauchyschen Integralsatz, dass (\ref{para11-nondiagonal-term-integral}) gilt mit der neuen Integrationslinie $(2-\frac{1}{k}-Re \{\rho\}-Re \{\overline{\rho'}\})$ anstelle von $(1)$. Dieses neue Integral schätzen wir nun ab. Zuerst gilt $\Gamma(s) \ll e^{-|t|}$ und $(X^s-(\frac{U}{\mathscr{L}^2})^s) \ll U^{-\frac{1}{k} + \varepsilon}$.
Die innere Summe in (\ref{para11-H(s)}) schätzen wir mit partieller Summation ab zu
\begin{eqnarray*}
  \sum_{j=1}^{\infty} \frac{w_1(\mathscr{L}^{-1} \log \{j \; kgV(d_1,d_2)\})^2 \chi \overline{\chi'} (j)}{j^{s'}} &\ll& L(s', \chi \overline{\chi'})+\int_1^{2X} \frac{\Big{|}\sum\limits_{1\leq j \leq x} \chi \overline{\chi'}(j) j^{-s'}\Big{|}}{\mathscr{L} y } \,dy \\
  &\ll& \sup_{y\geq 1} \Big{|}\sum\limits_{1 \leq j \leq y} \chi \overline{\chi'} (j) j^{-s'}\Big{|}.
\end{eqnarray*}
Insgesamt haben wir dann mit Lemma \ref{L-s-chi-abschaetzung} und $\varepsilon=\frac{1}{k^2}$
$$G(s) \ll (1+|t|) q^{\frac{1}{3k}+\frac{2}{k^2}} \sum_{d_1 \leq W, d_2 \leq W} kgV(d_1,d_2)^{-1+\frac{1}{k}}.$$
Nun ist
 $$\sum_{d_1 \leq W, d_2 \leq W} kgV(d_1,d_2)^{-1+\frac{1}{k}} \ll W^{\frac{2}{k}} \sum_{n \leq W^2} d(n)^2 n^{-1}$$
Dabei ist $d(n)=\sum\limits_{j|n} 1$ und wir verwendeten $$\sum_{d_1 \leq W,d_2 \leq W,\atop kgV(d_1,d_2)=n}1 \leq \sum_{d_1,\,d_2 \in \N , \atop d_1|n,\, d_2 | n} 1 \leq d(n)^2 .$$
Benutzt man jetzt $\sum\limits_{n\leq x} d(n)^2 \ll x \log^3 x \;$ (\cite[S.62]{Mur08}) dann folgt mit partieller Summation
$$ W^{\frac{2}{k}} \sum_{n \leq W^2} d(n)^2 n^{-1} \ll W^{\frac{2}{k}} q^{\frac{1}{k^2}}. $$

\noindent Insgesamt haben wir
$$S \ll (q^{2w+\frac{1}{3}-u})^{\frac{1}{k}} q^{\frac{3}{k^2}} $$
und da $u>2w+\frac{1}{3}$, erhalten wir für hinreichend großes $k$, dass
$$ S\ll \mathscr{L}^{-1} ,$$
was zu zeigen war.

\subsection{Beweis von (\ref{para11-zweite-ungleichung})}
Wir beweisen
$$\sum_{n=1}^{\infty} w_1(\mathscr{L}^{-1} \log n)^2 (\sum_{d|n} \theta(d))^2 \chi_0(n) n^{1-2Re \{\rho(\chi)\}} (e^{-n/X}-e^{-n\mathscr{L}^2/U})
$$
\begin{equation}
  \label{para11-zweite-ungleichung-ausgeschrieben}
  \leq  \frac{1+\varepsilon}{w} \int_u^x w_1(t)^2 e^{2 \lambda(\chi)t} \,dt .
\end{equation}
\noindent Wenn man das $\chi_0(n)$ in der Summe weglassen würde, dann könnte man entlang den nachfolgenden Zeilen sogar eine Gleichheit bis auf einen vernachlässigbaren Fehler beweisen.\\

\noindent Wir benötigen die folgende Abschätzung von Graham (\cite{Gra78}, wird auch in \cite[(11.13)]{Hea92} zitiert)
\begin{equation}\label{para11-abschaetzung-psi-d}
  \sum_{n \leq N} (\sum_{d|n} \psi_d)^2 = \left\lbrace \begin{array}{ll}
1 & 1 \leq N \leq U, \\
\frac{N \log (N/U)}{\log^2 (V/U)} +O(\frac{N}{\log^2 (V/U)}) &  U \leq N \leq V, \\
\frac{N }{\log (V/U)} +O(\frac{N}{\log^2 (V/U)}) &  V \leq N .\\
\end{array} \right.
\end{equation}
Setzt man $U=1$ und $V=W$, so wird daraus
\begin{equation}\label{para11-abschaetzung-theta-d}
  \sum_{n \leq N} (\sum_{d|n} \theta_d)^2 = \left\lbrace \begin{array}{ll}
\frac{N \log N}{\log^2 W} +O(\frac{N}{\log^2 W}) & 1 \leq N \leq W, \\
\frac{N }{\log W} +O(\frac{N}{\log^2 W}) &  W \leq N. \\
\end{array} \right.
\end{equation}\\

Wir teilen die Summe auf der linken Seite von (\ref{para11-zweite-ungleichung-ausgeschrieben}) in die Teilsummen $S_1$, $S_2$, $S_3$, $S_4$, $S_5$, $S_6$ auf, entsprechend den Intervallen
$$n \in [1,W],\; (W,U\mathscr{L}^{-3}],\;(U\mathscr{L}^{-3},U],\;(U,X],\;(X,X\mathscr{L}],\;(X\mathscr{L},\infty).$$ Wir erinnern außerdem daran, dass $ Re\{\rho(\chi) \} \geq 1-  (3\mathscr{L})^{-1} \log \log \mathscr{L}$. Also gilt für $n \leq q^y \; (y>0)$, dass $$n^{2-2\rho(\chi)} \leq \log^{\frac{2y}{3}} \mathscr{L}.$$ Ferner benutzen wir in den folgenden Abschätzungen die Gleichung (\ref{para11-abschaetzung-theta-d}), sowie dass $0<2w<u<v<x$, $1 \ll w_1(t) \ll 1$ und $w_1'(t) \ll 1$ gilt.\\

\noindent Es ist $\sum\limits_{d|n} \theta_d \ll n^{\varepsilon}$. Zusammen mit $e^x=1+x+O(x^2)$ für $|x|\leq 1$, folgt dass
$$S_1 \ll  W^{1+3\varepsilon} \frac{W\mathscr{L}^2}{U} \ll \mathscr{L}^{-1} $$
für $\varepsilon$ hinreichend klein.\\

\noindent Für $n \in (W,U\mathscr{L}^{-3}]$ gilt $e^{-n/X}-e^{-n\mathscr{L}^2/U} \ll \mathscr{L}^{-1}. $
Also folgt mit partieller Summation, dass
$$S_2 \ll \mathscr{L}^{-1} \log^u \mathscr{L}  \sum_{W<n\leq U} (\sum_{d|n} \theta_d)^2 n^{-1} \ll \mathscr{L}^{-1} \log^u \mathscr{L}  .$$\\

\noindent Mit partieller Summation folgt weiterhin
$$S_3 \ll \log^u \mathscr{L}  \sum_{U \mathscr{L}^{-3} <n\leq U} (\sum_{d|n} \theta_d)^2 n^{-1} \ll \mathscr{L}^{-1} \log^{u+1} \mathscr{L} $$
und
$$S_5 \ll \log^x \mathscr{L} \sum_{X<n\leq X \mathscr{L}} (\sum_{d|n} \theta_d)^2 n^{-1} \ll  \mathscr{L}^{-1} \log^{x+1} \mathscr{L} .$$\\

\noindent Für $S_6$ benutzen wir $n^{1-2Re\{\rho(\chi)\}} \sum\limits_{d|n} \theta_d \ll n^{-1+\varepsilon}$ und erhalten für hinreichend kleines $\varepsilon>0$
$$S_6 \ll (X \mathscr{L})^{-1+\varepsilon} \sum_{n > X \mathscr{L}} e^{-n/X} \ll \mathscr{L}^{-1} .  $$\\

\noindent Damit verbleibt der letzte Term $S_4$, welcher auch der einzige mit einem wirklichen Beitrag ist. Wir können durch partielle Summation alle differenzierbaren Terme aus der Summe rausnehmen und den Rest mit Graham's Abschätzung (\ref{para11-abschaetzung-theta-d}) behandeln. Setze $\lambda=\lambda(\chi)(>0)$:
\begin{eqnarray*}
 S_3&\leq& \sum_{U<n\leq X}(\sum_{d|n} \theta_d)^2 w_1(\mathscr{L}^{-1} \log n)^2 n^{-1+ 2\lambda \mathscr{L}^{-1}} \\
&=& O\left(\mathscr{L}^{-1} \log^x \mathscr{L}\right)   \\
& & - \frac{1}{w \mathscr{L}} \int_U^X (t-U) \Big{(}w_1(\mathscr{L}^{-1} \log t)^2 (-1+ 2\lambda \mathscr{L}^{-1}) t^{-2+2\lambda \mathscr{L}^{-1}} \\
& &+ 2 w_1(\mathscr{L}^{-1} \log t ) w_1'(\mathscr{L}^{-1} \log t ) \mathscr{L}^{-1} t^{-1} t^{-1+ 2\lambda \mathscr{L}^{-1}}  \Big{)}\, dt\\
\\
&\leq&O(\mathscr{L}^{-1} \log^x \mathscr{L}) \\
& &+ \frac{1+O(\mathscr{L}^{-1} \log \log \mathscr{L})}{w \mathscr{L}} \int_U^X w_1(\mathscr{L}^{-1} \log t)^2  t^{-1+2\lambda \mathscr{L}^{-1}}\, dt .
\end{eqnarray*}
\noindent Nach Substitution $s=\mathscr{L}^{-1} \log t$ folgt
$$S_3 \leq  \frac{1+\varepsilon}{w \mathscr{L}} \mathscr{L} \int_u^x w_1(s)^2 e^{2 \lambda s} \, ds $$
und wir sind fertig.

\subsection{Beweis von (\ref{para11-erste-ungleichung})}
Wir beweisen
$$\sum_{n=1}^{\infty} w_1(\mathscr{L}^{-1} \log n)^{-2} (\sum_{d|n} \psi_d)^2 n^{-1} (e^{-n/X}-e^{-n \mathscr{L}^2/U}) $$
$$\leq \frac{1+\varepsilon}{(v-u)^2} \int_u^x w_1(t)^{-2} \min\{t-u,v-u\} \,dt .$$

\noindent Auch hier könnten wir entlang der nachfolgenden Zeilen das $\leq$ durch ein $=$ ersetzen auf Kosten eines vernachlässigbaren Fehlers.\\
Für die partiellen Summationen benutzen wir diesmal (\ref{para11-abschaetzung-psi-d}) und teilen die Summe in die Intervalle $$n \in [1,U],\; (U,V],\; (V,X],\; (X,X\mathscr{L}],\; (X\mathscr{L},\infty)$$ auf. Die entsprechenden Summen seien der Reihe nach $S_1$, $S_2$, $S_3$, $S_4$, $S_5$.\\

\noindent Für $d\leq U$ ist $\psi_d=\mu(d)$. Also gilt für $n\leq U$, dass

$$(\sum_{d|n} \psi_d)^2= (\sum_{d|n} \mu(d))^2 = \left\lbrace \begin{array}{ll}1 & n=1,\\ 0 & n>1  \end{array} \right.$$
und damit
$$S_1 \ll  \mathscr{L}^{-1} . $$

\noindent Analog zur Vorgehensweise im letzten Unterabschnitt bekommen wir $$S_4 \ll \mathscr{L}^{-1} \log \mathscr{L} \;\; \text{ und }\;\; S_5 \ll \mathscr{L}^{-1} .$$

\noindent Weiterhin gilt
\begin{eqnarray*}
  S_2 &\leq& \sum_{U<n\leq V} (\sum_{d|n} \psi_d)^2 w_1(\mathscr{L}^{-1} \log n)^{-2} n^{-1} \\
  &=& \frac{1+\varepsilon}{(v-u)^2 \mathscr{L}^2} \int_U^V \frac{ w_1(\mathscr{L}^{-1} \log t)^{-2} \log \frac{t}{U}}{t} \,dt +O(\mathscr{L}^{-1}).
\end{eqnarray*}
Die Substitution $s=\mathscr{L}^{-1} \log t$ liefert
$$S_2\leq \frac{1+\varepsilon}{(v-u)^2} \int_u^v w_1(s)^{-2} (s-u) \,ds .$$

\noindent Schließlich haben wir
\begin{eqnarray*}
  S_3 &\leq& \sum_{V <n \leq X} (\sum_{d|n} \psi_d)^2 w_1(\mathscr{L}^{-1}\log n)^{-2} n^{-1} \\
  &=& \frac{1}{(v-u)\mathscr{L}} \int_V^X w_1(\mathscr{L}^{-1}\log t)^{-2} t^{-1} \, dt + O(\mathscr{L}^{-1})\\
  &\leq&\frac{1+\varepsilon}{(v-u)} \int_v^x w_1(s)^{-2}  \, ds.
\end{eqnarray*}

\noindent Damit ist der Beweis komplett.

\section{Abschätzungen der Nullstellendichte - kleine $\lambda$}\label{para12-neu}
Wir benutzen die gleichen Notationen wie zu Beginn von §\ref{para11-neu}. Heath-Brown beweist  in \cite[§12]{Hea92} Abschätzungen für $N(\lambda)$, die für kleine $\lambda$ (ungefähr $\lambda \leq 1.30$) besser sind als das, was \cite[§11]{Hea92} bzw. §\ref{para11-neu}  liefert.\\
In den Kommentaren zu Verbesserungspotential Nr. 9 \cite[S.336-337]{Hea92} wird gezeigt, dass im Argument von \cite[§12]{Hea92} mehr Potential steckt, als letztendlich benutzt wird. Dabei geht es darum, eine gewichtete Version von \cite[Lemma 12.1]{Hea92} zu finden und zwar in dem Stil wie (\ref{para11-abschaetzung-2}) eine gewichtete Version von (\ref{para11-abschaetzung-1-N-l}) ist. Eine sinnvolle gewichtete Version haben wir nun nicht finden können. Jedoch führen bereits die Ausführungen in \cite[S.336-337]{Hea92} zu Verbesserungen für die letztendlich zulässige Linniksche Konstante $L$. Speziell erhält man das folgende Lemma, welches \cite[Lemma 12.1]{Hea92} verallgemeinert (für $N_0=0$ bekommt man letzteres Lemma zurück).

\begin{lemma}[vergleiche S.336 in \cite{Hea92}]\label{para12-lemma-neu}
Bezeichnungen wie zu Beginn von §\ref{para11-neu}. Seien $\varepsilon>0, \; \lambda \in (0,2)$ und $\lambda_{11}>0$ beliebig gewählt. Es gelte $\lambda_{11} \leq \lambda_1 \leq \lambda$. Ferner sei $f$ eine Funktion, die Bedingung 1 und 2 erfüllt und für die gilt
$$F(\lambda-\lambda_{11})>\frac{f(0)}{6}.$$
Schließlich seien $\lambda_0 \geq 0$ und $N_0 \in \N_0$ beliebig gewählt mit $\lambda_0 \leq \lambda$ und $N_0 \leq 10000$.\\
Angenommen es ist $N(\lambda_0) \geq N_0$. Dann gilt für $q \geq q_0(\varepsilon,f)$ und mit der Abkürzung $N=N(\lambda)$ die Ungleichung
$$\left((N-N_0)F(\lambda-\lambda_{11})+N_0 F(\lambda_0-\lambda_{11})-N \frac{f(0)}{6}\right)^2$$
  \begin{equation}\label{para12-ungleichung-neu}\leq F(-\lambda_{11})\left(N F(-\lambda_{11})+(N^2-N)\frac{f(0)}{6}\right) +\varepsilon. \end{equation}
\end{lemma}
\begin{beweis}
  Heath-Brown erwähnt diese Aussage ohne Beweis (und mit einem Druckfehler) in \cite[S.336]{Hea92}. Der Beweis folgt analog zu demjenigen von \cite[Lemma 12.1]{Hea92}. Wir geben ihn der Vollständigkeit halber an. Seien also die Voraussetzungen des Lemmas erfüllt und $\rho^{(1)}, \rho^{(2)},\ldots, \rho^{(N)}$ die zu den $N=N(\lambda)$ verschiedenen Charakteren gehörigen Nullstellen. Man beachte, dass $\rho^{(k)}$ ($j=1,\ldots,N$) in der Menge $R$ liegt. \\
  Aus Lemma \ref{hb-lemma-5-2-extended} folgt nun mit $\beta_{11}=1-\mathscr{L}^{-1}\lambda_{11}$
  \begin{eqnarray*}
    K(\beta_{11} +i \gamma^{(j)}, \chi^{(j)}) \leq - \mathscr{L} F(\lambda^{(j)}-\lambda_{11}) + \mathscr{L} (\frac{f(0)}{6}+\varepsilon) .
  \end{eqnarray*}
  Unter den $N$ Charakteren gibt es $N_0$ Stück, für die $\lambda^{(j)} \leq \lambda_0$ gilt. Für die restlichen $N-N_0$ Charaktere benutzen wir $\lambda^{(j)} \leq \lambda$. Aus Monotoniegründen folgt also
  \begin{eqnarray*}
    \mathscr{L} \Big{(}(N-N_0) F(\lambda-\lambda_{11})&+&N_0 F(\lambda_0-\lambda_{11})-N (\frac{f(0)}{6}+\varepsilon)\Big{)}\\
    &\leq& \mathscr{L} (\sum_{j \leq N} F(\lambda^{(j)}-\lambda_{11})-N (\frac{f(0)}{6}+\varepsilon)) \\
    &\leq& - \sum_{j\leq N} K(\beta_{11} +i \gamma^{(j)},\chi^{(j)})\\
    &=& - \sum_{n=1}^{\infty} \Lambda(n) n^{-\beta_{11}} f(\mathscr{L}^{-1}\log n) Re \left\{\sum_{j \leq N} \frac{\chi^{(j)}(n)}{n^{i\gamma^{(j)}}}\right\} \\
    &\leq& \sum_{n=1}^{\infty} \Lambda(n) n^{-\beta_{11}} f(\mathscr{L}^{-1}\log n) \left|\sum_{j \leq N} \frac{\chi^{(j)}(n)}{n^{i\gamma^{(j)}}}\right| .
  \end{eqnarray*}

\noindent Nach Voraussetzung ist
\begin{eqnarray*}
  0&<&N (F(\lambda-\lambda_{11})- (\frac{f(0)}{6}+\varepsilon)) \\
   &\leq & (N-N_0) F(\lambda-\lambda_{11}) + N_0 F(\lambda_0-\lambda_{11})-N (\frac{f(0)}{6}+\varepsilon),
\end{eqnarray*}
wobei wir für die erste Ungleichung anmerken müssen, dass $N(\lambda) \leq N(2) $ beschränkt ist gemäß (\ref{para11-abschaetzung-1-N-l}) und $\varepsilon>0$ hinreichend klein gewählt wurde (beachte: wenn das Lemma für hinreichend kleines $\varepsilon$ gilt, dann automatisch auch für größeres $\varepsilon$). Also bleibt die vorherige Ungleichungskette nach Quadrieren richtig und wir erhalten mittels der Cauchy Schwarz Ungleichung
$$\mathscr{L}^2 \left((N-N_0) F(\lambda-\lambda_{11}) + N_0 F(\lambda_0-\lambda_{11})-N (\frac{f(0)}{6}+\varepsilon)\right)^2 \leq \Sigma_1 \Sigma_2, $$
wobei
\begin{eqnarray*}\Sigma_1=\sum_{n=1}^{\infty} \Lambda(n)  \chi_0(n) n^{-\beta_{11}} f(\mathscr{L}^{-1}\log n) = K(\beta_{11},\chi_0)
\end{eqnarray*}
und
\begin{eqnarray*}
\Sigma_2&=& \sum_{n=1}^{\infty} \Lambda(n) n^{-\beta_{11}} f(\mathscr{L}^{-1}\log n) \left|\sum_{j \leq N} \frac{\chi^{(j)}(n)}{n^{i\gamma^{(j)}}}\right|^2 \\
&=& \sum_{n=1}^{\infty} \Lambda(n) n^{-\beta_{11}} f(\mathscr{L}^{-1}\log n) \sum_{j,k \leq N} \frac{\chi^{(j)}(n)}{n^{i\gamma^{(j)}}} \frac{\overline{\chi^{(k)}}(n)}{n^{-i\gamma^{(k)}}} \\
&=&\sum_{j,k \leq N} K(\beta_{11} +i(\gamma^{(j)}-\gamma^{(k)}),\chi^{(j)}\overline{\chi^{(k)}}) .\end{eqnarray*}
Für die $N$ Terme in $\Sigma_2$ mit $j=k$ gilt nach Lemma \ref{hb-lemma-5-3} $$K(\beta_{11},\chi_0)\leq \mathscr{L}(F(-\lambda_{11})+\varepsilon),$$ während für die $N^2-N$ Terme mit $j\neq k$ nach Lemma \ref{hb-lemma-5-2-extended} $$K(\beta_{11} +i(\gamma^{(j)}-\gamma^{(k)}),\chi^{(j)}\overline{\chi^{(k)}}) \leq \mathscr{L}(\frac{f(0)}{6}+\varepsilon)$$ gilt. Man beachte, dass wegen der Voraussetzungen die Terme $F(z)$, $N$, $N_0$ alle beschränkt sind. Ersetzen wir also in den letzten Ungleichungen das $\varepsilon$ durch einen hinreichend kleinen positiven Wert, so folgt das Lemma.
\end{beweis}

\begin{bemerkung}
Weiter oben sprachen wir das Problem an, eine gewichtete Version von \cite[Lemma 12.1]{Hea92} zu finden. In diese Richtung liefert die Argumentation im Beweis die Ungleichung
  $$\sum_{j \leq N} F(\lambda^{(j)}-\lambda_{11}) \leq \sqrt{ F(-\lambda_{11})\left(N F(-\lambda_{11})+(N^2-N)\frac{f(0)}{6}\right)} + N \frac{f(0)}{6}+\varepsilon . $$
  Nun kann man für $N$ die obere Abschätzung, welche aus dem Lemma folgt (setze $N_0=0$) nehmen und damit die rechte Seite der letzten Ungleichung explizit (in Abhängigkeit von $\lambda_{11}$ und $f$) abschätzen. Jedoch bekamen wir aus der dadurch gewonnenen gewichteten Version keine Verbesserung für die Linniksche Konstante $L$.
\end{bemerkung}
\vspace*{0.5cm}
\noindent Man erinnere sich daran, dass wir an oberen Abschätzungen für $N(\lambda)$ interessiert sind. Angenommen, die Parameter $\lambda, \; \lambda_{11}, \; \lambda_0$ und $N_0$ sind fest gewählt. Dann wählen wir noch $\varepsilon=10^{-7}$ und $f$ mittels $$\gamma=0.975+0.525\lambda-0.550 \lambda_{11}-0.014 N_0 \cdot (\lambda-\lambda_0).$$
Letztere Wahl für $\gamma$ erweist sich als näherungsweise optimal. Aus (\ref{para12-ungleichung-neu}) folgt nun
\begin{equation}\label{para12lemma-deduction}0 \leq h(N), \end{equation}
wobei $h(N)$ eine gewisse Parabel in $N$ darstellt. Diese ist nach unten geöffnet, falls
$$\left(F(\lambda-\lambda_{11})-\frac{f(0)}{6}\right)^2 > F(-\lambda_{11})\frac{f(0)}{6} . $$
Wenn letzteres zutrifft und $h(N)$ zwei reelle einfache Nullstellen hat, sagen wir $h_1<h_2$ (Ermittlung überlassen wir dem Computer), dann folgt aus (\ref{para12lemma-deduction}) die Abschätzung
$$N=N(\lambda) \leq [h_2].$$

\noindent Die auf diese Weise erhaltenen Werte notieren wir in den beiden nachfolgenden Tabellen. Diese stellen teilweise eine Erweiterung von \cite[Table 13 (§12)]{Hea92} dar. In der obersten Zeile stehen dabei die Voraussetzungen - also die Werte für $\lambda_{11}$, $\lambda_0$ und $N_0$ - unter denen die Werte in der entsprechenden Spalte bewiesen wurden.
Zur genaueren Erklärung betrachte man beispielsweise die sechste Spalte der nachfolgenden Tabelle 12. In der obersten Zeile steht die Generalvoraussetzung $\lambda_1 \geq 0.60$, also setzt man zur Herleitung der entsprechenden Werte $\lambda_{11}=0.60$. Unterhalb der Generalvoraussetzung steht als erstes "`keine Vor."', was für "`keine Zusatzvoraussetzung"' steht. Das bedeutet, dass man $\lambda_0=0, \; N_0=0$ wählt und daraus allgemein gültige Abschätzungen für $N(\lambda)$  ($\lambda \in \{0.875, 0.900, \ldots 1.300 \}$) herleitet. Man erhält
$$N(\lambda)\leq \left\lbrace \begin{array}{ll}
9 & \lambda=0.875,\\
9 & \lambda=0.900,\\
\ldots & \\
89 & \lambda=1.275,\\
140 & \lambda=1.300.
\end{array} \right. $$
Unterhalb von "`keine Vor."' steht die Voraussetzung "`$N(1.075) \geq 7$"'. Das bedeutet man setzt im Lemma $\lambda_0=1.075$ und $N_0=7$. Dann bekommt man für $\lambda > \lambda_0$ Verbesserungen der Abschätzungen, nämlich
$$N(\lambda)\leq \left\lbrace \begin{array}{ll}
21 & \lambda=1.100,\\
23 & \lambda=1.125,\\
\ldots & \\
66 & \lambda=1.275,\\
101 & \lambda=1.300.
\end{array} \right. $$
Analog für "`$N(1.075) \geq 11$"' und "`$N(1.075) \geq 16$"'. \label{Tabelle-12}

\thinmuskip=0mu
\medmuskip=0mu
\thickmuskip=0mu
\begin{center}
\begin{small}
\begin{tabularx}{14.2cm}{|p{1cm}|p{1.25cm}p{0.85cm}p{1.3cm}p{1.5cm}p{2.3cm}X|}
\multicolumn{7}{c}{\begin{small}\textbf{Tabelle 12. Abschätzungen für $N(\lambda)$ unter verschiedenen Voraussetzungen}\end{small}}\\
\hline
  &
 $\lambda_1 \geq 0.40$\newline \begin{tiny}keine Vor.\end{tiny}&
 $ \geq 0.44$\newline \begin{tiny}keine Vor.\end{tiny}&
 $ \geq 0.54$\newline \begin{tiny}keine Vor.,\newline $N(1.125)\geq 9$\end{tiny}&
 $ \geq 0.58$\newline \begin{tiny}keine Vor.,\newline $N(1.100)\geq 11$\end{tiny}&
 $ \geq 0.60$\newline \begin{tiny}keine Vor.,\newline $N(1.075)\geq 7$, \newline $N(1.075)\geq 11$, \newline $N(1.075)\geq 16$ \end{tiny}&
 $ \geq 0.62$\newline \begin{tiny}keine Vor.,\newline $N(1.075)\geq 7$, \newline $N(1.075)\geq 10$, \newline $N(1.075)\geq 12$, \newline $N(1.075)\geq 15$ \end{tiny}\\ \hline
 $\lambda =$ & $N(\lambda)\geq$&$N(\lambda)\geq$&$N(\lambda)\geq$&$N(\lambda)\geq$&$N(\lambda)\geq$&$N(\lambda)\geq$\\
 & & & & & & \\
 0.875 &10  &10  &9       &9        &9               &9                  \\
 0.900 &11  &11  &10      &10       &9               &9                  \\
 0.925 &12  &12  &11      &11       &10              &10                 \\
 0.950 &14  &13  &12      &11       &11              &11                 \\
 0.975 &15  &15  &13      &13       &12              &12                 \\
 1.000 &17  &16  &15      &14       &14              &13                 \\
 1.025 &19  &18  &16      &15       &15              &15                 \\
 1.050 &22  &21  &18      &17       &17              &16                 \\
 1.075 &26  &24  &21      &19       &19              &18                 \\
 1.100 &30  &28  &23      &22       &21, 21, 20, 19  &21, 20, 20, 19, 19 \\
 1.125 &37  &33  &27      &25, 24   &25, 23, 21, 20  &24, 22, 21, 20, 20 \\
 1.150 &45  &40  &32, 30  &30, 26   &28, 25, 23, 21  &27, 24, 23, 22, 20 \\
 1.175 &58  &51  &38, 35  &35, 29   &33, 28, 25, 22  &32, 27, 25, 24, 22 \\
 1.200 &81  &67  &47, 41  &42, 33   &40, 33, 29, 23  &38, 31, 28, 26, 23 \\
 1.225 &127 &97  &61, 50  &53, 39   &50, 39, 33, 25  &47, 37, 33, 30, 25 \\
 1.250 &    &    &83, 66  &70, 48   &64, 49, 40, 29  &60, 46, 39, 35, 29 \\
 1.275 &    &    &128, 98 &99, 65   &89, 66, 53, 35  &81, 60, 51, 44, 34 \\
 1.300 &    &    &        &165, 103 &140, 101, 78, 48  &121, 87, 72, 62, 46 \\
 1.325 &    &    &        &         &~~~~-, ~~~-, 166, 97   &227, 161, 132, 112, 81 \\ \hline
\end{tabularx}
\vspace*{0.25cm}
\begin{tabularx}{14.2cm}{|p{1cm}|p{2.3cm}p{2.3cm}p{2.3cm}Xp{1cm}|}
\multicolumn{6}{c}{\begin{small}\textbf{Tabelle 13. Abschätzungen für $N(\lambda)$ unter verschiedenen Voraussetzungen}\end{small}}\\
\hline
  &
 $\lambda_1 \geq 0.64$\newline \begin{tiny}keine Vor.,\newline $N(1.050)\geq 5$, \newline $N(1.050)\geq 10$\end{tiny}&
 $ \geq 0.66$\newline \begin{tiny}keine Vor.,\newline $N(1.050)\geq 5$, \newline $N(1.050)\geq 8$\end{tiny}&
 $ \geq 0.68$\newline \begin{tiny}keine Vor.,\newline $N(1.075)\geq 5$, \newline $N(1.075)\geq 10$\end{tiny}&
 $ \geq 0.72$\newline \begin{tiny}keine Vor.,\newline $N(1.025)\geq 5$, \newline $N(1.025)\geq 8$, \newline $N(1.025)\geq 10$, \end{tiny}&
 $ \geq 0.78$\newline \begin{tiny}keine Vor.\end{tiny}\\ \hline
 $\lambda =$ & $N(\lambda)\geq$&$N(\lambda)\geq$&$N(\lambda)\geq$&$N(\lambda)\geq$&$N(\lambda)\geq$\\
 & & & & & \\
 0.875 &9           &8           &8           &8               &8  \\
 0.900 &9           &9           &9           &9               &8  \\
 0.925 &10          &10          &10          &9               &9  \\
 0.950 &11          &11          &11          &10              &10 \\
 0.975 &12          &12          &12          &11              &11 \\
 1.000 &13          &13          &13          &12              &12 \\
 1.025 &14          &14          &14          &13              &13 \\
 1.050 &16          &16          &15          &15, 14, 14, 14  &14 \\
 1.075 &18, 17, 17  &18, 17, 17  &17          &16, 15, 15, 14  &15 \\
 1.100 &20, 19, 18  &20, 19, 18  &19, 19, 18  &18, 17, 15, 15  &17 \\
 1.125 &23, 21, 19  &22, 20, 19  &22, 20, 19  &21, 18, 17, 15  &19 \\
 1.150 &26, 23, 20  &26, 23, 21  &25, 23, 21  &23, 20, 18, 16  &21 \\
 1.175 &31, 27, 22  &30, 26, 23  &29, 26, 22  &27, 22, 19, 17  &24 \\
 1.200 &37, 31, 25  &35, 29, 26  &34, 29, 25  &31, 25, 21, 19  &28 \\
 1.225 &45, 37, 28  &42, 35, 30  &40, 34, 28  &37, 29, 24, 21  &33 \\
 1.250 &56, 45, 34  &53, 42, 36  &50, 41, 32  &45, 34, 28, 24  &39 \\
 1.275 &74, 58, 42  &68, 54, 45  &64, 52, 39  &56, 42, 33, 28  &47 \\
 1.300 &107, 83, 58  &95, 74, 60  &86, 69, 51  &73, 54, 42, 34  &59 \\
 1.325 &182, 140, 96  &152, 116, 94  &131, 103, 74  &102, 75, 58, 46    &77 \\
 1.350 &            &            &~~~~-, ~~~-,  142  &167, 121, 92, 72   &110\\ \hline
\end{tabularx}
\end{small}
\end{center}
\thinmuskip=3mu
\medmuskip=4mu plus 2mu minus 4mu
\thickmuskip=5mu plus 5mu
\vspace*{0.25cm}

Beim Beweis von Theorem \ref{haupttheorem} werden wir diese Abschätzungen benutzen, indem wir für den Fall $\lambda_1 \in [0.60,0.62]$, beispielsweise, die vier unterschiedlichen Fälle
$N(1.075) \in [0,6],\;[7,10],\;[11,15],\;[16,\infty)$ betrachten. Sind wir im Fall $N(1.075) \in [0,6]$ so benutzen wir diejenigen Abschätzungen, die unter keinen Zusatzvoraussetzungen ("`keine Vor."') bewiesen wurden. Zusätzlich verwenden wir die Voraussetzung $N(1.075)\leq 6$ und bekommen die Werte
$$N(\lambda)\leq \left\lbrace \begin{array}{ll}
6 & \lambda=0.875,\\
6 & \lambda=0.900,\\
\ldots & \\
6 & \lambda=1.075,\\
21 & \lambda=1.100,\\
25 & \lambda=1.125,\\
\ldots & \\
140 & \lambda=1.300.
\end{array} \right. $$
Im Fall $N(1.075) \in [7,10]$ benutzen wir für $\lambda \leq 1.075$ die Abschätzungen, die unter keinen Zusatzvoraussetzungen bewiesen wurden, sowie die Voraussetzung $N(1.075) \leq 10$. Für $\lambda > 1.075$ nehmen wir die Abschätzungen, die unter der Voraussetzung $N(1.075) \geq 7$ bewiesen wurden und erhalten insgesamt
$$N(\lambda)\leq \left\lbrace \begin{array}{ll}
9 & \lambda=0.875,\\
9 & \lambda=0.900,\\
10 & \lambda=0.925\\
\ldots & \\
10 & \lambda=1.075,\\
21 & \lambda=1.100,\\
23 & \lambda=1.125,\\
\ldots & \\
101 & \lambda=1.300.
\end{array} \right. $$
Analog für die beiden anderen Fälle. \\
In jedem der vier Fälle bekommen wir also insgesamt bessere Werte, als wenn wir nur die Werte genommen hätten, die ohne Zusatzvoraussetzungen bewiesen wurden. Letzteres hätten wir tun müssen, wenn wir Verbesserungspotential Nr. 9 nicht berücksichtigt hätten.

\chapter{Beweis von Theorem \ref{haupttheorem}}

\section{Lemma 15.1 aus \cite{Hea92}}
Für den Beweis von Theorem \ref{haupttheorem} benutzen wir \cite[Lemma 15.1]{Hea92}. Jedoch verwenden wir bei der Herleitung des letzten Lemmas anstelle von \cite[Lemma 11.1]{Hea92} das geringfügig bessere Lemma \ref{para11-lemma1}. Wir werden diesen variierten Beweis nicht vollständig aufschreiben, da es genügt zu bemerken, dass an zwei Stellen in der Herleitung eine passende Monotonie vorliegt. Wir verwenden die folgenden in §\ref{para-beweisdeshaupttheorems-einfuehrung} definierten Objekte:
 $$\Sigma, \; L,\; K, \; H(z),\; \sum_{\rho}{'}, \; R_1, \;H_2(z).$$
Weiterhin benutzen wir die zu Beginn von §\ref{para11-neu} definierten Charaktere $\chi^{(k)}$ und deren Nullstellen $\rho^{(k)}$ mit $k \in \{1,2,\ldots,N(\lambda_0)\}$, wobei $\lambda_0=\frac{1}{3}\log \log \mathscr{L}$. In Kapitel 3 haben wir Abschätzungen für die Nullstellen $\rho_k$ und $\rho'$ in $R=R_0(L_1(q))$ bewiesen. Nun gilt für alle hinreichend großen $q$
\begin{eqnarray*}
  R &\supseteq& \{\sigma+it \in \C\,|\,  1-\mathscr{L}^{-1}\lambda_0 \leq \sigma \leq 1, \; |t|\leq 1 \} \\
   &\supseteq& \{\sigma+it \in \C \,|\,  1-\mathscr{L}^{-1}R_1 \leq \sigma \leq 1, \; |t|\leq \mathscr{L}^{-1}R_1 \}\\
   &=:&\tilde{R_1}.
\end{eqnarray*}
 Also haben wir insbesondere auch Abschätzungen für die $\rho^{(k)}$ und letztendlich auch für die Nullstellen $\rho \in \tilde{R_1}$ aus der Summe $\sum_{\rho}{'}$. Man achte jedoch darauf, dass bei der Definition der $\rho_k$ ein Charakter $\chi$ zusammen mit seiner komplexen Konjugation $\overline{\chi}$ betrachtet wurde im Gegensatz zur Definition der $\rho^{(k)}$. Ferner beachte man, dass jegliche Aussagen in diesem Abschnitt nur für $q \geq q_0$ gelten. Dabei ist $q_0$ eine Konstante, die von allen benutzten (später festen) Parametern abhängen kann.\\

\noindent Bevor wir zum Lemma kommen, möchten wir an die Ausgangslage für den Beweis von Theorem \ref{haupttheorem} erinnern. Das ist die Ungleichung (\ref{para13-Sigma-abzuschaetzender-Term}), wegen $H(0)=K^2$ also
\begin{equation}\label{para15-Ausgangspunkt}\frac{K^{-2}\varphi(q)}{\mathscr{L}}\Sigma \geq 1-K^{-2} \sum_{\chi \neq \chi_0} \sum_{\rho}{'} |H((1-\rho)\mathscr{L})| - \varepsilon . \end{equation}
Ziel ist es die letzte rechte Seite so in Abhängigkeit von $\lambda_1$, $\lambda_2$, $\lambda_3$ und $\lambda'$ nach unten abzuschätzen, dass man die in Kapitel 3 bewiesenen Resultate einspeisen und schließlich $\Sigma>0$ folgern kann. In diesem Zusammenhang sei darauf hingewiesen, dass im Folgenden die Existenz gewisser Nullstellen $\rho$ benutzt wird. Wenn diese nicht existieren würden, so wird ersichtlich sein, dass noch bessere Abschätzungen folgen würden. Wenn beispielsweise $\rho_1$ nicht existiert, so liefert letztere Ungleichung sofort die zulässige Linniksche Konstante $L=3+\varepsilon$.\\

\noindent Zuerst wird in (\ref{para15-Ausgangspunkt}) die Summation über "`viele"' Nullstellen für jeden Charakter reduziert auf eine Summe über jeweils eine Nullstelle für jeden Charakter. Das geschieht mittels Lemma \ref{hb-lemma-13-3} und führt zu \cite[(15.1)]{Hea92}:
\begin{equation}\label{para15-Ungl-15-1}
   K^{-2}\sum_{\chi \neq \chi_0} \sum_{\rho}{'} |H((1-\rho)\mathscr{L})| \leq \varepsilon+\sum_{k} e^{-(L-2K)\lambda^{(k)}} B(\lambda^{(k)}) - n(\chi_1) A(\chi_1),
\end{equation}
wobei die Summation über alle weiter oben erwähnten Nullstellen $\rho^{(k)}$ läuft und die Bezeichnungen
\begin{eqnarray*}
B(\lambda)&=&\left(\frac{1-e^{-2K\lambda}}{6K^2 \lambda}\right)+\frac{2K\lambda-1+e^{-2K\lambda}}{2K^2\lambda^2},  \\
  A(\chi_1) &=& \left\lbrace \begin{array}{ll}
(e^{-(L-2K)\lambda_1}-e^{-(L-2K)\lambda'})(B(\lambda_1)-\alpha(\chi_1) K^{-2}H_2(\lambda_1))& \text{ falls } \rho_1 \in \tilde{R_1}, \\
0 & \text{ sonst},
\end{array} \right.  \\
\alpha(\chi_1)&=&\left\lbrace \begin{array}{ll}
2 & \text{ wenn } \chi_1 \text{ reell und } \rho_1 \text{ komplex,} \\
1 & \text{ sonst,}\\
\end{array} \right.  \\
n(\chi_1)&=&\left\lbrace \begin{array}{ll}
2 & \text{ wenn } \chi_1 \text{ komplex,} \\
1 & \text{ wenn } \chi_1 \text{ reell,}\\
\end{array} \right.
\end{eqnarray*}
verwendet wurden. In \cite[§15]{Hea92} wird der Fall $\rho_1 \notin \tilde{R_1}$ nicht betrachtet. Zwar liefert dieser Fall letztendlich viel bessere Werte, jedoch erscheint es uns notwendig ihn so wie eben getan zu berücksichtigen, indem wir dann $A(\chi_1)=0$ wählen. Dies ist gemäß der Herleitung in \cite[S.328-329]{Hea92} offensichtlich möglich. Im Falle, dass $\rho'$ nicht existiert, erhält man die gleiche Abschätzung, wobei dann in $A(\chi_1)$ das $\lambda'$ mit "`$\infty$"' ersetzt wird (definiere dafür $e^{-\infty}:=0$).\\

\noindent Nun gilt es die rechte Seite von (\ref{para15-Ungl-15-1}) abzuschätzen. Sei $\Lambda>0$ eine später passend zu wählende Konstante. Die Terme mit $\lambda^{(k)}>\Lambda$ schätzt man mittels Lemma \ref{para11-lemma1} und diejenigen mit $\lambda^{(k)}\leq\Lambda$ mittels §\ref{para12-neu} ab.
Dazu definieren wir erstmal für ein $\theta$ in $\R$ die Funktion
\begin{equation}\label{para15-definition-w-1}w_1(t)=e^{-\frac{\theta}{2}t} \min \{t-u+10^{-7},v-u+10^{-7}\}^{\frac{1}{4}} \;\;\; (t \in [u,\infty)) \end{equation}
und wählen
\begin{equation*}\label{para15-definition-w} w(s)=\left(\int_u^x w_1(t)^2 e^{2 s t} \,dt \right)^{-1}\end{equation*}
anstelle der Funktion
$$w(s)=\frac{s}{e^{\frac{73}{30}s}-e^{\frac{16}{15}s}}, $$
die in \cite[S.329 (§15)]{Hea92} gewählt wird. Dabei seien $u$, $v$ und $x$ wie in Lemma \ref{para11-lemma1} mittels zweier positiver Konstanten $c_1$ und $c_2$ definiert. Durch eine kleine Zusatzüberlegung könnte man das "`$10^{-7}$"'  streichen, doch dies sparen wir uns. Die Funktion $w_1(t)$ wird im Zusammenhang mit Verbesserungspotential Nr. 7 vorgeschlagen (\cite[S.337]{Hea92}) und ist in einem leicht anderen Zusammenhang optimal. \\
Wir setzen nun voraus, dass $$L-2K>2x.$$ Dann ist $e^{-(L-2K)s} B(s)/w(s)$ monoton fallend in $s$, da beide Terme von $B(s)$ monoton fallend sind, genauso wie  $$\frac{e^{-(L-2K)s}}{w(s)}=\int_u^x w_1(t)^2 e^{s(2 t-(L-2K))} \,dt.$$
\noindent Mit Lemma \ref{para11-lemma1} bekommt man nun aus (\ref{para15-Ungl-15-1}) wegen der eben erwähnten Monotonie wortwörtlich analog zur Herleitung von \cite[(15.3)]{Hea92} die Aussage
\begin{eqnarray}
  & &K^{-2} \sum_{\chi \neq \chi_0} \sum_{\rho}{'} |H((1-\rho)\mathscr{L})| \label{para15-Ungl-15-3-neu} \\
  & &\leq \varepsilon + \frac{\int_u^x w_1(t)^{-2} \min \{t-u,v-u\}\, dt}{c_1 c_2^2} \cdot \frac{e^{-(L-2K)\Lambda} B(\Lambda)}{w(\Lambda)}+ \sum_{\lambda^{(k)}\leq \Lambda} C(\lambda^{(k)}) - n(\chi_1)A(\chi_1),  \nonumber
\end{eqnarray}
wobei
\begin{equation}\label{para15-definition-C-l}
 C(\lambda)=w(\lambda)\left(e^{-(L-2K)\lambda} \frac{B(\lambda)}{w(\lambda)}-e^{-(L-2K)\Lambda} \frac{B(\Lambda)}{w(\Lambda)} \right) .
\end{equation}
Die Funktion $w(\lambda)$ ist offensichtlich monoton fallend und nicht-negativ. Also gilt in $[0,\Lambda]$ das gleiche auch für $C(\lambda)$ als Produkt zweier nicht-negativer monoton fallender Funktionen. Damit erhält man - wieder mit identischer Herleitung - das folgende Analogon zu \cite[Lemma 15.1]{Hea92}. Dazu definieren wir erstmal
\begin{eqnarray*}
\lambda_3^{\star}&=&\min(\Lambda,\lambda_3), \\
\Lambda_r&=&\Lambda-0.025 r,  \\
s&=&[40(\Lambda-\lambda_3^{\star})], \\
   \lambda_1^{\star}&=&\left\lbrace \begin{array}{ll}
\lambda_1 & \rho_1 \in \tilde{R_1},\\
\lambda'& \text{sonst}, \\
\end{array} \right.  \\
T'&=& \max \big{\{}0, n(\chi_1) (C(\lambda_1^{\star})-A(\chi_1)) \big{\}}+(2-n(\chi_1))C(\lambda_3^{\star}),
\end{eqnarray*}
wobei wieder die Konvention $\lambda'=\infty$ und $C(\infty)=0$ benutzt wird, falls $\rho'$ gar nicht existiert. Das $s$ wurde so gewählt, dass $\Lambda_{s+1}<\lambda_3^{\star} \leq \Lambda_s$ gilt. Ferner habe man für die $s+1$ verschiedenen konkreten Werte $\Lambda_j$ ($j=0,\ldots,s$) gewisse Abschätzungen  \begin{equation}\label{Definition-N0-lambda}N(\Lambda_j)\leq N_0(\Lambda_j) \end{equation} mit konkreten $N_0(\Lambda_j) \;(\geq 4)$ (siehe Tabellen 12 und 13). Wir haben:
\vspace*{0.25cm}
\begin{lemma}[vergleiche Lemma 15.1 in \cite{Hea92}]\label{para15-lemma-15-1-variation}
Es gelten die im Laufe dieses Paragraphen angesprochenen und definierten Bezeichnungen. Seien $\lambda_1 \geq 0.348$ und $L$, $K$, $\theta$, $c_1$, $c_2$, $\varepsilon$  und $\Lambda$ positive reelle Konstanten. Es gelte $L-2K>\max \{ 3,2x\}$. Schließlich sei $q \geq q_0$, wobei $q_0$ von allen festen Parametern und Konstanten abhängt und explizit angebbar ist. Dann gilt
\begin{eqnarray}
   & &K^{-2} \sum_{\chi \neq \chi_0} \sum_{\rho}{'} |H((1-\rho)\mathscr{L})| \nonumber \\
   & &\leq \varepsilon + \frac{\int_u^x w_1(t)^{-2} \min \{t-u,v-u\} \,dt}{c_1 c_2^2} \cdot \frac{e^{-(L-2K)\Lambda} B(\Lambda)}{w(\Lambda)} + \max \{2 C(\lambda_2),0\} \nonumber \\
   & & +(N_0(\Lambda_s)-4)C(\lambda_3^{\star}) + \sum_{r=0}^{s-1}(N_0(\Lambda_r)-N_0(\Lambda_{r+1}))C(\Lambda_{r+1}) + T'. \label{para15-endlemma-rechte-seite}
\end{eqnarray}
\end{lemma}
\vspace*{0.25cm}
\noindent Wenn $\rho_3$ nicht existierte, so würde man (siehe (\ref{para15-Ungl-15-3-neu})) den Term $$(N_0(\Lambda_s)-4)C(\lambda_3^{\star}) + \sum_{r=0}^{s-1}(N_0(\Lambda_r)-N_0(\Lambda_{r+1}))C(\Lambda_{r+1}) \; (\geq 0)$$ im Lemma weglassen. Wenn auch $\rho_2$ nicht existierte, so würde man zusätzlich den Term "`$\max \{2 C(\lambda_2),0\}$"' $(\geq 0)$ weglassen. In beiden Fällen hat man also bessere Abschätzungen. Deswegen werden wir diese Fälle nicht mehr erwähnen. Übrigens steht in \cite[Lemma 15.1]{Hea92} nicht der Ausdruck $\max \{2 C(\lambda_2),0\}$, sondern $2 C(\lambda_2)$. Letzteres kann man schreiben, wenn man $\lambda_2 \leq \Lambda$ voraussetzt, was machbar ist. Analoges gilt für das "`$\max$"' in $T'$.

Um aus der Ungleichung (\ref{para15-Ausgangspunkt}) jetzt das Theorem \ref{haupttheorem} zu beweisen, müssen wir zeigen, dass für $L=5.2$ die rechte Seite von (\ref{para15-endlemma-rechte-seite}) echt kleiner als $1$ ist. Wir benutzen die Parameter
\begin{equation}\label{para15-werte-fuer-parameter}
L=5.2, \; K=0.32, \; \theta=1.15, \; c_1=0.11, \; c_2=0.27,
\end{equation}
welche sich als recht optimal erweisen. Jetzt geht man folgendermaßen vor. Betrachte den Fall \begin{equation}\label{para15-lambda-1-bereiche}\lambda_1 \in [\lambda_{11},\lambda_{12}] \;\;\; (0<\lambda_{11}\leq \lambda_{12} \leq \infty).\end{equation}
 Angenommen, wir haben vermöge Kapitel 3 und/oder \cite{Hea92} explizite Abschätzungen
\begin{equation}\label{para15-Abschaetzungen-lambdas}\lambda' \geq \lambda_{11}', \; \lambda_2 \geq \lambda_{21}, \; \lambda_3 \geq \lambda_{31}.\end{equation}
 Da $\lambda_3^{\star}$ monoton wachsend in $\lambda_3$ ist, sind $C(\lambda_3^{\star})$ und $s$ monoton fallend in $\lambda_3$. Also ist auch der Ausdruck auf der rechten Seite von (\ref{para15-endlemma-rechte-seite}) monoton fallend in $\lambda_3$, wenn $N_0(\Lambda_s) \geq 4$, was immer der Fall sein wird. Analog ist der betrachtete Ausdruck auch monoton fallend in $\lambda_2$. Also können wir $\lambda_2$ und $\lambda_3$ durch deren Abschätzungen nach unten ersetzen.  Wir setzen deswegen
 \begin{equation*}\lambda_{31}^{\star}=\min \{\lambda_{31},\Lambda \} .\end{equation*}
  Was die Monotonie in $\lambda_1$ angeht, so haben wir im Fall $\lambda_1^{\star}=\lambda_1$
\begin{eqnarray*}
  C(\lambda_1^{\star})-A(\chi_1)&=& e^{-(L-2K)\lambda'}B(\lambda_1) - e^{-(L-2K)\Lambda}B(\Lambda)\frac{w(\lambda_1)}{w(\Lambda)} \\
  & &+\alpha(\chi_1) K^{-2}H_2(\lambda_1) (e^{-(L-2K)\lambda_1}-e^{-(L-2K)\lambda'}) \\
  \\
  &\leq&  e^{-(L-2K)\lambda'}B(\lambda_{11}) - e^{-(L-2K)\Lambda}B(\Lambda)\frac{w(\lambda_{12})}{w(\Lambda)} \\
  & &+\alpha(\chi_1) K^{-2}H_2(\lambda_{11}) (e^{-(L-2K)\lambda_{11}}-e^{-(L-2K)\lambda'}),
\end{eqnarray*}
wobei wir im Falle $\lambda_{12}=\infty$ die Konvention $w(\infty)=0$ verwenden. Andererseits gilt im Fall $\lambda_1^{\star}=\lambda'$, dass $C(\lambda_1^{\star})-A(\chi_1) \leq \max \{0, C(\lambda_{11}')\}$, also gilt allgemein immer
\begin{eqnarray*}
  C(\lambda_1^{\star})-A(\chi_1) &\leq& \max \Big{\{}0, C(\lambda_{11}'), e^{-(L-2K)\lambda_{11}'}\max \big{\{}0,B(\lambda_{11})-\alpha(\chi_1)K^{-2}H_2(\lambda_{11})\big{\}}\\
  & &- e^{-(L-2K)\Lambda}B(\Lambda)\frac{w(\lambda_{12})}{w(\Lambda)} +\alpha(\chi_1) K^{-2}H_2(\lambda_{11}) e^{-(L-2K)\lambda_{11}} \Big{\}}\\
  &=:& C^{\star}.
\end{eqnarray*}
  Fassen wir alles zusammen, so bekommen wir mit $\varepsilon=10^{-7}$ die folgende konkrete obere Abschätzung der rechten Seite von (\ref{para15-endlemma-rechte-seite}):
\begin{eqnarray}
  W&:=&10^{-7}  + \frac{\int_u^x w_1(t)^{-2} \min \{t-u,v-u\} \,dt}{c_1 c_2^2} \cdot \frac{e^{-(L-2K)\Lambda} B(\Lambda)}{w(\Lambda)} + \max \{2 C(\lambda_{21}),0\} \label{Definition-W} \\
  &+&(N_0(\Lambda_s)-4)C(\lambda_{31}^{\star}) + \sum_{r=0}^{s-1}(N_0(\Lambda_r)-N_0(\Lambda_{r+1}))C(\Lambda_{r+1}) + (2-n(\chi_1))C(\lambda_{31}^{\star}) +n(\chi_1) C^{\star}. \nonumber
\end{eqnarray}
\noindent Wir betrachten in den nächsten drei Unterabschnitten die drei Fälle $\chi_1$ und $\rho_1$ beide reell, $\chi_1$ reell und $\rho_1$ komplex und $\chi_1$ komplex.
Dabei werden wir in möglichst selbsterklärenden Tabellen genau angeben, für welche Unterfälle (\ref{para15-lambda-1-bereiche}) wir welche Abschätzungen (\ref{Definition-N0-lambda}) und (\ref{para15-Abschaetzungen-lambdas}) benutzen, sowie wie wir den letzten verbleibenden freien Parameter $\Lambda$ jeweils wählen. Dann notieren wir welchen Wert $W$ wir für die jeweiligen Unterfälle bekommen. Für die Berechnung des $W$ benutzen wir dabei immer die in (\ref{para15-werte-fuer-parameter}) festgelegten Parameter. Gilt für alle (endlich viele) betrachteten Fälle, dass $W<1$, so haben wir Theorem \ref{haupttheorem} bewiesen.

\section{Fall 1: $\chi_1$ und $\rho_1$ beide reell}

Nach \cite[Lemma 14.2]{Hea92} können wir für diesen Fall annehmen, dass $\lambda_1 \geq 0.348$. Wir betrachten die sechs Fälle
$$\lambda_1 \in [\lambda_{11},\lambda_{12}]= [0.348,0.40],\; [0.40,0.44],\; [0.44,0.60],\; [0.60,0.68],\; [0.68,0.78],\; [0.78,\infty).$$
Für jeden dieser Fälle benötigen wir Abschätzungen von $\lambda'$, $\lambda_2$, $\lambda_3$ und Werte für $N_0(\Lambda_j)\; (j \in \{0,\ldots,s\})$.
Für die $\lambda'$- und $\lambda_2$-Abschätzungen benutzen wir \cite[Table 4 (§8)]{Hea92} und \cite[Table 7 (§8)]{Hea92}. Als $\lambda_3$-Abschätzungen benutzen wir das Maximum aus den eben genannten $\lambda_2$-Abschätzungen, der Tabelle 10 und der allgemein gültigen Abschätzung $\lambda_3 \geq 0.857$ (\cite[Lemma 10.3]{Hea92}). Werte für $N_0(\Lambda_j)$ lesen wir aus der Tabelle 12 bzw. 13 ab, wobei wir die zu der Bedingung $\lambda_1 \geq \lambda_{11}$ gehörigen allgemeinen Abschätzungen ("`keine Vor."') nehmen. Außerdem verwenden wir einmal auch die triviale Abschätzung $\lambda_2 \geq \lambda_{11}$. Wir bekommen (dabei stehe $W$ für den in (\ref{Definition-W}) definierten Ausdruck)
\begin{center}
\noindent\begin{tabularx}{\textwidth}{|Xp{1.0cm}|p{1.0cm}p{1.0cm}p{1.0cm}p{1.0cm}p{2.3cm}|}
\multicolumn{7}{c}{\textbf{Tabelle 14. $\chi_1$ und $\rho_1$ reell}}\\
\hline
 Fall &
 $W\leq$&
 $\lambda' \geq$&
 $\lambda_2 \geq$&
 $\lambda_3 \geq$&
 $\Lambda =$&
 \begin{small}benutzte Vor. in \newline Tabelle 12/13 \end{small}\\ \hline \hline
 $\lambda_1 \in [0.348,0.40]$ & 0.8250 & 2.108 & 1.29 & 1.29 & 1.29 & \begin{small}(Tabelle 12/13 wird nicht ben.) \end{small}\\
 $\lambda_1 \in [0.40,0.42]$ & 0.9991 & 2.030 & 1.18 & 1.18 & 1.225& $\lambda_1 \geq 0.40$\\
 $\lambda_1 \in [0.42,0.44]$ & 0.9869 & 2.030 & 1.18 & 1.18 & 1.225& $\lambda_1 \geq 0.40$\\
 $\lambda_1 \in [0.44,0.60]$ & 0.9967 & 1.832 & 0.92 & 1.175 & 1.225& $\lambda_1 \geq 0.44$\\
 $\lambda_1 \in [0.60,0.68]$ & 0.9711 & 1.724 & 0.79 & 1.078 & 1.300& $\lambda_1 \geq 0.60$\\
 $\lambda_1 \in [0.68,0.78]$ & 0.9987 & 1.630 & 0.745 & 0.971 & 1.325& $\lambda_1 \geq 0.68$\\
 $\lambda_1 \geq 0.78$ & 0.9869 & 1.294 & 0.78 & 0.857 & 1.350& $\lambda_1 \geq 0.78$\\
 \hline
\end{tabularx}
\end{center}

\vspace*{0.25cm}

\section{Fall 2: $\chi_1$ reell und $\rho_1$ komplex}\label{Beweis-Theorem1-1-Fall2}
Nach Lemma \ref{lemma-l1-abschaetzungen-werte} gilt für diesen Fall $\lambda_1 \geq 0.628$.
Für die $\lambda'$- und $\lambda_2$-Abschätzungen benutzen wir Tabelle 3 und 6. Als $\lambda_3$-Abschätzungen nehmen wir das Maximum aus der $\lambda_2$-Abschätzung und der Abschätzung $\lambda_3 \geq 0.857$. Wir bekommen
\begin{center}
\noindent\begin{tabularx}{\textwidth}{|Xp{1.0cm}|p{1.0cm}p{1.0cm}p{1.0cm}p{1.0cm}p{2.3cm}|}
\multicolumn{7}{c}{\textbf{Tabelle 15. $\chi_1$ reell und $\rho_1$ komplex}}\\
\hline
 Fall &
 $W\leq$&
 $\lambda' \geq$&
 $\lambda_2 \geq$&
 $\lambda_3 \geq$&
 $\Lambda =$&
 \begin{small}benutzte Vor. in \newline Tabelle 12/13 \end{small}\\ \hline \hline
 $\lambda_1 \in [0.628,0.74]$ & 0.9884 & 1.52 & 1.02 & 1.02 & 1.275& $\lambda_1 \geq 0.62$\\
 $\lambda_1 \in [0.74,0.78]$ & 0.9598 & 1.46 & 0.93 & 0.93 & 1.350& $\lambda_1 \geq 0.72$\\
 $\lambda_1 \geq 0.78$ & 0.9922 & 1.099 & 0.82 & 0.857 & 1.350& $\lambda_1 \geq 0.78$\\
 \hline
\end{tabularx}
\end{center}

\vspace*{0.25cm}

\section{Fall 3: $\chi_1$ komplex}
Nach Lemma \ref{lemma-l1-abschaetzungen-werte} gilt für diesen Fall $\lambda_1 \geq 0.44$.
Für die $\lambda'$- und $\lambda_2$-Abschätzungen benutzen wir Tabelle 2 und 7. Als $\lambda_3$-Abschätzungen nehmen wir das Maximum aus Tabelle 8/9 und $\lambda_3 \geq 0.857$. \\
Exemplarisch wollen wir die siebte Zeile der Tabelle 16 näher erklären. In dieser Zeile betrachten wir den Fall ($\lambda_1 \in [0.64,0.66]$ und $\lambda_2 \leq 0.86$). Für diesen Fall betrachten wir die drei Unterfälle $N(1.050) \in [0,4],\; [5,9],\; [10,\infty)$. Befinden wir uns im ersten Unterfall dann wählen wir $\Lambda=1.300$ und benutzen für $N_0(\Lambda_j)$ die Werte aus Tabelle 13, welche unter der Bedingung $\lambda_1 \geq 0.64$ und ohne weitere Voraussetzungen bewiesen wurden ("`keine Vor."'). Zusätzlich benutzen wir auch die Bedingung $N(\lambda)\leq 4$ für $\lambda \leq 1.050$. Vergleiche diesbezüglich auch den Kommentar nach Tabelle 13. Weiterhin liefern Tabelle 2, 7 und 9 jeweils $\lambda'\geq 1.04$, $\lambda_2 \geq 0.79$ und $\lambda_3 \geq 0.938$. Mit diesen Werten bekommen wir $W \leq 0.9921$. Analog für die anderen Unterfälle und Fälle.
\begin{center}
\begin{small}
\noindent\begin{tabularx}{\textwidth}{|Xp{1.3cm}|p{0.75cm}p{0.75cm}p{0.75cm}p{0.75cm}p{3.3cm}|}
\multicolumn{7}{c}{\textbf{Tabelle 16. $\chi_1$ komplex}}\\
\hline
 Fall &
 $W\leq$&
 $\lambda' \geq$&
 $\lambda_2 \geq$&
 $\lambda_3 \geq$&
 $\Lambda =$&
 \begin{small}benutzte Vor. in \newline Tabelle 12/13 \end{small}\\ \hline \hline
 $\lambda_1 \in [0.44,0.54]$ & 0.9577 & 1.34 & 1.19 & 1.243 & 1.243 &  \begin{small}(Tabelle 12/13 wird nicht ben.) \end{small}\\ \hline
 $\lambda_1 \in [0.54,0.58]$ & & 1.23 & 1.04 & 1.079   &  & $\lambda_1 \geq 0.54$\\
  & 0.9956 &  &  &  & 1.250 &  $\;\;\;\; N(1.125) \leq 8$\\
  & 0.9954 &  &  &  & 1.275 &  $\;\;\;\; 9 \leq N(1.125)$\\ \hline
 $\lambda_1 \in [0.58,0.60]$ &  & 1.18 & 0.97 & 1.001 &  & $\lambda_1 \geq 0.58$\\
  & 0.9977 &  &  &  & 1.275 &  $\;\;\;\; N(1.100) \leq 10$\\
  & 0.9969 &  &  &  & 1.275 &  $\;\;\;\;\ 11 \leq N(1.100) $\\ \hline
 $\lambda_1 \in [0.60,0.62]$ &  & 1.13 & 0.91 & 0.933 & & $\lambda_1 \geq 0.60$\\
  & 0.9957 &  &  &  & 1.275 &  $\;\;\;\; N(1.075) \leq 6$\\
  & 0.9935 &  &  &  & 1.300 &  $\;\;\;\; 7 \leq N(1.075) \leq 10$\\
  & 0.9996 &  &  &  & 1.300 &  $\;\;\;\; 11 \leq  N(1.075) \leq 15$\\
  & 0.9634 &  &  &  & 1.325 &  $\;\;\;\; 16 \leq N(1.075) $\\ \hline
 $\lambda_1 \in [0.62,0.64]$ &  & 1.09 & 0.85 & 0.902 & & $\lambda_1 \geq 0.62$\\
  & 0.9983 &  &  &  & 1.300 &  $\;\;\;\; N(1.075) \leq 6$\\
  & 0.9952 &  &  &  & 1.300 &  $\;\;\;\; 7 \leq N(1.075) \leq 9$\\
  & 0.9998 &  &  &  & 1.325 &  $\;\;\;\; 10 \leq N(1.075) \leq 11$\\
  & 0.9989 &  &  &  & 1.325 &  $\;\;\;\; 12 \leq N(1.075) \leq 14$\\
  & 0.9813 &  &  &  & 1.325 &  $\;\;\;\; 15 \leq N(1.075)$\\  \hline
 $\lambda_1 \in [0.64,0.66]$&  & 1.04 & 0.79 & 0.938 &  & $\lambda_1 \geq 0.64$\\
 und $\lambda_2 \leq 0.86$ & 0.9921 &  &  &  & 1.300 &  $\;\;\;\; N(1.050) \leq 4$\\
  & 0.9971 &  &  &  & 1.300 &  $\;\;\;\; 5 \leq N(1.050) \leq 9$\\
  & 0.9846 &  &  &  & 1.325 &  $\;\;\;\; 10 \leq N(1.050)$ \\ \hline
  $\lambda_1 \in [0.64,0.66]$ &  & 1.04 & 0.86 & 0.898 &  & $\lambda_1 \geq 0.64$\\
  und $\lambda_2 \geq 0.86$& 0.9688 &  &  &  & 1.300 &  $\;\;\;\; N(1.050) \leq 4$\\
  & 0.9942 &  &  &  & 1.300 &  $\;\;\;\; 5 \leq N(1.050) \leq 9$\\
  & 0.9869 &  &  &  & 1.325 &  $\;\;\;\; 10 \leq N(1.050) $\\  \hline
 $\lambda_1 \in [0.66,0.68]$ &  & 1.00 & 0.74 & 0.960 & & $\lambda_1 \geq 0.66$\\
  und $\lambda_2 \leq 0.83$& 0.9995 &  &  &  & 1.300 &  $\;\;\;\; N(1.050) \leq 4$\\
  & 0.9800 &  &  &  & 1.325 &  $\;\;\;\; 5 \leq N(1.050) \leq 7$\\
  & 0.9939 &  &  &  & 1.325 &  $\;\;\;\; 8 \leq N(1.050) $\\  \hline
 $\lambda_1 \in [0.66,0.68]$ &  & 1.00 & 0.83 & 0.893 & & $\lambda_1 \geq 0.66$\\
 und $\lambda_2 \geq 0.83$ & 0.9997 &  &  &  & 1.300 &  $\;\;\;\; N(1.050) \leq 7$\\
  & 0.9998 &  &  &  & 1.325 &  $\;\;\;\; 8 \leq N(1.050) $\\  \hline
$\lambda_1 \in [0.68,0.72]$ &  & 0.93 & 0.702 & 0.962 & & $\lambda_1 \geq 0.68$\\
und $\lambda_2 \leq 0.81$ & 0.9955 &  &  &  & 1.325 &  $\;\;\;\; N(1.075) \leq 4$\\
  & 0.9977 &  &  &  & 1.325 &  $\;\;\;\; 5 \leq N(1.075) \leq 9 $\\
  & 0.99984 &  &  &  & 1.325 &  $\;\;\;\; 10 \leq N(1.075) $\\ \hline
  $\lambda_1 \in [0.68,0.72]$ &  & 0.93 & 0.81 & 0.883 & & $\lambda_1 \geq 0.68$\\
  und $\lambda_2 \geq 0.81$& 0.9438 &  &  &  & 1.325 &  $\;\;\;\; N(1.075) \leq 4$\\
  & 0.9859 &  &  &  & 1.325 &  $\;\;\;\; 5 \leq N(1.075) \leq 9 $\\
  & 0.9987 &  &  &  & 1.350 &  $\;\;\;\; 10 \leq N(1.075) $\\ \hline
   $\lambda_1 \in [0.72,0.84]$ &  & 0.827 & 0.72 & 0.857 &  & $\lambda_1 \geq 0.72$\\
 & 0.9963 &  &  &  & 1.325 &  $\;\;\;\; N(1.025) \leq 4$\\
  & 0.9960 &  &  &  & 1.350 &  $\;\;\;\; 5 \leq N(1.025) \leq 7 $\\
  & 0.9923 &  &  &  & 1.350 &  $\;\;\;\; 8 \leq N(1.025) \leq 9 $\\
  & 0.9874 &  &  &  & 1.350 &  $\;\;\;\; 10 \leq N(1.025) $\\ \hline
 $\lambda_1 \geq 0.84$ & 0.9871 & 0.84 & 0.84 & 0.857 & 1.350  & $\lambda_1 \geq 0.78$\\
 \hline
\end{tabularx}
\end{small}
\end{center}
\vspace*{0.5cm}
In allen obigen Fällen ist $W<1$. Damit ist der Beweis von Theorem \ref{haupttheorem} erbracht.

\chapter{Verbesserungspotentiale, die nicht benutzt wurden}
Wir diskutieren in diesem Kapitel kurz jene Verbesserungspotentiale aus \cite[S.332-337]{Hea92}, die wir letztendlich nicht verwendet haben. Das liegt daran, dass die resultierenden Verbesserungen\footnote{Damit meinen wir natürlich die aus \emph{unserem} Vorgehen resultierenden Verbesserungen. Es mag andere Vorgehensweisen geben, die zu größeren Verbesserungen führen.} für Theorem \ref{haupttheorem} zu gering sind (darunter verstehen wir Verbesserungen von $\leq 0.005$ für die zulässige Konstante $L$). Bei einigen Potentialen ist es dabei so, dass man zwar hilfreiche Verbesserungen bekommt, dies aber nur für einen Teil der später auftauchenden Fälle (z.B. nur für ca. $\lambda_1 \leq 0.60$ bei Potential Nr. 8), für den wir aber keine Verbesserungen benötigen.\\
Es sei nochmal darauf hingewiesen, dass wir in diesem Kapitel kein Interesse an ausführlichen und rigorosen Beweisen haben. Stattdessen möchten wir nur verständlich aufzeigen, dass die erzielten Verbesserungen zu gering sind - zumindest für den Fall, dass man keine größeren Veränderungen bzw. Verfeinerungen an den Potentialen vornimmt.

\section{Verbesserungspotential Nr. 1}\label{para-vb1}
Für diesen Abschnitt verweisen wir auf die Bezeichnungen und Kommentare in §\ref{hb-para-7}.\\
Den Funktionen $f$, die wir in Lemma \ref{hb-lemma-5-3} und \ref{hb-lemma-5-2-extended} einsetzen, fällt eine große Bedeutung zu. Die Konstanten, die wir für unsere Abschätzungen bekommen hängen Eins zu Eins damit zusammen, wie die benutzte Funktion $f$ aussieht. Insofern ist es bestimmt sinnvoll auch nach Funktionen $f$ zu schauen, die Bedingung 1 und 2 erfüllen, aber nicht aus dem speziellen Ansatz in \cite[§7]{Hea92} stammen.\\
Als "`Verbesserungspotential Nr. 1"' schlägt Heath-Brown vor, beim Aufstellen der "`Variationsbedingung"' etwas allgemeiner vorzugehen, als es in \cite[§7]{Hea92} gemacht wurde. Dies führt zu einer Differentialgleichung 4. Ordnung, die die optimale Funktion $g_1$ erfüllen muss. Dieser Ansatz erfordert recht lange Rechnungen und verschiedene Fallbesprechungen. Letztendlich bekommt man im Wesentlichen folgende Funktionstypen für $g$ ($c_1,c_2,c_3,x_1,x_2 \in \R$ Konstanten):\\
$$g(t)= c_1 co_1(x_1 t) + c_2 co_2(x_2 t) - c_3, $$
$$g(t)= c_1 \cosh( x_1 t) \cos(x_2 t) + c_2 \sinh(x_1 t) \sin(x_2 t) - c_3 ,$$
wobei $co_i(t)=\cos(t)$ oder $co_i(t)=\cosh(t)$ gewählt wird $(i \in \{1,2\})$.
Was die optimale Wahl der Konstanten angeht, so ergeben sich verschiedene Bedingungen mit dessen Hilfe man in den meisten Fällen konkrete Werte für die Konstanten bekommt. Alternativ kann man aber auch für eine obige Funktion $g$ per Hand (bzw. Computer) durch "`numerisches Experimentieren"' annähernd optimale Werte für die Konstanten finden. Letzteres ist einfacher.\\

\subsubsection{Welche Verbesserungen bekommt man mit den neuen $g$?}

Mit diesem Verbesserungspotential erhalten wir nur minimale, nicht nennenswerte Verbesserungen:\\
Beispielsweise können wir in \cite[Table 9-10]{Hea92} nur ein paar der Abschätzungen $\lambda_2 \geq C$ verbessern und dies nur um $0.001-0.002$ für die Konstante $C$. In \cite[Table 2-5]{Hea92} bekommen wir für kleine $\lambda_1$ etwas bessere Verbesserungen im "`Hundertstel-Bereich"'. Doch auch diese lohnen sich nicht, einerseits weil wir für den in \cite[§8]{Hea92} diskutierten Fall keine Verbesserungen nötig haben, andererseits weil diese Verbesserungen nur einen minimalen Effekt auf die zulässige Endkonstante $L$ in diesem Fall haben.\\
Letztendlich sollte man beachten, dass die oben erwähnten neuen Funktionen $g$ qualitativ "`genauso aussehen"' wie die alten: Zum Beispiel kann man für die Verbesserung des ersten Eintrags in \cite[Table 9 (§9)]{Hea92} um $0.001$ auf $\lambda_j \geq 1.284$ die Funktion
$g_1(t)=c_1 \cosh(x_1 t)+c_2 \cos(x_2 t) - c_3$ mit den Konstanten $c_1=0.3$, $x_1=0.5$, $c_2=2.5$, $x_2=1.0$, $c_3=c_1 \cosh(x_1 \gamma)+c_2 \cos(x_2 \gamma)$, $\gamma=1.48$ benutzen. Setze andererseits $g_2(t)=\gamma^2-t^2$ mit dem gleichen $\gamma$. Letztere Funktion liefert $\lambda_j \geq 1.281$ (und durch $\gamma=1.44$ bekäme man $\lambda_j \geq 1.282$). Schließlich sehen die Graphen der beiden Funktionen sehr ähnlich aus:\\

\begin{center}
\textbf{Neues und altes $g$}\\
\includegraphics[scale=0.75]{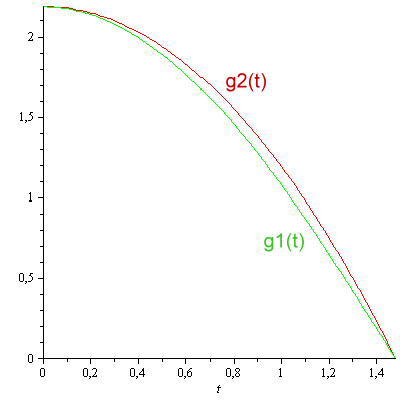}\\
\end{center}
\vspace*{0.5cm}

\textbf{Bemerkungen zu weiteren fehlgeschlagenen Verbesserungsansätzen:}

\begin{itemize}
\item Beginnend mit dem Ansatz (\ref{hb-para7-ansatz-fuer-f}) haben wir eine Fülle von verschiedenen zulässigen Funktionen $g$ ausprobiert. Dies ergab keine Verbesserungen im Vergleich zu den uns bereits bekannten Funktionen $g$.

\item Wir haben versucht andere sinnvolle Funktionen $f$ zu finden, die Bedingung 1 und 2 erfüllen und nicht aus dem speziellen Ansatz (\ref{hb-para7-ansatz-fuer-f}) herkommen. Dies ist uns nicht gelungen.

\item Weiterhin haben wir auch Funktionen ausprobiert, die nicht vollständig die Bedingung 2 erfüllen. Für gewisse Funktionen, die teilweise negativ sind, bekommt man einerseits bessere Werte. Andererseits verursachen sie in den entsprechenden Herleitungen einen Fehler, der den Gewinn wieder (mehr als) aufzehrt.\\

In den Fällen, in denen wir im "`zweiten Schritt"' (vergleiche §\ref{hb-para-6})  $\beta =\beta_1$ nehmen und wo wir garantieren können, dass $\min (\lambda_2,\lambda')-\lambda_1 \geq c>0$, könnten wir Funktionen $f$ nehmen, dessen Laplace-Transformierte
nur $$Re \{F(z)\} \geq 0  \text{ für } Re \{z\} \geq c$$ erfüllt (Beispiel: $f(t)=\gamma^4-t^4$). Diese Funktionen $f$ können dann bessere Werte liefern. Jedoch sind die genannten Anfangsvoraussetzungen in unserer Arbeit nicht erfüllt, da wir vermöge Verbesserungspotential Nr. 2 eher $\beta=\max \{\beta', \beta_2\}$ wählen.
\end{itemize}

\section{Verbesserungspotential Nr. 3}
Heath-Brown beschreibt in \cite[S.332-333]{Hea92} eine vergleichsweise einfache Vorgehensweise für dieses Potential, welche minimale Verbesserungen für alle Abschätzungen liefert, die mit der Methode aus §\ref{hb-para-6} bewiesen werden. Die Frage ist, ob man durch eine verfeinerte Vorgehensweise vielleicht hilfreiche Verbesserungen bekäme? Als Antwort auf diese Frage skizzieren wir in diesem Abschnitt ein "`worst case"'-Szenario, welches sich nicht einfach ausschließen lässt und bei welchem dieses Potential keine ausreichenden Verbesserungen liefert (solange man keine stark verbesserten Brun-Titchmarsh Ungleichungen zur Verfügung hat, was aber aus verschiedenen Gründen nicht zu erwarten ist). Da wir dieses Potential ja schlussendlich nicht benutzen, werden wir die Argumentation nicht ausführlich beweisen. Auch gehen wir nicht darauf ein, inwieweit das unten Besprochene wirklich der "`worst case"' ist.\\

Wir setzen voraus, dass $\gamma_1=\gamma_2=0$ und $ord \; \chi_1, ord \; \chi_2 > \mathscr{L}$. Aufgrund der letzten Voraussetzung sind $\chi_1(n)$ und $\chi_2(n)$ für $n \in \N \;  ((q,n)=1)$ "`einigermaßen gleichverteilt"' auf dem Einheitskreis\footnote{Beweisskizze: Definiere zuerst, was bewiesen werden soll, d.h. konkretisiere den Begriff "`gleichverteilt"'. Wegen $ord \, \chi \leq kgV_{n \, mod \, q} \, (ord_{\C^{\star}} \chi (n))$ gibt es zu jedem $M \in \N$ für $q\geq q_0(M)$ ein $n_1$ mit $\chi_j(n_1)$ ist primitive $N$-te Einheitswurzel, wobei $N\geq M$. Benutze dann, dass $(\Z / q\Z)^{\star}$ endliches Produkt endlicher zyklischer Gruppen ist.}.\\
Wir zeigen repräsentativ anhand der folgenden Ungleichung, dass dieses Potential zu geringe Verbesserungen liefert. Wir verweisen an dieser Stelle auf die Kommentare und Notationen, die in §\ref{hb-para-6} festgelegt wurden. Es gilt
\begin{eqnarray}
0 &\leq& \chi_0(n) (1+Re \{\chi_1(n)\})(1+Re \{\chi_2(n)\}) \label{vb3-anfangs-abschaetzung}\\
 &=&\chi_0(n) + Re  \{\chi_1(n)\}+Re  \{\chi_2(n)\} + \frac{1}{2}Re \{\chi_1\chi_2(n)\}+\frac{1}{2}Re \{\chi_1\overline{\chi_2}(n)\} \label{vb3-anfangs-gleichheit}
. \end{eqnarray}

\noindent Der nächste Schritt wäre mit $\Lambda(n)n^{-\beta_1} f(\mathscr{L}^{-1}\log n)$ zu multiplizieren und über $n$ zu summieren. Dabei bezeichne $f$ eine Funktion die Bedingung 1 und 2 erfüllt. Heath-Brown bemerkt in diesem Verbesserungspotential, dass man hier noch genauer arbeiten kann, denn: Für unbekanntes $n$ ist (\ref{vb3-anfangs-abschaetzung}) zwar optimal, da es durchaus sein kann, dass die rechte Seite in dieser Ungleichung beliebig nah bei Null ist. Wenn man jedoch über alle $n \in \N$ aufsummiert, so kann man benutzen, dass $\chi_j(n) \; (j\in \{1,2\})$ einigermaßen gleichverteilt ist auf dem Einheitskreis, also ist auch für verschiedene $n \; ((n,q)=1)$ der Wert $$r(n)=\chi_0(n) (1+Re \{\chi_1(n)\})(1+Re \{\chi_2(n)\})$$ einigermaßen kontinuierlich verteilt auf $[0,4]$ (jedoch verstärkt bei $0$ und $4$; wir halten hier auch fest, dass $r(n)$ Periode $q$ hat). Das heißt, dass wir für einen großen Anteil der $n$ auf der linken Seite von (\ref{vb3-anfangs-abschaetzung}) einen positiven Wert $c$ nehmen können anstelle von $0$. Wir benutzen also statt (\ref{vb3-anfangs-abschaetzung}) die Gleichheit (\ref{vb3-anfangs-gleichheit}) und bekommen
\begin{eqnarray*}
S&:=&\sum_{n=1}^{\infty} \Lambda (n) n^{-\beta_1} f(\mathscr{L}^{-1} \log n) r(n)\\
&\leq& K(\beta_1,\chi_0)+K(\beta_1,\chi_1)+K(\beta_1,\chi_2)+\frac{1}{2}K(\beta_1,\chi_1 \chi_2) +\frac{1}{2} K(\beta_1, \chi_1 \overline{\chi_2}).
\end{eqnarray*}
\noindent Nach Definition der $\chi_j$ gilt, dass $\chi_1 \chi_2, \chi_1 \overline{\chi_2} \neq \chi_0$. Wir schätzen jeden der fünf Terme auf der rechten Seite der letzten Ungleichung mit Hilfe der Lemmata \ref{hb-lemma-5-3} (erster Term) und \ref{hb-lemma-5-2-extended} (die vier anderen Terme) auf die übliche Art und Weise ab und bekommen\footnote{Wir erinnern daran, dass wir $\gamma_1=\gamma_2=0$ voraussetzten. Diese Voraussetzung vereinfacht die Rechnungen, ist aber nicht zwingend notwendig.}
\begin{equation}\label{vb3-inequality}\mathscr{L}^{-1} S \leq F(-\lambda_1)- F(0)-F(\lambda_2-\lambda_1) + \frac{f(0)}{2} + \varepsilon . \end{equation}

\textbf{Triviale Abschätzung von $S$}\\
\noindent Wir benutzen die triviale Abschätzung $S\geq 0$. Setzen wir dann z.B. $\lambda_1 \leq 0.54$ voraus, so ist $\gamma=1.07$ recht optimal und wir bekommen
$$\lambda_1 \leq 0.54 \Rightarrow \lambda_2 \geq 0.656. $$
Wir möchten sehen welche Verbesserung wir von $0.656$ ungefähr erhalten würden, wenn wir $S$ mit dem Verfahren aus Verbesserungspotential Nr. 3 abschätzten.\\

\textbf{Nicht-triviale Abschätzung von $S$ (=Verbesserungspotential Nr.3)}\\
Vorneweg weisen wir darauf hin, dass man analog zum Beweis von Lemma \ref{hb-lemma-5-3} hat, dass für $0 \leq \alpha_1<\alpha_2$
\begin{equation}\label{vb3-lemma-5-3-abschaetzung}\sum_{q^{\alpha_1} < n \leq q^{\alpha_2}} \chi_0(n)\Lambda (n) n^{-\beta_1} f(\mathscr{L}^{-1} \log n) \leq \mathscr{L} \int_{\alpha_1}^{\alpha_2} f(t) e^{\lambda_1 t} \,dt + \varepsilon \mathscr{L}.\end{equation}

\noindent Wir möchten nun $S$ so gut es geht nach unten abschätzen, möchten also gerne, dass die Terme in der Summe $S$ so groß wie möglich sind. Wir wissen, dass der Wert von $r(n)$ für die $\varphi(q)$ verschiedenen $n$  ($(n,q)=1, n \leq q$) recht kontinuierlich verteilt ist in $[0,4]$. Die zentrale Problematik bei der Abschätzung von $S$ liegt nun in der Frage, ob diejenigen Restklassen $n$, für die $r(n)$ groß ist (z.B. $\geq \frac{1}{2}$), nur verhältnismäßig wenig Primzahlpotenzen enthalten. Dann wäre $\Lambda(n)$ dort verstärkt gleich $0$ und würde die großen $r(n)$ wieder zunichte machen.\\
Man mag sich zwar jetzt vielleicht daran erinnern, dass die Primzahlen bzw. Primzahlpotenzen gleichverteilt sind in den primen Restklassen. Dies gilt aber nur wenn $n$ sehr groß gegenüber $q$ ist, z.B. $n \gg e^q$. Für $n \leq q^{x_0}$ jedoch (nur dann ist $f(\mathscr{L}^{-1}\log n)\neq 0$) ist dies nicht bekannt. Schließlich geht es in dieser Diplomarbeit ja gerade darum, dass für gewisse prime Restklassen $a \; mod \; q$ nicht einmal garantiert werden kann, dass \emph{eine} Primzahl in dieser Restklasse für $n\leq q^5$, sagen wir, enthalten ist.\\

\noindent Wir betrachten in $S$ zuerst den Beitrag der $n$ mit $1\leq n \leq q,\; ((n,q)=1)$ (ist $(n,q)>1$, so ist der entsprechende Term in der Summe wegen $\chi(n)=0$ auch Null). Für passende $q$ ist
\begin{equation}\label{vb3-phi-q-lambda-n}\varphi(q) \geq \frac{1}{100} q \; \text{ und } \; \sum_{n \leq q \atop \Lambda(n) \neq 0} 1 =o(q).\end{equation}
Weiterhin kann man für solch kleine $n$ nix über die Verteilung der Primzahlpotenzen in den primen Restklassen $mod \; q$ sagen (die Restklassen haben nur ein Element im Bereich $[1,q]$).  Wegen (\ref{vb3-phi-q-lambda-n}) und der kontinuierlichen Verteilung von $r(n)$ in $[0,4]$  kann man also nicht ausschließen, dass für jegliche $n$, für die $\Lambda(n) \neq 0$ gilt, wir $r(n) < \varepsilon$ haben. Der Beitrag dieser betrachteten $n$ zur Summe $S$ unter diesen Umständen wäre vermöge (\ref{vb3-lemma-5-3-abschaetzung})
$$\leq \varepsilon \mathscr{L} \int_0^1 f(t) e^{\lambda_1 t} \,dt \ll_f \varepsilon \mathscr{L},$$
bei Betrachtung von (\ref{vb3-inequality}) also vernachlässigbar. Der Beitrag der $n \in (q,q^{1+\varepsilon}]$ ist auch vernachlässigbar nach (\ref{vb3-lemma-5-3-abschaetzung}).\\

\noindent Sei jetzt $n \in (q^{1+\varepsilon},q^{x_0}]$. Dann kann man das Problem der Gleichverteilung der Primzahlen/Primzahlpotenzen in den primen Restklassen teilweise lösen mit Hilfe einer Brun-Titchmarsh Ungleichung\footnote{Es gibt verschiedene Versionen der Brun-Titchmarsh Ungleichung. Diese ist für unsere Zwecke die beste, die wir gefunden haben.}:\\
\begin{theorem}[Y. Motohashi, 1974, \cite{Mot74}]
Seien $a,q \in \N$ mit $(a,q)=1$ und $a$ liege nicht in einer gewissen Ausnahmemenge mit $\ll q^{1-0.2\varepsilon}$ Klassen $mod \; q$. Dann gilt für $x \geq q^{1+\varepsilon}$
$$\pi (x;q,a) \leq \frac{2 (1+2\varepsilon) x}{\varphi(q) \log (xq^{-\frac{1}{2}})}=:C \frac{x}{\phi (q) \log x} .$$
Dabei ist $C=C(t)=2(1+2\varepsilon) \frac{t}{t-\frac{1}{2}}$, wenn man $x=q^t$ setzt.\\
\end{theorem}

\noindent Wir nehmen im Folgenden an, dass wir letzteres Theorem für alle $a$ ohne Ausnahmen verwenden dürfen (was auch gehen würde).\\
Um den Beitrag der betrachteten $n$ möglichst genau abzuschätzen, teilen wir die Summation auf in endlich viele Teile $n \in (q^{\alpha_l},q^{\alpha_{l+1}}]$, wobei $1+\varepsilon=\alpha_0<\alpha_1<\ldots<\alpha_M=x_0$. Also bekommen wir $M$ Teilsummen, die wir jeweils mit $S_l$ $(l\in \{1,\ldots,M \})$ bezeichnen wollen. Halte für nachfolgende Überlegungen ein $S_l$ fest und lasse dort die Summation über Primzahlen anstelle von Primzahlpotenzen gehen mit vernachlässigbarem Fehler.\\
Wir können nicht ausschließen, dass die "`schlechtesten"' Restklassen (d.h. diejenigen, für die $r(n)$ am kleinsten ist) alle extrem viele Primzahlen beinhalten im Intervall $(q^{\alpha_l},q^{\alpha_{l+1}}]$. Nach der Brun-Titchmarsh Ungleichung wären das maximal ca. $C x (\phi (q) \log x)^{-1}$ Stück ($x=q^{\alpha_{l+1}}, \;C=C(\alpha_{l+1})$). Alle anderen (besseren) Restklassen könnten keine Primzahl beinhalten. Gemäß diesen Annahmen haben wir dann nach Primzahlsatz ca. $C^{-1} \varphi(q)$ (dies sei oBdA $\in \N$) Restklassen $a \; mod \; q$ (sei $A$ die Menge dieser $a$), welche alle im Intervall $(q^{\alpha_l},q^{\alpha_{l+1}}]$ ca. $C x (\phi (q) \log x)^{-1}$ Primzahlen enthalten. Die anderen Restklassen besitzen in diesem Intervall keine Primzahlen. Weiterhin ist $r(a)$ für $a \in A$ so klein wie möglich, mit anderen Worten besteht also $\{r(a)\,|\, a \in A \}$ genau aus den $C^{-1} \varphi(q)$ kleinsten Werten der Menge $\{r(a) \,| \, 1 \leq a \leq q, \; (a,q)=1 \}$.\\

Nun ist $(1+Re \{\chi_j(a)\})=1+\cos(\theta_j(a))$, wobei $\theta_j(a)$ auf $[0,2\pi]$ näherungsweise gleichverteilt ist für verschiedene $a$ $((a,q)=1)$. Wegen dieser Gleichverteilung und der Konstruktion der Menge $A$ gilt dann für $a \in A$ ungefähr
$$\theta_j(a) \in \pi + \left[-\frac{C^{-1}}{2}2 \pi,\frac{C^{-1}}{2}2 \pi\right].$$
Nimmt man zusätzlich an, dass die Verteilung von $\theta_1(a)$ unabhängig von $\theta_2(a)$ ist, so erhält man im "`Mittel"'\footnote{Nehme hier an, dass $(1+Re \{\chi_j(a)\})$ im Mittel gleich $(1+\cos(\pi +\frac{C^{-1}}{4} 2\pi))$ ist. Gehe dann grob vor und multipliziere die beiden Mittel einfach miteinander, um das Mittel für $r(n)$ zu schätzen. Man beachte, dass aufgrund der Struktur der Funktion $(1+\cos(\pi+t))$  eine genauere Analyse einen noch schlechteren Wert für das sogenannte Mittel liefern würde.}
\begin{equation}\label{vb3-abschaetzung-r}r(n) \approx (1+\cos(\pi +\frac{C(\alpha_{l+1})^{-1}}{4} 2\pi))^2 .\end{equation}

\noindent Setze nun in $S_l$ für $r(n)$ den Wert aus (\ref{vb3-abschaetzung-r}) ein und streiche das $\chi_0(n)$ mit vernachlässigbarem Fehler. Teile die verbleibende Summe in $C^{-1} \varphi(q)$ Summen $S_{l,a}$ über die verschiedenen Restklassen $a \in A$ ein; die anderen Restklassen liefern nach Annahme keinen Beitrag. Benutzt man die obige Annahme zu der Anzahl der Primzahlen in diesen Restklassen, so erhält man insgesamt

\begin{eqnarray*}
\mathscr{L}^{-1} S&  \leq& 2\varepsilon + \mathscr{L}^{-1}\sum_{l=1}^M \sum_{a \in A} \sum_{q^{\alpha_l}<p \leq q^{\alpha_{l+1}}, \atop p \equiv a \; mod \; q}  (\log p) n^{-\beta_1} f(\mathscr{L}^{-1} \log p) r(a) \\
&\lessapprox & 3\varepsilon + \sum_{l=1}^M (1+\cos(\pi + \frac{C(\alpha_{l+1})^{-1}}{4} 2 \pi))^2 \int_{\alpha_l}^{\alpha_{l+1}} f(t) e^{\lambda_1 t} \,dt\\
& \leq & 4\varepsilon + \int_{1}^{x_0} f(t) e^{\lambda_1 t}(1+\cos(\pi + \frac{C(t)^{-1}}{4} 2 \pi))^2 \,dt
.
\end{eqnarray*}
Bei der letzten Ungleichung benutzt man, dass $(\alpha_{l+1}-\alpha_l)$ klein genug ist (in Abhängigkeit von $\varepsilon$ und $f$).\\

\noindent Ersetzt man nun $\mathscr{L}^{-1} S$ in (\ref{vb3-inequality}) durch den rechten Wert der letzten Ungleichung (das wäre das bestmögliche, was man als untere Abschätzung für $\mathscr{L}^{-1} S$ beweisen könnte), so erhält man

\begin{equation}\label{vb3-inequality-v2}0 \leq F(-\lambda_1)- F(0)-F(\lambda_2-\lambda_1) + \frac{f(0)}{2} - \int_1^{x_0} f(t) e^{\lambda_1 t}(1+\cos(\pi + \frac{\pi}{2} C(t)^{-1}))^2 \,dt+\varepsilon .\end{equation}

\noindent Optimiert man letzte Ungleichung nach $f$, so bekommt man für $\gamma=1.09$
$$\lambda_1 \leq 0.54 \Rightarrow \lambda_2 \geq 0.664 . $$

\noindent Im Vergleich zu der trivialen Abschätzung von $S$ ist dies eine Verbesserung um $0.008$. Verbesserungen in dieser Größenordnung helfen uns nicht weiter (ab ca. $0.02$ schon). Außerdem bemerken wir, dass wir an einer Stelle zu gut rechneten (nämlich bei (\ref{vb3-abschaetzung-r})). Weiterhin benutzen wir später anstelle von (\ref{vb3-anfangs-abschaetzung}) etwas ungünstigere Ungleichungen  (mit drei anstatt zwei Faktoren), so dass die Verbesserung dann nochmals kleiner wird. Letztendlich kann man bei rigorosem Vorgehen Verbesserungen von ca. $0.001$ erwarten.

\section{Verbesserungspotential Nr. 4}
In den Ungleichungen vom Typ (\ref{vb3-anfangs-gleichheit}) tauchen immer verschiedene Charaktere $\chi$ auf. In der Regel mussten wir für alle $\chi$ die Abschätzung $\phi(\chi)\leq \frac{1}{3}$ benutzen. Durch eine genauere Analyse der Größe der verschiedenen $\phi(\chi)$ erhält Heath-Brown Verbesserungen für die $\lambda_2$-Abschätzungen im Fall, dass $\chi_1$ und $\rho_1$ beide reell sind.\\ Konkret sieht die Verbesserung so aus, dass man in \cite[Lemma 8.5]{Hea92} den Wert $\frac{11}{24}=0.4583...$ durch $\frac{97}{216}=0.4490...$ ersetzen kann. Dies führt zu entsprechenden Verbesserungen an den Stellen, an denen letzteres Lemma angewendet wurde. Also in \cite[Lemma 8.6]{Hea92} (ersetze $\frac{12}{11}$ durch $\frac{108}{97}$) und damit auch in \cite[Theorem 4 (S.268)]{Hea92}. Außerdem bekommt man bessere Werte in \cite[Table 5,6]{Hea92} (in "`Table 5"' bekommt man Verbesserungen der $\lambda_2$-Abschätzungen um ca. $0.05-0.15$). Wie man am Beweis von Theorem \ref{haupttheorem} erkennen kann, haben wir diese Verbesserungen nicht nötig.

\section{Verbesserungspotential Nr. 6}
Vergleiche für diesen Abschnitt die Kommentare und Notationen aus (\ref{para11-neu}). Die Konstanten $w,u,v,x$ wählen wir dabei so wie in \cite[S.321 und Lemma 11.1]{Hea92}. Also $w=c_1-\varepsilon$, $u=\phi+2c_1$, $v=\phi+2c_1+c_2$, $x=2\phi+3c_1+c_2$ mit $c_1=\frac{1}{10}$ und $c_2=\frac{1}{4}$. \\

Heath-Brown beschreibt in \cite[S.335]{Hea92} eine grobe Vorgehensweise für dieses Potential, bei der man eine Verbesserung der Konstante auf der rechten Seite der Ungleichung \cite[(11.4)]{Hea92} bekommt und zwar von $\frac{67}{6}=11.1666...$ auf  $\frac{67}{6}-\eta$ mit einem gewissen $\eta < 0.0001$. Für eine Verbesserung der zulässigen Linnikschen Konstante $L$ um $0.005$ bräuchten wir grob gesagt $\eta \approx 0.10$. Es bleibt die Frage, ob man durch eine Verfeinerung der Argumente dies erreichen kann?\\
Wir führten eine gewisse Verfeinerung durch und bekamen nur $\eta \approx 0.01$. Im Folgenden möchten wir diese Verfeinerung skizzieren und dabei gleichzeitig zeigen, um was es sich bei diesem Potential handelt.\\

Wir erinnern an die Ungleichung (\ref{para11-ungleichung-am-anfang}), also
\begin{equation}\label{vb6-anfangsungleichung}\sum_{\chi} w_{\chi} \leq (1+ O(\mathscr{L}^{-1})) \sum_{\chi} \left|\sum_n a_{n \chi} b_n \right|^2, \end{equation}
wobei $a_{n \chi}$ und $b_n$ Parameter sind, die
$$a_{n \chi} b_n =w_{\chi}^{\frac{1}{2}} (\sum_{d|n} \psi_d )(\sum_{d|n} \theta_d ) \chi(n) n^{-\rho(\chi)} (e^{-n/X}-e^{-n\mathscr{L}^2/U}) $$
erfüllen. Wir lassen jetzt den Parameter $w_1(t)$ der Einfachheit halber weg und wählen wie in \cite[S.319]{Hea92}
\begin{equation}\label{vb6-definition-a-n}a_{n \chi}=w_{\chi}^{\frac{1}{2}} (\sum_{d|n} \theta_d ) \chi(n) n^{\frac{1}{2}-\rho(\chi)} (e^{-n/X}-e^{-n\mathscr{L}^2/U})^{\frac{1}{2}}\end{equation}
und
\begin{equation}\label{vb6-definition-b-n}b_n=(\sum_{d|n} \psi_d )n^{-\frac{1}{2}}(e^{-n/X}-e^{-n\mathscr{L}^2/U})^{\frac{1}{2}} . \end{equation}

\noindent Heath-Brown bemerkt in Verbesserungspotential Nr. 6, dass für einige $n \in \N$ gilt $$h_1(n):=\sum_{d|n} \theta_d=0.$$ Gleichung (\ref{vb6-anfangsungleichung}) bleibt also wahr, wenn man für diese $n$ dann $b_n=0$ anstelle von (\ref{vb6-definition-b-n}) setzt. Diese leichte Änderung in den $b_n$ führt dann zu einer Verbesserung bei der Konstanten auf der rechten Seite der Gleichung \cite[(11.14)]{Hea92} um
\begin{equation}\label{vb6-summe} S:=\sum_{n=1, \atop h_1(n)=0}^{\infty} (\sum_{d|n} \psi_d )^2 n^{-1} (e^{-n/X}-e^{-n\mathscr{L}^2/U}) .\end{equation}
Kann man $S \geq \eta'>0$ beweisen, so resultiert dies in eine Verbesserung von $\eta=\frac{1}{2w} \eta'=5\eta'$ für die Konstante auf der rechten Seite von \cite[(11.4)]{Hea92}. Für den oben angesprochenen Wert $\eta=0.20$ bräuchten wir also $\eta'=0.04$. Dies ist das Ziel, also machen wir uns jetzt an die Arbeit und schätzen $S$ nach unten ab.\\

Wir müssen einerseits bestimmen für welche $n$ wir $h_1(n)=0$ haben und andererseits wie $$h_2(n):=\sum_{d|n} \psi_d$$ für diese $n$ ausfällt, um schließlich die Summe über diese $n$, also $S$, abschätzen zu können. Zuerst bemerken wir, dass $h_1(n)$ und $h_2(n)$ nur vom quadratfreien Teil von $n$ (also $\prod\limits_{p|n} p$) abhängen.\\

\noindent Hat $n$ keinen oder genau einen Primfaktor, so ist $h_1(n) \neq 0$. Wir betrachten im folgenden nur diejenigen $n$, die genau drei Primfaktoren haben, also $n=p_1^{\alpha_1} p_2^{\alpha_2} p_3^{\alpha_3}$. Die Analyse für $4,5,\ldots$ Primfaktoren geht analog und wird "`überproportional komplizierter"'. Wir haben keinen Weg gefunden alle $n$ gleichzeitig zu betrachten und dies erscheint auch sehr schwer, da $h_1(n)$ und $h_2(n)$ sehr stark von der Anzahl der Primfaktoren von $n$ und deren Lage zueinander abhängen. \\
Weiterhin liefert eine genauere Betrachtung der Summe (\ref{vb6-summe}) die Vermutung, dass die nicht quadratfreien $n$ nur einen vernachlässigbaren Beitrag produzieren\footnote{Diese Vermutung basiert auf der folgenden Ungleichung, die für $p_1,p_2,p_3 \geq M(\varepsilon)\in \N$ gilt: $$ \sum_{\alpha_1 \geq 1}\sum_{\alpha_2 \geq 1}\sum_{\alpha_3 \geq 1}\frac{1}{p_1^{\alpha_1}}\frac{1}{p_2^{\alpha_2}}\frac{1}{p_3^{\alpha_3}}-\frac{1}{p_1 p_2 p_3}= \frac{1}{(p_1-1)(p_2-1)(p_3-1)}-\frac{1}{p_1 p_2 p_3} \leq \frac{\varepsilon}{p_1 p_2 p_3}.$$}, also beschränken wir uns auf quadratfreie $n$. Wir betrachten ab jetzt also nur noch $n=p_1 p_2 p_3 \; (p_1<p_2<p_3 \text{ Primzahlen})$. \\

\noindent Für diese $n$ gilt
$$h_1(n)=0 \; \Longleftrightarrow \; p_1 p_2 p_3 \leq W \text{ oder } (p_1 p_2 \leq W \text{ und } p_3>W) . $$
Der Beitrag aller $n$ mit $1<n\leq W \,(\leq U)$ zu $S$ ist wegen $h_2(n)=\sum\limits_{d|n} \mu(d)=0$ auch gleich Null. Also bleibt jetzt die Aufgabe
\begin{equation}\label{vb6-summe-spezialfall} \sum_{n=p_1 p_2 p_3, \atop p_1 p_2 \leq W, \; p_3>W} h_2(n)^2 n^{-1} (e^{-n/X}-e^{-n\mathscr{L}^2/U}) \end{equation}
zu berechnen. Es gilt $2 w<u$ und $w+u<v$. Damit bekommt man für $p_1<p_2<p_3, \,p_1 p_2 \leq W$ und $p_3>W$ die Gleichung
\begin{equation}\label{vb6-h2} \left(\log \frac{V}{U}\right) h_2(p_1 p_2 p_3)=
\left\lbrace \begin{array}{ll}
\log \frac{p_1 p_2 p_3}{U} & p_3 \in (\frac{U}{p_1 p_2},\frac{U}{p_2}], \\
\log p_1 & p_3 \in (\frac{U}{p_2},\frac{U}{p_1}], \\
\log \frac{U}{p_3} & p_3 \in (\frac{U}{p_1},U], \\
\log \frac{V}{p_1 p_2 p_3}& p_3 \in (\frac{V}{p_1 p_2},\frac{V}{p_2}], \\
-\log p_1 & p_3 \in (\frac{V}{p_2},\frac{V}{p_1}], \\
\log \frac{p_3}{V} & p_3 \in (\frac{V}{p_1},V], \\
0& \text{sonst.} \\
\end{array} \right.
\end{equation}

\noindent Wir zeigen beispielhaft, welchen Beitrag man zur Summe (\ref{vb6-summe-spezialfall}) bekommt für $p_3 \in (\frac{U}{p_1},U]$ (dann ist $  U \leq n \leq V \leq X \mathscr{L}^{-1}$):
$$\sum_{n=p_1 p_2 p_3, \, p_1 p_2 \leq W, \atop p_3 \in [\frac{U}{p_1},U]} h_2(n)^2 n^{-1} (e^{-n/X}-e^{-n\mathscr{L}^2/U})
$$
\begin{eqnarray*}&=&\frac{(1+o(1))}{(v-u)^2 \mathscr{L}^2}\sum_{n=p_1 p_2 p_3,\atop p_1 p_2 \leq W, \, p_3 \in [\frac{U}{p_1},U]} (u^2\mathscr{L}^2 -2u\mathscr{L}\log p_3+ \log^2 p_3)\frac{1}{p_1 p_2 p_3} \\
&=:& (S_1+S_2+S_3) .\end{eqnarray*}
Dabei entspricht $S_i$ der Summe mit dem i-ten Term in der inneren Klammer. Mit $$\sum_{p \leq x} p^{-1} = \log \log x+C + O(\frac{1}{\log x})$$ bekommen wir für $q\rightarrow \infty$, dass
\begin{eqnarray*}
  S_1 &=& \frac{(1+o(1))u^2\mathscr{L}^2}{(v-u)^2 \mathscr{L}^2} \sum_{p_1 < W^{\frac{1}{2}}} \frac{1}{p_1} \sum_{p_1<p_2 \leq \frac{W}{p_1}} \frac{1}{p_2} \sum_{\frac{U}{p_1}< p_3 \leq U} \frac{1}{p_3} \\
  &=& \frac{(1+o(1))u^2}{(v-u)^2} \sum_{p_1 < W^{\frac{1}{2}}} \frac{1}{p_1} \left(\log \frac{\log (W/p_1)}{\log p_1} \right) \left( \log \frac{\log U}{\log (U/p_1)} \right)\\
  & = & (1+o(1))\frac{u^2}{(v-u)^2} \int_0^{\frac{w}{2}} \frac{1}{s} \left(\log \frac{w-s}{s}\right) \left(\log \frac{u}{u-s}\right) \,ds  .
\end{eqnarray*}
Beim letzten Gleichheitszeichen benutzt man partielle Summation und den Primzahlsatz und führt im anschließenden Integral die Substitution $t=q^s$ durch (wichtig: wähle als zu differenzierende Funktion $h(t)=t^{-1} \log \frac{\log (W/t)}{\log t} \log \frac{\log U}{\log (U/t)}$ aus).

Die gleiche Prozedur müssen wir auch für $S_2$ und $S_3$ machen und erhalten\footnote{Die Rechnungen lassen sich vereinfachen, denn es gilt $$S_1+S_2+S_3=\frac{1+o(1)}{(v-u)^2}\int_0^{w/2} s_1^{-1} \int_{s_1}^{w-s_1} s_2^{-1} \int_{u-s_1}^{u}\frac{u^2-2us_3+s_3^2}{s_3} \,ds_3\,ds_2\,ds_1 .$$} $S_1 + S_2 + S_3 \approx 0.0002$. Dieselbe Prozedur muss für die anderen Fälle in (\ref{vb6-h2}) gemacht werden. Insgesamt bekommen wir für den Beitrag der $n=p_1 p_2 p_3$ in (\ref{vb6-summe}) den Wert $\approx 0.002$, bekommen also ungefähr $\eta' \geq 0.002$. Der entsprechende Beitrag für die angesprochene Konstante $11.13\ldots$ ist dann $\eta \approx 0.01 \, (=\frac{1}{2w} 0.002)$ und damit viel zu klein. \\

Wir haben auch mit 4 und 5 Primfaktoren die Analyse vorangetrieben (man beachte, dass die Kompliziertheit stark ansteigt, einerseits bei der Analyse von $h_1(t)$ und $h_2(t)$, andererseits bei der Analyse der auftretenden Summen). Bei vier Primfaktoren haben wir beispielsweise für eine zu $S_1+S_2+S_3$ analoge Summe den Wert $\approx 0.00004$ berechnet. Dazu muss man sagen, dass im Fall von vier Primfaktoren es
mehr von diesen Summen geben wird, die zusammenzuaddieren sind. Wir haben auch ein paar Berechnungen im Fall mit fünf Primfaktoren durchgeführt. Zusammenfassend ist aber festzuhalten, dass wir durch unsere Berechnungen kein hinreichend großes $\eta'$ beweisen konnten, noch bekamen wir einen Hinweis darauf, dass durch eine viel längere Analyse ein solches zu erreichen gewesen wäre. \\
Ähnliches gilt für die Verwendung von Verbesserungspotential Nr. 6 in Verbindung mit dem Verfahren aus \cite[§4]{Hea92}. Dieses Verfahren haben wir kurz in §\ref{para-methoden-fuer-nullstellenfreie-regionen} unter der Bezeichnung "`erste spezielle Methode"' kommentiert, aber in unserer Arbeit nicht verwendet.

\section{Verbesserungspotential Nr. 8}
Dieses Potential lieferte gute Verbesserungen für den Fall $\lambda_1 \leq 0.60$, für größeres $\lambda_1$ gab es sehr geringe bis gar keine Verbesserungen. Wir hatten aber nur für letzteren Bereich Verbesserungen nötig, wollten wir Theorem \ref{haupttheorem} verbessern. Deswegen benutzen wir dieses Potential nicht. Man beachte, dass, wenn wir dieses Potential angewendet hätten, wir an einigen anderen Stellen der Diplomarbeit etwas weniger genau hätten arbeiten müssen, um $L=5.2$ zu bekommen. Jedoch erschien uns das Auslassen dieses Potentials und das etwas genauere Arbeiten an anderer Stelle als insgesamt günstiger.\\

Mit diesem Potential kann man \cite[§12]{Hea92} erweitern und die Abschätzungen für die Anzahl $N(\lambda)$ der Nullstellen im Bereich $\sigma \geq 1 - \mathscr{L}^{-1}\lambda , \; |t|\leq 1$ verbessern (vergleiche Tabellen 12 und 13). Die Güte der Verbesserung hängt dabei davon ab, wie gute Abschätzungen man von $\lambda_2$ in Abhängigkeit von $\lambda_1$ hat.\\
Weiterhin muss man bei diesem Potential - wie auch bei Verbesserungspotential Nr. 4 - keine weitere Arbeit mehr vornehmen. Aus \cite[§12, (16.7) und (16.8)]{Hea92} folgt sofort

\begin{lemma}[S.336 in \cite{Hea92}]
Seien $\varepsilon>0$, $\lambda \in (0,2)$, $\lambda_{21}>0$ Konstanten mit $\lambda_{21}<\lambda$. Weiterhin sei $f$ eine Funktion, die Bedingung 1 und 2 erfüllt und für die gilt
$$F(\lambda-\lambda_{21})>\frac{f(0)}{6} \;\; \text{ und } \;\; \left(F(\lambda-\lambda_{21})-\frac{f(0)}{6}\right)^2 > F(-\lambda_{21})\frac{f(0)}{6} .$$
Ferner gelte $\lambda_2 \geq \lambda_{21}$. Dann haben wir für $q \geq q_0(\varepsilon,f)$ die Abschätzung
$$N(\lambda) \leq 3 S_1 + 2, $$
wobei
$$S_1 = \frac{F(-\lambda_{21})(F(-\lambda_{21})-\frac{f(0)}{6})}{(F(\lambda-\lambda_{21})-\frac{f(0)}{6})^2 - F(-\lambda_{21})\frac{f(0)}{6}} +\varepsilon .  $$
\end{lemma}

\end{document}